\documentclass[11pt,a4paper,reqno]{amsart}

\usepackage[T1]{fontenc}
\usepackage[utf8]{inputenc}
\usepackage[
  a4paper,
  textwidth=15.5cm,
  textheight=23.4cm,
  centering,
  headsep=7mm,
  footskip=11mm
]{geometry}
\usepackage{microtype}
\usepackage{amsmath}
\usepackage{amssymb}
\let\leq\leqslant
\let\geq\geqslant

\usepackage{amsthm}
\usepackage[scr=boondox,cal=esstix]{mathalpha}
\usepackage{graphicx,tikz,enumitem,xurl,placeins,needspace}
\usetikzlibrary{decorations.pathreplacing}
\usepackage{booktabs}
\usepackage{array}
\renewcommand{\arraystretch}{1.1}
\newcommand{\paddedmathcell}[1]{\vcenter{\kern3pt\hbox{$\displaystyle#1$}\kern3pt}}
\usepackage{xcolor}
\usepackage[
  font=small,
  labelfont=sc,
  labelsep=period,
  justification=centering
]{caption}
\usepackage{cite}
\usepackage[hidelinks]{hyperref}
\usepackage[nameinlink,noabbrev]{cleveref}

\allowdisplaybreaks[2]
\newtheorem{thm}{Theorem}[section]

\newtheorem{lem}[thm]{Lemma}
\newtheorem{prop}[thm]{Proposition}

\theoremstyle{definition}
\newtheorem{defn}[thm]{Definition}
\theoremstyle{remark}
\newtheorem{rem}[thm]{Remark}

\numberwithin{equation}{section}
\newcommand{\B}{\mathcal{B}}

\makeatletter
\newcommand{\BibAlias}[2]{%
  \@ifundefined{b@#2}{}{\expandafter\let\csname b@#1\expandafter\endcsname\csname b@#2\endcsname}}
\makeatother
\AtBeginDocument{%
\BibAlias{Ale23}{Alexandrov2023}%
\BibAlias{And70(1)}{Andreev1970Compact}%
\BibAlias{And70(2)}{Andreev1970Finite}%
\BibAlias{Bou68}{Bourbaki1968}%
\BibAlias{bugaenko-compact-2}{Bugaenko1992}%
\BibAlias{Coxeter:1934}{Coxeter1934}%
\BibAlias{Ess94}{Esselmann1994}%
\BibAlias{Ess96}{EsselmannCompact1996}%
\BibAlias{Kap74}{Kaplinskaya1974}%
\BibAlias{Lanner:1950}{Lanner1950}%
\BibAlias{P1882}{Poincare1882}%
\BibAlias{prokhorov-1986}{Prokhorov1987}%
\BibAlias{Tum04b}{Tumarkin2004Plus3}%
\BibAlias{Tum07}{Tumarkin2007}%
\BibAlias{vinberg-1985}{Vinberg1985}%
\BibAlias{tumarkin-nplustwo}{Tumarkin2004}%
}

\usepackage{longtable,algorithm,algpseudocode}

\newcommand{\wt}{\widetilde}

\DeclareMathOperator{\dist}{dist}
\crefname{algorithm}{Algorithm}{Algorithms}
\crefname{lem}{Lemma}{Lemmas}
\Crefname{lem}{Lemma}{Lemmas}
\crefname{prop}{Proposition}{Propositions}
\Crefname{prop}{Proposition}{Propositions}
\crefname{thm}{Theorem}{Theorems}
\Crefname{thm}{Theorem}{Theorems}

\usepackage{etoolbox}
\AtBeginEnvironment{tabular}{\small\linespread{1}\selectfont}
\AtBeginEnvironment{longtable}{\small\linespread{1}\selectfont}
\crefname{part}{Part}{Parts}
\Crefname{part}{Part}{Parts}
\begin{document}
\title[Dimension bounds for hyperbolic Coxeter polytopes]
{Dimension bounds for finite-volume hyperbolic Coxeter polytopes}
\author{Lizi Guo}
\address{Mathematics Department, University of Wisconsin--Madison,
480 Lincoln Dr, Madison, WI 53706, USA}
\email{llguo@wisc.edu}
\author{Fangting Zheng}
\address{Department of Mathematical Sciences, Xi'an Jiaotong-Liverpool University,
Suzhou 215123, China}
\email{Fangting.Zheng@xjtlu.edu.cn}
\date{October 31, 2026}
\thanks{Fangting Zheng was supported by the NSFC (No.~12471067).}
\begin{abstract}

	We prove that a finite-volume hyperbolic Coxeter polytope has dimension
	at most $81$, substantially improving Prokhorov's bound of $995$.
	A refinement of Prokhorov's method gives the intermediate bound $162$,
	while an extension of Vinberg's planar-angle comparison to the
	finite-volume setting yields the sharper bound.
	Both arguments combine improved weighting functions with estimates
	using cusp geometry to control weight contributions of angles at ideal
	vertices or edges incident to them.

\end{abstract}
\maketitle
\tableofcontents
\medskip
\section{Introduction}\label{sec:introduction}

Let $\mathbb{H}^n$ be hyperbolic $n$-space. Among the discrete subgroups of $\mathrm{Isom}(\mathbb{H}^n)$, the finite-covolume reflection groups form a class that is geometrically explicit and surprisingly rich. The fundamental domain of such a group $\Gamma$ is a Coxeter polytope $P \subset \mathbb{H}^n$, i.e., a convex polytope whose dihedral angles are submultiples of $\pi$. Conversely, any Coxeter polytope $P$ gives rise to a discrete group $\Gamma(P)$, generated by the reflections in its bounding hyperplanes, for which $P$ is a fundamental domain.
We may thus pass freely between a Coxeter polytope and the group it defines. This correspondence converts questions about the reflection groups --- existence, classification, arithmeticity, commensurability --- into combinatorial and metric questions about the Coxeter polytopes.  Let $s_i$ be the
reflection in the supporting hyperplane $F_i$ of a Coxeter polytope $P$. The group $\Gamma(P)$ admits a presentation
\[
\left\langle s_1,\ldots,s_m\ \middle|\ s_i^2=1,\quad (s_i s_j)^{k_{ij}}=1\ \ (i<j,\ k_{ij}<\infty)\right\rangle,
\]
where $k_{ij}\in\{2,3,\ldots\}$ when $F_i$ and $F_j$ meet along a
ridge at a dihedral angle $\pi/k_{ij}$. Otherwise $k_{ij}=\infty$ and no relation on $s_i s_j$ is imposed.

In dimension $2$, the theory of hyperbolic polygons and their reflection groups goes back to the work of Poincar\'{e} \cite{Poincare1882} on discrete isometry groups of the hyperbolic plane. The complete characterization of $3$-dimensional acute-angled hyperbolic polytopes was obtained much later by Andreev \cite{Andreev1970Compact,Andreev1970Finite}; see also \cite{RoederHubbardDunbar2007}. This characterization became a cornerstone of Thurston's geometrization of $3$-dimensional Haken manifolds \cite{Thurston1979}. More generally, finite-covolume reflection groups give rise to finite-volume hyperbolic manifolds after passing to torsion-free subgroups of finite index. In higher dimensions, reflection groups provide concrete constructions of hyperbolic manifolds alongside arithmetic and hybrid constructions.

The spherical and Euclidean Coxeter groups were completely classified nearly a hundred years ago
\cite{Coxeter1934,Bourbaki1968}. In the hyperbolic setting, however, even the maximal dimension admitting a finite-volume
Coxeter polytope remains unknown.
Vinberg proved that compact hyperbolic Coxeter polytopes, as well
as finite-volume hyperbolic Coxeter polytopes whose reflection
groups are arithmetic, have dimension at most $29$
\cite{Vinberg1981,Vinberg1984}.
For arbitrary finite-volume Coxeter polytopes, Prokhorov gave the bound \(n\leq995\)
\cite{Prokhorov1987}, using a combinatorial inequality for
edge-simple polytopes due to Khovanskii \cite{Khovanskii1986}, which
generalizes an earlier estimate of Nikulin for simple polytopes
\cite{Nikulin1982}.
Nikulin announced the stronger bound $n<300$ in his ICM address,
where he only outlined an argument for this estimate
\cite[p.~659]{Nikulin1986}.
The highest-dimensional known compact and noncompact finite-volume
examples have dimensions \(8\) and \(21\), respectively, and are due
to Bugaenko \cite{Bugaenko1992} and Borcherds \cite{Borcherds1987}. The gap between the known examples and the general upper bound remains large. Our main result improves the latter for the entire finite-volume class.

\begin{thm}\label{thm:dimension-bound}
Let $P\subset\mathbb H^n$ be a finite-volume hyperbolic Coxeter
polytope. Then $n\leq81$.
\end{thm}

We develop following two weighting approaches for this dimension problem. Improving Prokhorov's method gives $n\leq162$, as stated in Theorem~\ref{edge:thm:dimension-bound}); extending Vinberg's
planar-angle method gives $n\leq81$. We retain both arguments for
their distinct geometric observations and interesting intermediate estimates involved.
The first compares weights summed over each 3-face with the same
weights summed around each edge. Square cusp sections reduce the
edgewise upper bound and constrain the average number of facets
of a 3-face. The second compares weights summed over each 2-face
with their sums at vertices. At an ideal vertex, the product structure
of the link determines the number of contributing angles, while
Timorin's truncation identities supply a term to subtract from their
total weight. The remaining contribution can be bounded per incident
edge and used in the same double-counting comparison.

Prokhorov's proof of the bound \(n\leq995\) adapts Vinberg's
distance-weighting method to the finite-volume case. Compact Coxeter polytopes
are simple, and their elliptic vertex diagrams control the sums of
planar-angle weights. As Prokhorov explains \cite[p.~403]{Prokhorov1987},
this estimate does not extend directly to ideal vertices, where many
planar angles can contribute to a single vertex sum. He therefore
passes to dihedral angles along the edges of 3-faces. Finite-volume
Coxeter polytopes need not be simple, but they are edge-simple, and
their edge diagrams are elliptic of order $n-1$.
The weights are chosen to be sparse around each edge but sufficient
on each 3-face with few facets. Elliptic diagrams control the edgewise
sums, while Lann\'er and quasi-Lann\'er diagrams supply the short paths
needed for the cellwise lower bounds.
After the exceptional right-angled triangular bipyramids are controlled
through four-dimensional incidence estimates, double counting gives an
average weight \(a_w\) bounded below by \(1/4-O(n^{-1})\) and above by
\(O(n^{-1})\). These two bounds are incompatible in sufficiently high
dimension, which yields the bound.

Our improvement combines a refinement of these weights with
cusp geometry. The refined weighting rule reduces the leading edgewise coefficient
from \(29/3\) to \(73/12\). The more substantial gain comes from the
exceptional 3-faces, whose square links at ideal vertices
force a reduction in the edgewise upper bound. Retaining cusp terms in
Timorin's truncation identities gives another restriction on the average
number of facets of a 3-face. For \(n\geq163\), the resulting
cellwise lower bound for the average weight exceeds its edgewise
upper bound for every proportion of exceptional 3-faces, including
zero. Hence \(n\leq162\).

For the planar-angle method, we assign a main weight and an auxiliary
weight to each vertex of a 2-face. The main weights sum to at least
two on a triangle and one on a quadrilateral; the auxiliary weights
strengthen the triangle estimate. Two classes of long angles receive
main weight zero: those in affine cycles and those whose nonspecial
paths extend sufficiently far at both ends. Gram-matrix arguments
for the full quadrilateral diagram show that these reductions preserve
the required facewise lower bounds.

At an ideal vertex $u$, an angle of a 2-face corresponds to an edge
of the link $L_u=\prod_i\Delta_i$. A fixed edge of $\Delta_i$
gives a distinct angle for every choice of a vertex in each other
factor, hence contributes $\prod_{j\ne i}(d_j+1)$ times, where
$d_j=\dim\Delta_j$. This counts distinct angles represented by the
same marked pair in the factor diagram. Even after division by the
number of incident edges, their total weight need not be uniformly
bounded as the number of factors grows.

The facewise weight inequality and the inequality obtained from
Timorin's truncation identities involve the same number of 2-faces.
Eliminating that number gives a lower bound for the total weight
after an explicit cusp term has been subtracted at each ideal vertex.
We estimate this same difference locally, dividing by the number
of incident edges. Product-link identities bound it in terms of the
simplex-factor dimensions; concavity gives a uniform estimate and
reduces the remaining dimension-dependent estimates to finite checks.
Together with the elliptic estimates at finite vertices and the
edge-endpoint incidence identity, these bounds contradict the
facewise lower bound for every $n\geq82$.

Stronger bounds are available when
arithmeticity, dihedral angles, or vertex structure are prescribed. In the
noncompact arithmetic case, Esselmann proved that $n\leq21$ and
$n\neq20$ \cite{EsselmannArithmetic1996}. Borcherds' example shows
that the upper bound is sharp.
For right-angled polytopes, the compact bound $n\leq4$ is sharp
\cite{PotyagailoVinberg2005}, and the finite-volume bound is
$n\leq12$ \cite{Dufour2010}, improving the earlier bound $n\leq14$
of Potyagailo and Vinberg \cite{PotyagailoVinberg2005}.
The ideal right-angled subclass satisfies $n\leq6$
\cite{Kolpakov2012}, while simple ideal Coxeter polytopes satisfy
$n\leq8$ \cite{FeliksonTumarkin2008Ideal}.
Furthermore, Jacquemet and Tschantz proved that finite-volume hyperbolic
Coxeter cubes have dimension at most $5$
\cite{JacquemetTschantz2018}.
For polytopes with pairwise intersecting facets, Felikson and
Tumarkin proved that $n\leq4$ in the compact case and $n\leq9$
in the simple noncompact finite-volume case
\cite{FeliksonTumarkin2008Intersect}.
They also showed that a compact hyperbolic Coxeter polytope
with exactly one pair of disjoint facets has at most $n+3$
facets and satisfies $n\leq6$ or $n=8$
\cite{FeliksonTumarkin2009OnePair}.
For compact $3$-free polytopes, in which every inclusion-minimal
collection of facets with empty intersection has two members,
Burcroff obtained
$n\leq13$ \cite{Burcroff2024}. Alexandrov subsequently improved
this bound to $n\leq12$ \cite{Alexandrov2023}.

Dimension bounds are also known when the difference between the number of facets $m$ and the dimension $n$ of the polytope is small. Let $D_{\mathrm c}(k)$ and $D_{\mathrm{nc}}(k)$
denote the largest dimensions of compact and noncompact finite-volume
hyperbolic Coxeter $n$-polytopes with exactly $n+k$ facets, respectively.
The known sharp compact values are
$(D_{\mathrm c}(1),D_{\mathrm c}(2),D_{\mathrm c}(3),D_{\mathrm c}(4))
=(4,5,8,7)$, by the work of Lann\'er, Kaplinskaya, Esselmann,
Tumarkin, and Felikson--Tumarkin
\cite{Lanner1950,Kaplinskaya1974,EsselmannCompact1996,Esselmann1994,Tumarkin2007,FeliksonTumarkin2008Few}.
In the noncompact case, the simplex classification and Tumarkin's
results give $(D_{\mathrm{nc}}(1),D_{\mathrm{nc}}(2),D_{\mathrm{nc}}(3))
=(9,17,16)$ \cite{Chein1969,Tumarkin2004,Tumarkin2004Plus3}.
For facet excess $k=5,6,7,8,9,10$, Burcroff obtained the respective
upper bounds $9,12,15,18,22,26$ for $D_{\mathrm c}(k)$
\cite{Burcroff2024}.\footnote{Burcroff also reports the stronger
bound $D_{\mathrm c}(5)\leq8$ from unpublished work of Felikson and
Tumarkin, citing a personal communication from Tumarkin in 2021
\cite{Burcroff2024}.}

The paper is organized as follows. Section~\ref{section:cchp} contains
the common preliminaries. In Part~\ref{part:edge},
Section~\ref{edge:sec:comparison} recalls Prokhorov's
cellwise--edgewise comparison, and Section~\ref{ext:sec:weights}
refines the weights. Sections~\ref{ext:sec:cusp} and~\ref{ext:sec:h}
use cusp geometry to bound, respectively, the weight around an edge
and the average number of facets of a 3-face.
Section~\ref{ext:sec:conclusion} combines these bounds to obtain
$n\leq162$. In Part~\ref{part:vertex},
Section~\ref{two:sec:comparison} sets up the facewise--vertexwise
comparison. Sections~\ref{two:sec:weights} and~\ref{two:sec:deletions}
define the planar-angle weights and prove their facewise lower bounds.
Section~\ref{two:sec:budgets} counts angle multiplicities and bounds
the vertex contributions after subtracting the retained cusp terms.
Section~\ref{two:sec:conclusion} proves $n\leq81$.
The appendices give the classification diagrams and the finite checks.

\textbf{Acknowledgements}.
The second author would like to express her deep gratitude to Jiming Ma for introducing her to the problem of dimension bounds about ten years ago, and for many stimulating discussions on hyperbolic Coxeter polytopes over the years of supervision and collaboration.

\textbf{AI disclosure}.
The overall conception and strategy of this work are due to the
authors. Substantial AI assistance was essential to its completion.
ChatGPT (OpenAI) assisted in developing mathematical arguments and
writing the computational code. It helped formulate the inequalities
relating equal-length edge multiplicities to cusp weight deficits in
Lemma~\ref{ext:lem:factor}.
All extensive finite computations in this work were carried out with
AI assistance. In particular, computations on reduced
Lann\'er and quasi-Lann\'er diagrams and related marked configurations
informed the design and testing of the weighting rules in
\eqref{ext:eq:weight} and
\eqref{two:eq:weight-table}--\eqref{two:eq:final-weight}.
ChatGPT also developed the proof of Lemma~\ref{two:lem:tree-delete},
which justifies assigning zero weight to certain long tree angles
while preserving the required lower bound for the total weight
of a quadrilateral.
The computational code is available at
\url{https://github.com/GeoTopChristy/CoxeterDimBound}.
Claude (Anthropic) assisted with drafting the manuscript.
The authors take full responsibility for the accuracy and integrity
of the paper.

\section{Preliminaries}\label{section:cchp}
We begin with the combinatorics of convex polytopes,
then list some essential facts about hyperbolic Coxeter polytopes.

\subsection{Combinatorics of convex polytopes}\label{subsec:face-numbers}
We recall some combinatorial definitions and results used below.
We write \emph{$n$-polytope} for an $n$-dimensional polytope and
\emph{$i$-face} for an $i$-dimensional face. A \emph{facet} and a
\emph{ridge} are faces of codimension one and two, respectively.
An $n$-polytope is \emph{simple} if every vertex is incident with
exactly $n$ edges, and \emph{edge-simple} if every edge belongs to
exactly $n-1$ facets.
The \emph{vertex link} at a vertex $u$ of a convex $n$-polytope $D$
is the $(n-1)$-polytope, considered up to combinatorial equivalence,
obtained by intersecting $D$ with a hyperplane that separates $u$
from all its other vertices. Its combinatorial type is independent
of the choice of this hyperplane. Its $i$-faces correspond to the
$(i+1)$-faces of $D$ containing $u$.

For a convex $n$-polytope $D$, let $f_i(D)$ denote the number of its
$i$-faces, including $D$ itself as $f_n(D)=1$.
The vector $(f_0(D),\ldots,f_n(D))$ is its $f$-vector.
Following Timorin \cite[Introduction]{Timorin2001}, 
the $h$-numbers are defined as
\begin{equation}\label{ext:eq:h}
h_j(D):=\sum_{i=j}^{n}(-1)^{i-j}\binom ij f_i(D),
\qquad 0\leq j\leq n.
\end{equation}
The inverse transformation is
\begin{equation}\label{ext:eq:h-inverse}
f_k(D)=\sum_{j=k}^{n}\binom jk h_j(D),
\qquad 0\leq k\leq n.
\end{equation}
Thus the $f$-vector and the $h$-vector $(h_0(D),\ldots,h_n(D))$
determine one another. We take $\binom jk=0$ when $j<k$.
The associated \emph{$h$-polynomial}, with indeterminate $t$, is
\begin{equation}\label{ext:eq:h-polynomial}
h_D(t):=\sum_{j=0}^{n}h_j(D)t^j
       =\sum_{i=0}^{n}f_i(D)(t-1)^i.
\end{equation}
The second equality follows by expanding $(t-1)^i$. Replacing $t$
by $t+1$ and expanding $(t+1)^j$ by the binomial theorem, we obtain
\eqref{ext:eq:h-inverse} by comparing the coefficients of $t^k$.
Euler's relation gives $h_0(D)=1$, and $h_n(D)=f_n(D)=1$.
Setting $t=1$ gives
\begin{equation}\label{ext:eq:h-at-one}
h_D(1)=\sum_{j=0}^n h_j(D)=f_0(D).
\end{equation}

For a simple $n$-polytope $D$, the $h$-numbers satisfy the
Dehn--Sommerville relations $h_j(D)=h_{n-j}(D)$. Although this
symmetry can fail for an edge-simple polytope, truncating its nonsimple
vertices relates its $h$-numbers to those of a simple polytope.

Let $D$ be an edge-simple $n$-polytope. Truncate each nonsimple
vertex $u$ by a sufficiently close hyperplane, choosing the removed
neighborhoods to be disjoint. Denote the resulting polytope by $Q$
and the new facet at $u$ by $L_u$, which is a vertex link of $D$ at $u$.
Each vertex of $L_u$ corresponds to an edge of $D$ incident
to $u$ and belongs to exactly $n-1$ facets of $L_u$. Thus each
$L_u$ is simple, and $Q$ is simple as well.

Each positive-dimensional face of polytope $D$ persists after sufficiently
small truncations, and the additional positive-dimensional faces are
exactly those in the inserted facets.  Consequently
\[
\begin{aligned}
f_i(Q)&=f_i(D)+\sum_{u\text{ nonsimple}}f_i(L_u) &&(1\leq i\leq n),\\
f_0(Q)&=f_0(D)+\sum_{u\text{ nonsimple}}\bigl(f_0(L_u)-1\bigr).
\end{aligned}
\]
Here face numbers and $h$-numbers beyond the dimension of a polytope
are taken to be zero. Substitution in \eqref{ext:eq:h-polynomial} gives
\[
h_Q(t)=h_D(t)+\sum_{u\text{ nonsimple}}\bigl(h_{L_u}(t)-1\bigr),
\]
and hence
\begin{equation}\label{ext:eq:truncation}
h_j(D)=h_j(Q)-\sum_{u\text{ nonsimple}}h_j(L_u),\qquad 1\leq j\leq n.
\end{equation}
This is the truncation identity in \cite[Lemma~2.1]{Timorin2001}.
Timorin also proved the upper-half nonnegativity for edge-simple
polytopes \cite[Theorem~2.4]{Timorin2001}:
\begin{equation}\label{ext:eq:h-upper-nonnegative}
h_j(D)\geq0,\qquad n/2<j\leq n.
\end{equation}
We use these facts in
Section~\ref{ext:sec:h} and
Subsection~\ref{two:subsec:retained-cusp}.  Write $a_i=f_i(D)$ and define
the average number of $i$-faces of a $k$-face by
\begin{equation}\label{two:eq:average-faces}
 a_k^i:=\frac1{a_k}\sum_{\substack{F\subseteq D\\\dim F=k}}f_i(F),
 \qquad 0\leq i<k\leq n.
\end{equation}
Thus $a_2^0=a_2^1$ is the average numbers of vertices or
edges of a 2-face, while $a_3^2$ and $a_3^1$ are the average numbers
of facets and edges of a 3-face. 

Every positive-dimensional $i$-face belongs to exactly
$\binom{n-i}{k-i}$ $k$-faces. 
Choosing $n-k$ of these facets determines a containing $k$-face.
Counting incidences gives
\begin{equation}\label{two:eq:positive-face-incidences}
 a_ka_k^i=\binom{n-i}{k-i}a_i,
 \qquad 1\leq i<k\leq n.
\end{equation}
In particular,
\begin{equation}\label{two:eq:edge-face-incidence}
 a_2a_2^0=a_2a_2^1=(n-1)a_1,
\end{equation}
and
\begin{equation}\label{eq:incidences}
 a_3a_3^2=(n-2)a_2,\qquad
 a_3a_3^1=\binom{n-1}{2}a_1.
\end{equation}
For $i=0$, the number of containing faces need not be constant,
since an edge-simple polytope can have nonsimple vertices.
The equality for $a_2^0$ follows from the polygon identity
$a_2^0=a_2^1$, rather than from a uniform vertex count.

The Nikulin--Khovanskii inequality for edge-simple polytopes is
\begin{equation}\label{eq:nikulin-khovanskii}
 a_k^i<\binom{n-i}{n-k}
 \frac{\binom{\lfloor n/2\rfloor}{i}+\binom{\lceil n/2\rceil}{i}}
      {\binom{\lfloor n/2\rfloor}{k}+\binom{\lceil n/2\rceil}{k}},
 \qquad 0\leq i<k,\quad 2\leq k\leq\lfloor n/2\rfloor.
\end{equation}
Nikulin proved the estimate for simple polytopes, and Khovanskii
extended it to edge-simple polytopes
\cite{Nikulin1982,Khovanskii1986}. For $n\geq4$, the cases
$(k,i)=(2,0),(2,1)$ give
\begin{equation}\label{two:eq:average-two-face-bound}
 a_2^0=a_2^1<u_{n,2}:=\begin{cases}
 4+\dfrac4{n-2},&n\text{ even},\\[4pt]
 4+\dfrac4{n-1},&n\text{ odd}.
 \end{cases}
\end{equation}
For $n\geq6$, the cases $(k,i)=(3,2),(3,1)$ give
\begin{equation}\label{eq:A-E-bounds}
\begin{aligned}
 a_3^2&<
 \begin{cases}6+\dfrac{12}{n-4},&n\text{ even},\\[3pt]
               6+\dfrac{12}{n-3},&n\text{ odd},\end{cases}
 &a_3^1&<
 \begin{cases}12+\dfrac{36}{n-4},&n\text{ even},\\[3pt]
               12+\dfrac{36}{n-3},&n\text{ odd}.
 \end{cases}
\end{aligned}
\end{equation}
The upper bounds for $a_3^2$ and $a_2^0$ tend to six and four,
respectively. These constants motivate the later cellwise target
$7-f_2(F)$ for 3-faces and the facewise target $5-f_0(F)$ for 2-faces.
The bound for
$a_3^1$, which tends to twelve, which together with
the incidence identity \eqref{eq:incidences}, it converts a linear
bound on the weight sum at an edge into an upper bound of order
$n^{-1}$ for the average total weight of a 3-face; see
\eqref{eq:normalized-edge-bound}.
Finite-volume hyperbolic Coxeter
polytopes are edge-simple, see Theorem \ref{Vinberg:thm3.1}, therefore these estimates apply.

\subsection{Hyperbolic Coxeter polytopes}\label{subsec:hyperbolic-polytopes}
A hyperbolic polytope is \emph{acute-angled} if all its dihedral
angles are at most $\pi/2$. Right angles are included
\cite{Vinberg1985}. It is a \emph{Coxeter polytope} if all its
dihedral angles are of the form $\pi/m$ for integers $m\geq2$.
Every face of a Coxeter polytope that meets $\mathbb H^n$ is
acute-angled, although it need not
itself be a Coxeter polytope.
For an acute-angled hyperbolic $n$-polytope
$P=\bigcap_{i\in I}H_i^-\subset\mathbb H^n$, where $H_i^-$ are its
bounding half-spaces, choose outward unit normals $e_i$ in the Lorentz
space $\mathbb E^{n,1}$. Its \emph{Gram matrix} is
$G(P)=(\langle e_i,e_j\rangle)_{i,j\in I}$.
The entries record the relative positions of the supporting hyperplanes:
\[
\langle e_i,e_j\rangle=\begin{cases}
1,&\text{if }i=j,\\
-\cos\alpha_{ij},&\text{if }H_i,H_j\text{ intersect at dihedral angle }\alpha_{ij},\\
-\cosh\rho_{ij},&\text{if }H_i,H_j\text{ are ultraparallel at distance }\rho_{ij},\\
-1,&\text{if }H_i,H_j\text{ are parallel.}
\end{cases}
\]
Here parallel hyperplanes meet at infinity, whereas ultraparallel
hyperplanes are disjoint even at infinity.

For a finite-volume polytope $P$, let $\overline P$ be its closure in
$\overline{\mathbb H^n}=\mathbb H^n\cup\partial\mathbb H^n$.
Its vertices in $\mathbb H^n$ and on $\partial\mathbb H^n$ are called
\emph{finite} and \emph{ideal} vertices, respectively. Such a polytope
is compact if and only if all its vertices are finite.
We count the faces of $\overline P$, including its ideal vertices,
when writing $f_i(P)$. In the Klein model, $\overline P$ is a convex
polytope, so the definitions of $f_i$ and $h_j$ in
Subsection~\ref{subsec:face-numbers} apply.

For a Coxeter polytope, the \emph{Coxeter diagram}
$\Sigma=\Sigma(P)$ has one node $i$ for each supporting hyperplane
$H_i$. Two nodes $i,j$ are joined by a solid edge labeled $k_{ij}$
if the hyperplanes intersect at angle $\pi/k_{ij}$, by a bold edge
if they are parallel, and by a dotted edge, sometimes labeled
$\cosh\rho_{ij}$, if they are ultraparallel at distance $\rho_{ij}$.
We omit edges with $k_{ij}=2$ and omit the label when $k_{ij}=3$.
An edge with label three is called \emph{ordinary} or \emph{simple};
an edge with label at least four is called \emph{heavy}. A bold edge
is also treated as heavy in the path rules.
The \emph{order} $|\Sigma|$ is the number of nodes. The
\emph{signature} and \emph{rank} of $\Sigma$ are those of its Gram
matrix $G(\Sigma)$, recovered from the edge labels by the formula above.
Here the signature records the numbers of positive and negative
eigenvalues, with zero eigenvalues omitted. A subdiagram will always
mean the induced subdiagram on its node set. Every principal submatrix
of $G(P)$ has at most one negative eigenvalue.

For a diagram $D$, we also write $G_D=G(D)$ for its Gram matrix.
Whenever $G_D$ is nonsingular, $G_D^{-1}$ denotes its inverse,
and $(G_D^{-1})_{uv}$ denotes the entry indexed by $u,v$.
The cofactor formula gives
\[
(G_D^{-1})_{uu}=\frac{\det G_{D\setminus\{u\}}}{\det G_D}.
\]
In $G_{D\setminus\{x\}}^{-1}$, the row and column indexed by $x$
are deleted before the resulting matrix is inverted.

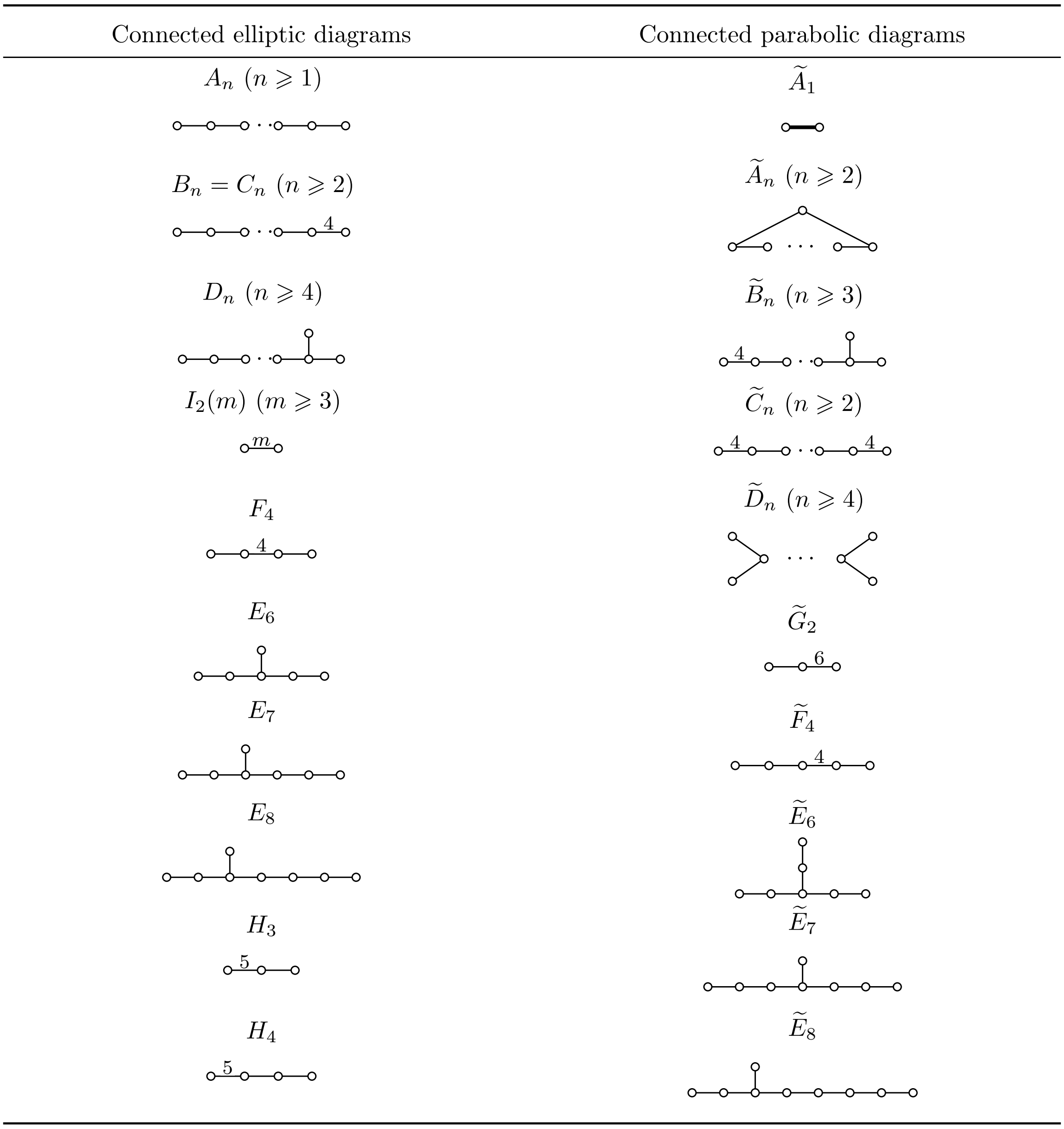
\begin{figure}[!htb]\centering\small
\begin{tabular}{@{}cc@{}}\toprule
Connected elliptic diagrams&Connected parabolic diagrams\\\midrule
\begin{minipage}[c]{.475\linewidth}\centering $A_n\ (n\geq1)$\par\begin{tikzpicture}[line width=.55pt,baseline=(current bounding box.center)]
\path[use as bounding box](-1.75,-0.45)rectangle(1.75,0.45);
\coordinate (n0) at (-1.2000,0.0000);
\coordinate (n1) at (-0.7200,0.0000);
\coordinate (n2) at (-0.2400,0.0000);
\coordinate (n3) at (0.2400,0.0000);
\coordinate (n4) at (0.7200,0.0000);
\coordinate (n5) at (1.2000,0.0000);
\draw[](n0)--(n1);
\draw[](n1)--(n2);
\node at (0.0000,0.0000){$\cdots$};
\draw[](n3)--(n4);
\draw[](n4)--(n5);
\draw[fill=white](n0)circle(1.6pt);
\draw[fill=white](n1)circle(1.6pt);
\draw[fill=white](n2)circle(1.6pt);
\draw[fill=white](n3)circle(1.6pt);
\draw[fill=white](n4)circle(1.6pt);
\draw[fill=white](n5)circle(1.6pt);
\end{tikzpicture}\end{minipage}&\begin{minipage}[c]{.475\linewidth}\centering $\widetilde A_1$\par\begin{tikzpicture}[line width=.55pt,baseline=(current bounding box.center)]
\path[use as bounding box](-1.75,-0.45)rectangle(1.75,0.45);
\coordinate (n0) at (-0.2400,0.0000);
\coordinate (n1) at (0.2400,0.0000);
\draw[line width=1.5pt](n0)--(n1);
\draw[fill=white](n0)circle(1.6pt);
\draw[fill=white](n1)circle(1.6pt);
\end{tikzpicture}\end{minipage}\\[.5mm]
\begin{minipage}[c]{.475\linewidth}\centering $B_n=C_n\ (n\geq2)$\par\begin{tikzpicture}[line width=.55pt,baseline=(current bounding box.center)]
\path[use as bounding box](-1.75,-0.45)rectangle(1.75,0.45);
\coordinate (n0) at (-1.2000,0.0000);
\coordinate (n1) at (-0.7200,0.0000);
\coordinate (n2) at (-0.2400,0.0000);
\coordinate (n3) at (0.2400,0.0000);
\coordinate (n4) at (0.7200,0.0000);
\coordinate (n5) at (1.2000,0.0000);
\draw[](n0)--(n1);
\draw[](n1)--(n2);
\node at (0.0000,0.0000){$\cdots$};
\draw[](n3)--(n4);
\draw[](n4)-- node[midway,above,fill=white,inner sep=.7pt,font=\scriptsize]{$4$}(n5);
\draw[fill=white](n0)circle(1.6pt);
\draw[fill=white](n1)circle(1.6pt);
\draw[fill=white](n2)circle(1.6pt);
\draw[fill=white](n3)circle(1.6pt);
\draw[fill=white](n4)circle(1.6pt);
\draw[fill=white](n5)circle(1.6pt);
\end{tikzpicture}\end{minipage}&\begin{minipage}[c]{.475\linewidth}\centering $\widetilde A_n\ (n\geq2)$\par\begin{tikzpicture}[line width=.55pt,baseline=(current bounding box.center)]
\path[use as bounding box](-1.75,-0.62)rectangle(1.75,0.62);
\coordinate (n0) at (-1.0000,-0.1700);
\coordinate (n1) at (-0.5000,-0.1700);
\coordinate (n2) at (0.5000,-0.1700);
\coordinate (n3) at (1.0000,-0.1700);
\coordinate (n4) at (0.0000,0.3500);
\draw[](n0)--(n1);
\node at (0.0000,-0.1700){$\cdots$};
\draw[](n2)--(n3);
\draw[](n3)--(n4);
\draw[](n4)--(n0);
\draw[fill=white](n0)circle(1.6pt);
\draw[fill=white](n1)circle(1.6pt);
\draw[fill=white](n2)circle(1.6pt);
\draw[fill=white](n3)circle(1.6pt);
\draw[fill=white](n4)circle(1.6pt);
\end{tikzpicture}\end{minipage}\\[.5mm]
\begin{minipage}[c]{.475\linewidth}\centering $D_n\ (n\geq4)$\par\begin{tikzpicture}[line width=.55pt,baseline=(current bounding box.center)]
\path[use as bounding box](-1.75,-0.54)rectangle(1.75,0.54);
\coordinate (n0) at (0.6750,-0.1850);
\coordinate (n1) at (0.2250,-0.1850);
\coordinate (n2) at (-0.2250,-0.1850);
\coordinate (n3) at (-0.6750,-0.1850);
\coordinate (n4) at (-1.1250,-0.1850);
\coordinate (n5) at (1.1250,-0.1850);
\coordinate (n6) at (0.6750,0.1850);
\draw[](n0)--(n1);
\node at (0.0000,-0.1850){$\cdots$};
\draw[](n2)--(n3);
\draw[](n3)--(n4);
\draw[](n0)--(n5);
\draw[](n0)--(n6);
\draw[fill=white](n0)circle(1.6pt);
\draw[fill=white](n1)circle(1.6pt);
\draw[fill=white](n2)circle(1.6pt);
\draw[fill=white](n3)circle(1.6pt);
\draw[fill=white](n4)circle(1.6pt);
\draw[fill=white](n5)circle(1.6pt);
\draw[fill=white](n6)circle(1.6pt);
\end{tikzpicture}\end{minipage}&\begin{minipage}[c]{.475\linewidth}\centering $\widetilde B_n\ (n\geq3)$\par\begin{tikzpicture}[line width=.55pt,baseline=(current bounding box.center)]
\path[use as bounding box](-1.75,-0.54)rectangle(1.75,0.54);
\coordinate (n0) at (0.6750,-0.1850);
\coordinate (n1) at (0.2250,-0.1850);
\coordinate (n2) at (-0.2250,-0.1850);
\coordinate (n3) at (-0.6750,-0.1850);
\coordinate (n4) at (-1.1250,-0.1850);
\coordinate (n5) at (1.1250,-0.1850);
\coordinate (n6) at (0.6750,0.1850);
\draw[](n0)--(n1);
\node at (0.0000,-0.1850){$\cdots$};
\draw[](n2)--(n3);
\draw[](n3)-- node[midway,above,fill=white,inner sep=.7pt,font=\scriptsize]{$4$}(n4);
\draw[](n0)--(n5);
\draw[](n0)--(n6);
\draw[fill=white](n0)circle(1.6pt);
\draw[fill=white](n1)circle(1.6pt);
\draw[fill=white](n2)circle(1.6pt);
\draw[fill=white](n3)circle(1.6pt);
\draw[fill=white](n4)circle(1.6pt);
\draw[fill=white](n5)circle(1.6pt);
\draw[fill=white](n6)circle(1.6pt);
\end{tikzpicture}\end{minipage}\\[.5mm]
\begin{minipage}[c]{.475\linewidth}\centering $I_2(m)\ (m\geq3)$\par\begin{tikzpicture}[line width=.55pt,baseline=(current bounding box.center)]
\path[use as bounding box](-1.75,-0.45)rectangle(1.75,0.45);
\coordinate (n0) at (-0.2400,0.0000);
\coordinate (n1) at (0.2400,0.0000);
\draw[](n0)-- node[midway,above,fill=white,inner sep=.7pt,font=\scriptsize]{$m$}(n1);
\draw[fill=white](n0)circle(1.6pt);
\draw[fill=white](n1)circle(1.6pt);
\end{tikzpicture}\end{minipage}&\begin{minipage}[c]{.475\linewidth}\centering $\widetilde C_n\ (n\geq2)$\par\begin{tikzpicture}[line width=.55pt,baseline=(current bounding box.center)]
\path[use as bounding box](-1.75,-0.45)rectangle(1.75,0.45);
\coordinate (n0) at (-1.2000,0.0000);
\coordinate (n1) at (-0.7200,0.0000);
\coordinate (n2) at (-0.2400,0.0000);
\coordinate (n3) at (0.2400,0.0000);
\coordinate (n4) at (0.7200,0.0000);
\coordinate (n5) at (1.2000,0.0000);
\draw[](n0)-- node[midway,above,fill=white,inner sep=.7pt,font=\scriptsize]{$4$}(n1);
\draw[](n1)--(n2);
\node at (0.0000,0.0000){$\cdots$};
\draw[](n3)--(n4);
\draw[](n4)-- node[midway,above,fill=white,inner sep=.7pt,font=\scriptsize]{$4$}(n5);
\draw[fill=white](n0)circle(1.6pt);
\draw[fill=white](n1)circle(1.6pt);
\draw[fill=white](n2)circle(1.6pt);
\draw[fill=white](n3)circle(1.6pt);
\draw[fill=white](n4)circle(1.6pt);
\draw[fill=white](n5)circle(1.6pt);
\end{tikzpicture}\end{minipage}\\[.5mm]
\begin{minipage}[c]{.475\linewidth}\centering $F_4$\par\begin{tikzpicture}[line width=.55pt,baseline=(current bounding box.center)]
\path[use as bounding box](-1.75,-0.45)rectangle(1.75,0.45);
\coordinate (n0) at (-0.7200,0.0000);
\coordinate (n1) at (-0.2400,0.0000);
\coordinate (n2) at (0.2400,0.0000);
\coordinate (n3) at (0.7200,0.0000);
\draw[](n0)--(n1);
\draw[](n1)-- node[midway,above,fill=white,inner sep=.7pt,font=\scriptsize]{$4$}(n2);
\draw[](n2)--(n3);
\draw[fill=white](n0)circle(1.6pt);
\draw[fill=white](n1)circle(1.6pt);
\draw[fill=white](n2)circle(1.6pt);
\draw[fill=white](n3)circle(1.6pt);
\end{tikzpicture}\end{minipage}&\begin{minipage}[c]{.475\linewidth}\centering $\widetilde D_n\ (n\geq4)$\par\begin{tikzpicture}[line width=.55pt,baseline=(current bounding box.center)]
\path[use as bounding box](-1.75,-0.62)rectangle(1.75,0.62);
\coordinate (n0) at (-1.0000,-0.3200);
\coordinate (n1) at (-1.0000,0.3200);
\coordinate (n2) at (-0.5500,0.0000);
\coordinate (n3) at (0.5500,0.0000);
\coordinate (n4) at (1.0000,-0.3200);
\coordinate (n5) at (1.0000,0.3200);
\draw[](n0)--(n2);
\draw[](n1)--(n2);
\node at (0.0000,0.0000){$\cdots$};
\draw[](n3)--(n4);
\draw[](n3)--(n5);
\draw[fill=white](n0)circle(1.6pt);
\draw[fill=white](n1)circle(1.6pt);
\draw[fill=white](n2)circle(1.6pt);
\draw[fill=white](n3)circle(1.6pt);
\draw[fill=white](n4)circle(1.6pt);
\draw[fill=white](n5)circle(1.6pt);
\end{tikzpicture}\end{minipage}\\[.5mm]
\begin{minipage}[c]{.475\linewidth}\centering $E_6$\par\begin{tikzpicture}[line width=.55pt,baseline=(current bounding box.center)]
\path[use as bounding box](-1.75,-0.54)rectangle(1.75,0.54);
\coordinate (n0) at (0.0000,-0.1850);
\coordinate (n1) at (-0.4500,-0.1850);
\coordinate (n2) at (-0.9000,-0.1850);
\coordinate (n3) at (0.4500,-0.1850);
\coordinate (n4) at (0.9000,-0.1850);
\coordinate (n5) at (0.0000,0.1850);
\draw[](n0)--(n1);
\draw[](n1)--(n2);
\draw[](n0)--(n3);
\draw[](n3)--(n4);
\draw[](n0)--(n5);
\draw[fill=white](n0)circle(1.6pt);
\draw[fill=white](n1)circle(1.6pt);
\draw[fill=white](n2)circle(1.6pt);
\draw[fill=white](n3)circle(1.6pt);
\draw[fill=white](n4)circle(1.6pt);
\draw[fill=white](n5)circle(1.6pt);
\end{tikzpicture}\end{minipage}&\begin{minipage}[c]{.475\linewidth}\centering $\widetilde G_2$\par\begin{tikzpicture}[line width=.55pt,baseline=(current bounding box.center)]
\path[use as bounding box](-1.75,-0.45)rectangle(1.75,0.45);
\coordinate (n0) at (-0.4800,0.0000);
\coordinate (n1) at (0.0000,0.0000);
\coordinate (n2) at (0.4800,0.0000);
\draw[](n0)--(n1);
\draw[](n1)-- node[midway,above,fill=white,inner sep=.7pt,font=\scriptsize]{$6$}(n2);
\draw[fill=white](n0)circle(1.6pt);
\draw[fill=white](n1)circle(1.6pt);
\draw[fill=white](n2)circle(1.6pt);
\end{tikzpicture}\end{minipage}\\[.5mm]
\begin{minipage}[c]{.475\linewidth}\centering $E_7$\par\begin{tikzpicture}[line width=.55pt,baseline=(current bounding box.center)]
\path[use as bounding box](-1.75,-0.54)rectangle(1.75,0.54);
\coordinate (n0) at (-0.2250,-0.1850);
\coordinate (n1) at (-0.6750,-0.1850);
\coordinate (n2) at (-1.1250,-0.1850);
\coordinate (n3) at (0.2250,-0.1850);
\coordinate (n4) at (0.6750,-0.1850);
\coordinate (n5) at (1.1250,-0.1850);
\coordinate (n6) at (-0.2250,0.1850);
\draw[](n0)--(n1);
\draw[](n1)--(n2);
\draw[](n0)--(n3);
\draw[](n3)--(n4);
\draw[](n4)--(n5);
\draw[](n0)--(n6);
\draw[fill=white](n0)circle(1.6pt);
\draw[fill=white](n1)circle(1.6pt);
\draw[fill=white](n2)circle(1.6pt);
\draw[fill=white](n3)circle(1.6pt);
\draw[fill=white](n4)circle(1.6pt);
\draw[fill=white](n5)circle(1.6pt);
\draw[fill=white](n6)circle(1.6pt);
\end{tikzpicture}\end{minipage}&\begin{minipage}[c]{.475\linewidth}\centering $\widetilde F_4$\par\begin{tikzpicture}[line width=.55pt,baseline=(current bounding box.center)]
\path[use as bounding box](-1.75,-0.45)rectangle(1.75,0.45);
\coordinate (n0) at (-0.9600,0.0000);
\coordinate (n1) at (-0.4800,0.0000);
\coordinate (n2) at (0.0000,0.0000);
\coordinate (n3) at (0.4800,0.0000);
\coordinate (n4) at (0.9600,0.0000);
\draw[](n0)--(n1);
\draw[](n1)--(n2);
\draw[](n2)-- node[midway,above,fill=white,inner sep=.7pt,font=\scriptsize]{$4$}(n3);
\draw[](n3)--(n4);
\draw[fill=white](n0)circle(1.6pt);
\draw[fill=white](n1)circle(1.6pt);
\draw[fill=white](n2)circle(1.6pt);
\draw[fill=white](n3)circle(1.6pt);
\draw[fill=white](n4)circle(1.6pt);
\end{tikzpicture}\end{minipage}\\[.5mm]
\begin{minipage}[c]{.475\linewidth}\centering $E_8$\par\begin{tikzpicture}[line width=.55pt,baseline=(current bounding box.center)]
\path[use as bounding box](-1.75,-0.54)rectangle(1.75,0.54);
\coordinate (n0) at (-0.4500,-0.1850);
\coordinate (n1) at (-0.9000,-0.1850);
\coordinate (n2) at (-1.3500,-0.1850);
\coordinate (n3) at (0.0000,-0.1850);
\coordinate (n4) at (0.4500,-0.1850);
\coordinate (n5) at (0.9000,-0.1850);
\coordinate (n6) at (1.3500,-0.1850);
\coordinate (n7) at (-0.4500,0.1850);
\draw[](n0)--(n1);
\draw[](n1)--(n2);
\draw[](n0)--(n3);
\draw[](n3)--(n4);
\draw[](n4)--(n5);
\draw[](n5)--(n6);
\draw[](n0)--(n7);
\draw[fill=white](n0)circle(1.6pt);
\draw[fill=white](n1)circle(1.6pt);
\draw[fill=white](n2)circle(1.6pt);
\draw[fill=white](n3)circle(1.6pt);
\draw[fill=white](n4)circle(1.6pt);
\draw[fill=white](n5)circle(1.6pt);
\draw[fill=white](n6)circle(1.6pt);
\draw[fill=white](n7)circle(1.6pt);
\end{tikzpicture}\end{minipage}&\begin{minipage}[c]{.475\linewidth}\centering $\widetilde E_6$\par\begin{tikzpicture}[line width=.55pt,baseline=(current bounding box.center)]
\path[use as bounding box](-1.75,-0.54)rectangle(1.75,0.54);
\coordinate (n0) at (0.0000,-0.3700);
\coordinate (n1) at (-0.4500,-0.3700);
\coordinate (n2) at (-0.9000,-0.3700);
\coordinate (n3) at (0.4500,-0.3700);
\coordinate (n4) at (0.9000,-0.3700);
\coordinate (n5) at (0.0000,0.0000);
\coordinate (n6) at (0.0000,0.3700);
\draw[](n0)--(n1);
\draw[](n1)--(n2);
\draw[](n0)--(n3);
\draw[](n3)--(n4);
\draw[](n0)--(n5);
\draw[](n5)--(n6);
\draw[fill=white](n0)circle(1.6pt);
\draw[fill=white](n1)circle(1.6pt);
\draw[fill=white](n2)circle(1.6pt);
\draw[fill=white](n3)circle(1.6pt);
\draw[fill=white](n4)circle(1.6pt);
\draw[fill=white](n5)circle(1.6pt);
\draw[fill=white](n6)circle(1.6pt);
\end{tikzpicture}\end{minipage}\\[.5mm]
\begin{minipage}[c]{.475\linewidth}\centering $H_3$\par\begin{tikzpicture}[line width=.55pt,baseline=(current bounding box.center)]
\path[use as bounding box](-1.75,-0.45)rectangle(1.75,0.45);
\coordinate (n0) at (-0.4800,0.0000);
\coordinate (n1) at (0.0000,0.0000);
\coordinate (n2) at (0.4800,0.0000);
\draw[](n0)-- node[midway,above,fill=white,inner sep=.7pt,font=\scriptsize]{$5$}(n1);
\draw[](n1)--(n2);
\draw[fill=white](n0)circle(1.6pt);
\draw[fill=white](n1)circle(1.6pt);
\draw[fill=white](n2)circle(1.6pt);
\end{tikzpicture}\end{minipage}&\begin{minipage}[c]{.475\linewidth}\centering $\widetilde E_7$\par\begin{tikzpicture}[line width=.55pt,baseline=(current bounding box.center)]
\path[use as bounding box](-1.75,-0.54)rectangle(1.75,0.54);
\coordinate (n0) at (0.0000,-0.1850);
\coordinate (n1) at (-0.4500,-0.1850);
\coordinate (n2) at (-0.9000,-0.1850);
\coordinate (n3) at (-1.3500,-0.1850);
\coordinate (n4) at (0.4500,-0.1850);
\coordinate (n5) at (0.9000,-0.1850);
\coordinate (n6) at (1.3500,-0.1850);
\coordinate (n7) at (0.0000,0.1850);
\draw[](n0)--(n1);
\draw[](n1)--(n2);
\draw[](n2)--(n3);
\draw[](n0)--(n4);
\draw[](n4)--(n5);
\draw[](n5)--(n6);
\draw[](n0)--(n7);
\draw[fill=white](n0)circle(1.6pt);
\draw[fill=white](n1)circle(1.6pt);
\draw[fill=white](n2)circle(1.6pt);
\draw[fill=white](n3)circle(1.6pt);
\draw[fill=white](n4)circle(1.6pt);
\draw[fill=white](n5)circle(1.6pt);
\draw[fill=white](n6)circle(1.6pt);
\draw[fill=white](n7)circle(1.6pt);
\end{tikzpicture}\end{minipage}\\[.5mm]
\begin{minipage}[c]{.475\linewidth}\centering $H_4$\par\begin{tikzpicture}[line width=.55pt,baseline=(current bounding box.center)]
\path[use as bounding box](-1.75,-0.45)rectangle(1.75,0.45);
\coordinate (n0) at (-0.7200,0.0000);
\coordinate (n1) at (-0.2400,0.0000);
\coordinate (n2) at (0.2400,0.0000);
\coordinate (n3) at (0.7200,0.0000);
\draw[](n0)-- node[midway,above,fill=white,inner sep=.7pt,font=\scriptsize]{$5$}(n1);
\draw[](n1)--(n2);
\draw[](n2)--(n3);
\draw[fill=white](n0)circle(1.6pt);
\draw[fill=white](n1)circle(1.6pt);
\draw[fill=white](n2)circle(1.6pt);
\draw[fill=white](n3)circle(1.6pt);
\end{tikzpicture}\end{minipage}&\begin{minipage}[c]{.475\linewidth}\centering $\widetilde E_8$\par\begin{tikzpicture}[line width=.55pt,baseline=(current bounding box.center)]
\path[use as bounding box](-1.75,-0.54)rectangle(1.75,0.54);
\coordinate (n0) at (-0.6750,-0.1850);
\coordinate (n1) at (-1.1250,-0.1850);
\coordinate (n2) at (-1.5750,-0.1850);
\coordinate (n3) at (-0.2250,-0.1850);
\coordinate (n4) at (0.2250,-0.1850);
\coordinate (n5) at (0.6750,-0.1850);
\coordinate (n6) at (1.1250,-0.1850);
\coordinate (n7) at (1.5750,-0.1850);
\coordinate (n8) at (-0.6750,0.1850);
\draw[](n0)--(n1);
\draw[](n1)--(n2);
\draw[](n0)--(n3);
\draw[](n3)--(n4);
\draw[](n4)--(n5);
\draw[](n5)--(n6);
\draw[](n6)--(n7);
\draw[](n0)--(n8);
\draw[fill=white](n0)circle(1.6pt);
\draw[fill=white](n1)circle(1.6pt);
\draw[fill=white](n2)circle(1.6pt);
\draw[fill=white](n3)circle(1.6pt);
\draw[fill=white](n4)circle(1.6pt);
\draw[fill=white](n5)circle(1.6pt);
\draw[fill=white](n6)circle(1.6pt);
\draw[fill=white](n7)circle(1.6pt);
\draw[fill=white](n8)circle(1.6pt);
\end{tikzpicture}\end{minipage}\\[.5mm]
\bottomrule\end{tabular}
\caption{Connected elliptic and parabolic Coxeter diagrams.}
\label{figure:coxeter}\end{figure}

A symmetric matrix is \emph{indecomposable} if no simultaneous
permutation of its rows and columns puts it into two nonempty diagonal
blocks. For a Gram matrix this means that its diagram is connected.
For symmetric matrices, indecomposability is equivalent to irreducibility
in the Perron--Frobenius theorem.
The Gram matrix of a finite-volume acute-angled hyperbolic polytope is
indecomposable~\cite[Section~2]{Vinberg1985}.
A diagram is \emph{elliptic} if its Gram matrix is positive definite,
and \emph{parabolic} if each connected component has a positive
semidefinite singular Gram matrix. A diagram is \emph{hyperbolic}
if its Gram matrix has exactly one negative eigenvalue.
Connected elliptic and parabolic
diagrams correspond to spherical and Euclidean Coxeter simplices,
respectively. Their classification is shown in
Figure~\ref{figure:coxeter}~\cite{Coxeter1934}.

Every connected elliptic diagram is a tree. A connected parabolic
diagram is a tree unless it is $\widetilde A_d$ with $d\geq2$,
a cycle whose edges all have label three. The diagram
$\widetilde A_1$ has two nodes joined by a bold edge. 
The following form of the Perron--Frobenius theorem is applied to  Gram
matrices.

    \begin{thm}[Perron--Frobenius, \cite{G:1959}]\label{thm:PF}
    	A nonzero irreducible real matrix $A=(a_{ij})$ with nonnegative entries has a positive, simple eigenvalue $r$ equal to its spectral radius, with a corresponding eigenvector whose coordinates are all positive.
    \end{thm}

For a connected Coxeter diagram with at least two nodes, its Gram
matrix $G$ has diagonal entries one and nonpositive off-diagonal
entries, so $I-G$ is nonzero, irreducible, and nonnegative.
Theorem~\ref{thm:PF} therefore gives a simple smallest eigenvalue of
$G$ with a strictly positive eigenvector. In particular, if $G$ is
positive semidefinite and singular, its kernel is spanned by a
strictly positive vector. A null vector of a proper principal
submatrix, extended by zeros, would lie in this kernel, which is
impossible. Thus every proper principal submatrix is positive definite.
The one-node case is immediate.

A Coxeter diagram is \emph{minimal indefinite} if its Gram matrix
is indefinite and every proper induced subdiagram has a positive
semidefinite Gram matrix.
A connected diagram is a \emph{Lann\'er diagram} if its Gram matrix
is indefinite and every proper subdiagram is elliptic. These diagrams
correspond to compact hyperbolic Coxeter simplices and were classified
by Lann\'er \cite{Lanner1950}. A \emph{quasi-Lann\'er diagram} is a
connected indefinite diagram whose proper subdiagrams are positive
semidefinite, with at least one parabolic proper subdiagram. These
correspond to noncompact finite-volume hyperbolic Coxeter simplices.
A dotted edge is a Lann\'er diagram of order two; it corresponds to
a compact segment in $\mathbb H^1$.
The classification gives orders at most five and ten, respectively,
Gram signature $(|L|-1,1)$ for a minimal indefinite diagram $L$,
and diameter at most eight
\cite{Lanner1950,Chein1969,Prokhorov1987}.
In particular, these Gram matrices are nonsingular.
The diagrams of order at least three, which correspond to simplices
of dimension at least two, are listed in Appendix~\ref{app:lanner}.

The following results are frequently used throughout this paper.

    \begin{thm}[Theorem 2.1 in \cite{Vinberg1985}] \label{thm:signature}
        Let $G=(g_{ij})$ be an indecomposable symmetric matrix of signature $(n,1)$, where $g_{ii}=1$ and $g_{ij}\leq 0$ if $i\ne j$. Then there exists a convex hyperbolic polytope $P\subseteq\mathbb H^n$, unique up to isometry, whose Gram matrix coincides with $G$.
    \end{thm}
    \FloatBarrier

\begin{thm}[Theorems 2.1 and 3.2 in \cite{Vinberg1985}]\label{Vinberg:thm3.1}
Let $P=\bigcap_{i\in I}H_i^-\subseteq\mathbb H^n$ be an acute-angled
polytope, and let $G_J$ be the principal submatrix of its Gram matrix
indexed by $J\subseteq I$.
\begin{enumerate}[label=\textnormal{(\roman*)}]
\item If $G_J$ is positive definite, then
$F_J=P\cap\bigcap_{j\in J}H_j$ is a face of codimension $|J|$.
Conversely, for a face meeting $\mathbb H^n$, the diagram of all
facets containing it is elliptic, with order equal to its codimension.
\item Ideal vertices correspond to parabolic subdiagrams of rank
$n-1$, taking all facets through the vertex.
\end{enumerate}
\end{thm}

For a face $F$ meeting $\mathbb H^n$, let $S_F$ denote the elliptic
diagram of its containing facets. If $\dim F\geq2$ and distinct
nodes $a,b\notin S_F$ make $S_F\cup\{a,b\}$ elliptic, the facets
$F\cap H_a$ and $F\cap H_b$
meet along a ridge at angle at most $\pi/2$, with equality precisely
when $a,b$ lie in different connected components of
$S_F\cup\{a,b\}$ \cite{Vinberg1985}.

    \begin{thm}[Theorem 4.2 in \cite{Vinberg1985}] \label{Vinberg:thm4.2}
    	Let $P$ be an acute-angled polytope and $G$ its Gram matrix. Let $\mathcal{F}$ (resp. $\widetilde{\mathcal{F}}$) be the collection of elliptic principal submatrices (resp. elliptic principal submatrices and parabolic principal submatrices of rank $n-1$) of $G$, including the empty elliptic submatrix. Partially order these collections by reverse inclusion, and adjoin a formal minimal element representing the empty face.
    	Then $P$ is compact (resp. of finite volume) if and only if either of the following equivalent conditions holds:
    	\begin{enumerate}[label=\textnormal{(\roman*)}]
    		\item The polytope $P$ contains at least one finite vertex (resp. finite or ideal vertex). For every finite vertex (resp. vertex) of $P$ and any edge of $P$ emanating from it, there is precisely one other finite vertex (resp. vertex) of $P$ on that edge.
    		\item The partially ordered set $\mathcal{F}$ (resp. $\widetilde{\mathcal{F}}$) is isomorphic to the poset of some $n$-dimensional abstract combinatorial polytope.
    	\end{enumerate}
    \end{thm}

\begin{prop}[Proposition~2 in \cite{Prokhorov1987}]
	\label{prop:prokhorov-minimal-indefinite}
Every hyperbolic Coxeter diagram contains a Lann\'er or
quasi-Lann\'er subdiagram.
\end{prop}

From now on, let $P$ be a finite-volume hyperbolic Coxeter $n$-polytope.
We use the following geometric realizations of its vertex links.
At a finite vertex $u$, intersect the polytope $P$ with a sufficiently small metric
sphere centered at $u$. Radial projection of this intersection to the
unit tangent sphere at $u$ gives the spherical vertex link at $u$.
At an ideal vertex $u$, choose a horosphere $\mathcal H_u$ centered at $u$
that bounds a sufficiently small horoball disjoint from the facets
not containing $u$. The intersection $P\cap \mathcal H_u$, with the Euclidean
metric induced on $\mathcal H_u$, is the Euclidean vertex link $L_u$ at $u$,
also called the \emph{cusp section} of the polytope $P$ at $u$.
Changing $\mathcal H_u$ only rescales this Euclidean polytope.
Both constructions realize the combinatorial vertex link defined
in Subsection~\ref{subsec:face-numbers}.
The finite vertex link is a spherical Coxeter simplex, whereas
the ideal vertex link is isometric to a product of
irreducible Euclidean Coxeter simplices:
\begin{equation}\label{ext:eq:link}
L_u=P\cap \mathcal H_u\cong\Delta_1\times\cdots\times\Delta_r,
\qquad d_i=\dim\Delta_i\geq1,\qquad \sum_{i=1}^r d_i=n-1.
\end{equation}
\noindent Both types of vertex
link are simple, so each edge of the polytope $P$ lies in exactly $n-1$ facets. Thus
$P$ is edge-simple, and the incidence identities and
Nikulin--Khovanskii estimates of Subsection~\ref{subsec:face-numbers}
apply to it.

For a positive-dimensional face $F$ containing an ideal vertex $u$,
its \emph{cusp section} is $F\cap\mathcal H_u$, using the same
horosphere as in \eqref{ext:eq:link}. This construction gives an
inclusion-preserving bijection between the positive-dimensional faces
of $P$ containing $u$ and the nonempty faces of $L_u$, decreasing
dimension by one. In particular, link vertices correspond to edges
of $P$ incident to $u$, and link edges correspond to 2-faces through $u$.

Let $V_{\mathrm f}$ denote the number of finite vertices of $P$.
For an ideal vertex $u$, let $V_u$ denote the number of vertices of
its link $L_u$:
\begin{equation}\label{two:eq:link-vertices}
V_u:=f_0(L_u)=\prod_i(d_i+1).
\end{equation}
Equivalently, $V_u$ is the number of edges of $P$ incident to $u$.
A finite vertex is incident to $n$ edges.
Counting the two endpoints of every edge gives
\begin{equation}\label{two:eq:endpoint-incidence}
nV_{\mathrm f}+\sum_{u\text{ ideal}}V_u=2a_1.
\end{equation}

The vertices of the 2-faces can also be counted at each
vertex of $P$. At a finite vertex their number is
$\binom n2$, the number of edges of its simplex link.
At an ideal vertex $u$ it is $f_1(L_u)$.
Since the $(n-1)$-dimensional link $L_u$ is simple, each
of its $V_u$ vertices is incident to $n-1$ link edges,
and hence $f_1(L_u)=(n-1)V_u/2$.
Thus
\begin{equation}\label{two:eq:corner-incidences}
a_2a_2^0
=\binom n2V_{\mathrm f}
 +\frac{n-1}{2}\sum_{u\text{ ideal}}V_u.
\end{equation}
Multiplying \eqref{two:eq:endpoint-incidence} by $(n-1)/2$
recovers \eqref{two:eq:edge-face-incidence}.
The number $V_u$ will normalize ideal-vertex totals in both
comparisons, although the original ideal vertices need not be simple.

\subsection{\texorpdfstring{$h$-polynomials of product links}{h-polynomials of product links}}
\label{common:subsec:product-links}
\label{common:subsec:asymmetry}
We retain the notation of \eqref{ext:eq:link}, writing
$L=L_u$ and $V=V_u$.
Every face of a product is a product of faces, with dimensions adding.
The defining face sum in \eqref{ext:eq:h-polynomial} therefore gives
$h_{D\times D'}(t)=h_D(t)h_{D'}(t)$.
Since $f_k(\Delta_i)=\binom{d_i+1}{k+1}$, the binomial theorem gives
\begin{equation}\label{ext:eq:link-product-h}
 h_L(t)=\prod_{i=1}^r(1+t+\cdots+t^{d_i}),\qquad
 V=h_L(1)=\prod_i(d_i+1).
\end{equation}
For the lower half of this polynomial, define
\begin{equation}\label{ext:eq:g-convention}
 g_0(L)=h_0(L)=1,\qquad
 g_j(L)=h_j(L)-h_{j-1}(L)
 \quad(1\leq j\leq\lfloor (n-1)/2\rfloor).
\end{equation}
Equivalently, take $h_{-1}(L)=0$ and use the difference formula also
at $j=0$.

\begin{lem}
\label{ext:lem:link-polynomial}
For a product of simplices of total dimension $n-1$,
\begin{equation}\label{ext:eq:strings}
 h_L(t)=\sum_{j=0}^{\lfloor (n-1)/2\rfloor}
 g_j(L)t^j(1+t+\cdots+t^{n-1-2j}),\qquad g_j(L)\geq0.
\end{equation}
In particular,
$h_0(L)\leq h_1(L)\leq\cdots\leq h_{\lfloor (n-1)/2\rfloor}(L)$.
\end{lem}
\begin{proof}
For integers $a,b\geq0$,
\begin{equation}\label{ext:eq:two-factor-intervals}
 (1+t+\cdots+t^a)(1+t+\cdots+t^b)
 =\sum_{p=0}^{\min(a,b)}t^p(1+t+\cdots+t^{a+b-2p}).
\end{equation}
The identity follows by comparing coefficients.

Applying this identity successively to the factors in
\eqref{ext:eq:link-product-h} yields a decomposition of the form
\eqref{ext:eq:strings}, with nonnegative integers $\gamma_j$ in place
of $g_j(L)$. For $0\leq k\leq\lfloor(n-1)/2\rfloor$, comparing
coefficients gives $h_k(L)=\sum_{j=0}^k\gamma_j$, and hence
$\gamma_j=h_j(L)-h_{j-1}(L)=g_j(L)$.
\end{proof}

The polynomial $t^j(1+t+\cdots+t^{n-1-2j})$ has $n-2j$
monomials. Evaluating \eqref{ext:eq:strings} at $t=1$ gives
$V=\sum_j(n-2j)g_j(L)$. Define the normalized weights and two
weighted averages by
\begin{equation}\label{ext:eq:normalized-g}
 \pi_j:=\frac{(n-2j)g_j(L)}V,\qquad
 \pi_j\geq0,\qquad\sum_j\pi_j=1,
\end{equation}
\begin{equation}\label{common:eq:weighted-averages}
 m=m(d_1,\ldots,d_r):=\sum_jj\pi_j,\qquad
 N:=\sum_jj(j-1)\pi_j,\qquad
 e_2:=\sum_{i<\ell}d_id_\ell.
\end{equation}
Unless another range is stated, sums over $j$ for this fixed product
run from $0$ to $\lfloor (n-1)/2\rfloor$.
Each $\pi_j$ is the proportion of $h_L(1)=V$ contributed by the
summand $g_j(L)t^j(1+t+\cdots+t^{n-1-2j})$
in \eqref{ext:eq:strings}.
The quantities $m$ and $N$ are ordinary finite weighted averages.
They depend only on the factor dimensions and are symmetric in them.
For one simplex, $m=N=0$.

\begin{lem}
\label{ext:lem:moments}
For every product of simplices of total dimension $n-1$,
\begin{equation}\label{two:eq:centered-identity}
 \sum_j(n-2j)^2\pi_j=1+\sum_i d_i(d_i+2).
\end{equation}
Equivalently,
\begin{equation}\label{ext:eq:moments}
 \sum_jj(n-j)\pi_j=\frac{e_2}{2},\qquad
 m\geq\frac{e_2}{2n},\qquad N\geq m^2-m.
\end{equation}
In addition,
\begin{equation}\label{ext:eq:first-second-control}
 m\leq N+1,\qquad m^2\leq2N+1.
\end{equation}
\end{lem}
\begin{proof}
We compute the centered second moment of the coefficients of $h_L$
in two ways. For $a\geq0$,
\begin{equation}\label{ext:eq:elementary-centered-sums}
 \frac1{a+1}\sum_{k=0}^a(k-a/2)=0,\qquad
 \frac1{a+1}\sum_{k=0}^a(k-a/2)^2=\frac{a(a+2)}{12}.
\end{equation}
In the product expansion \eqref{ext:eq:link-product-h}, the exponents
are sums of independent choices $0\leq k_i\leq d_i$. The mixed
centered terms vanish by the first identity, so
\begin{equation}\label{ext:eq:product-centered-square}
 \frac1V\sum_{k=0}^{n-1}h_k(L)(k-(n-1)/2)^2
 =\frac1{12}\sum_i d_i(d_i+2).
\end{equation}
Each summand in \eqref{ext:eq:strings} has exponents from $j$ to
$n-1-j$, with the same center $(n-1)/2$. Its centered second moment
is $((n-2j)^2-1)/12$. Hence
\begin{equation}\label{ext:eq:interval-centered-square}
 \frac1V\sum_{k=0}^{n-1}h_k(L)(k-(n-1)/2)^2
 =\frac1{12}\sum_j\pi_j\bigl((n-2j)^2-1\bigr).
\end{equation}
Comparing the two expressions proves \eqref{two:eq:centered-identity}.
Since $\sum_i d_i(d_i+2)=n^2-1-2e_2$, expansion gives
$\sum_j j(n-j)\pi_j=e_2/2$; the inequality $j(n-j)\leq nj$
then gives $m\geq e_2/(2n)$.

Finally, Cauchy--Schwarz yields $m^2\leq\sum_jj^2\pi_j=N+m$.
Summing $j\leq j(j-1)+1$ and $j^2\leq2j(j-1)+1$ against
$\pi_j$ proves \eqref{ext:eq:first-second-control}.
\end{proof}

For example, let $L=\Delta^3\times\Delta^2$, where $\Delta^k$
denotes a simplex of dimension $k$. Then $n=6$, and
\begin{equation}\label{ext:eq:link-example-polynomial}
\begin{aligned}
 h_L(t)&=(1+t+t^2+t^3)(1+t+t^2)\\
 &=1+2t+3t^2+3t^3+2t^4+t^5\\
 &=(1+t+\cdots+t^5)+(t+t^2+t^3+t^4)+(t^2+t^3).
\end{aligned}
\end{equation}
Thus $g_0=g_1=g_2=1$. The three summands in the last line have
six, four, and two monomials, so $V=12$ and
\[
 (\pi_0,\pi_1,\pi_2)=\left(\frac12,\frac13,\frac16\right),
 \qquad m(3,2)=\frac23,\qquad N=\frac13,\qquad e_2=6.
\]
The centered-square identity reads
\[
 \frac12\,6^2+\frac13\,4^2+\frac16\,2^2
 =24=1+3(3+2)+2(2+2),
\]
and $\sum_jj(6-j)\pi_j=3=e_2/2$.

We now apply the truncation construction to the closure of $P$ in the Klein
model. Its nonsimple vertices are ideal, and the inserted facet has
the combinatorial type of the corresponding cusp link $L_u$.
The simple polytope $Q$ satisfies $h_j(Q)=h_{n-j}(Q)$, while each
simple $(n-1)$-dimensional link satisfies
$h_{n-j}(L_u)=h_{j-1}(L_u)$.
Subtracting \eqref{ext:eq:truncation} in degrees $j$ and $n-j$ gives
\begin{equation}\label{ext:eq:asymmetry}
\begin{aligned}
 h_{n-j}(P)-h_j(P)
 &=\sum_{u\text{ ideal}}\bigl(h_j(L_u)-h_{n-j}(L_u)\bigr)\\
 &=\sum_{u\text{ ideal}}g_j(L_u),
 \qquad 1\leq j\leq\lfloor(n-1)/2\rfloor.
\end{aligned}
\end{equation}
Initially the truncation sum runs over nonsimple vertices. It may
be extended to all ideal vertices, since a simplex link has
$g_j=0$ for every $j\geq1$.
The identity is not asserted at $j=0$, where the truncation formula
has a different form. Both later comparisons retain ideal-vertex
terms from this same identity, with different coefficients.

\part{Improving Prokhorov's dimension bound}\label{part:edge}
\section{The cellwise--edgewise comparison}\label{edge:sec:comparison}
In this part, an incidence consists of a 3-face of $P$ and one of
its edges. The term \emph{cellwise} refers to sums on a fixed 3-face.
The weights $w$ and their total $W$ refer to these incidences until
Part~\ref{part:vertex}. We first establish the following bound.
\begin{thm}\label{edge:thm:dimension-bound}
A finite-volume hyperbolic Coxeter polytope has dimension at most $162$.
\end{thm}
The proof retains Prokhorov's four local configurations and sharpens
both the edgewise weight estimate and the estimates for average face numbers.
\subsection{Three-dimensional faces and square cusp sections}
Every 2-face of the product link \eqref{ext:eq:link} has, up to a
permutation of factors, one of the forms
\begin{equation}\label{eq:link-two-face-types}
	\Delta_i^2\times\prod_{k\ne i}\{v_k\},
	\qquad
	\Delta_i^1\times\Delta_j^1\times\prod_{k\ne i,j}\{v_k\},\quad i<j,
\end{equation}
where $\Delta_i^2$ and $\Delta_i^1$ denote a two-dimensional face and
an edge of the simplex $\Delta_i$, respectively, and $v_k$ is a vertex
of $\Delta_k$.

Every 3-face $F$ of $P$ is a finite-volume acute-angled
hyperbolic polyhedron, although it need not itself be a Coxeter polyhedron.
Its finite vertices are three-valent. At an ideal vertex, its link
is a 2-face of a product of Euclidean simplices, hence a triangle
or a rectangle. Thus an ideal vertex has valence three or four,
and every four-valent vertex is ideal. Let $v_j(F)$ denote the number of vertices of
valence $j$ in $F$. From
$3v_3(F)+4v_4(F)=2f_1(F)$ and Euler's identity
$v_3(F)+v_4(F)-f_1(F)+f_2(F)=2$, we obtain
\begin{equation}\label{eq:euler-face}
 f_1(F)=3f_2(F)-6-v_4(F).
\end{equation}

\FloatBarrier

\subsubsection{The right-angled triangular bipyramid}

One particular 3-face type, denoted by $\mathcal{P}$, will
be exceptional for the weighting rule introduced below. Combinatorially, it is the triangular bipyramid shown in
Figure~\ref{fig:small-faces}(i), with three four-valent
vertices and two three-valent vertices. In any finite-volume acute-angled
hyperbolic realization of this combinatorial type, the four-valent
vertices are ideal and have rectangular horospherical links.
Since every edge is incident to one of these vertices, all dihedral
angles are right. The two three-valent vertices are finite: an
ideal one would have a Euclidean triangular link with three right
angles, which is impossible. Thus the realization is right-angled,
with three ideal vertices and two finite vertices. Uniqueness up to isometry
follows from \cite{Andreev1970Finite}. Its $f$-vector is
\begin{equation}\label{eq:bad-data}
 (f_0(\mathcal{P}),f_1(\mathcal{P}),f_2(\mathcal{P}))=(5,9,6).
\end{equation}

\begin{figure}[!htb]
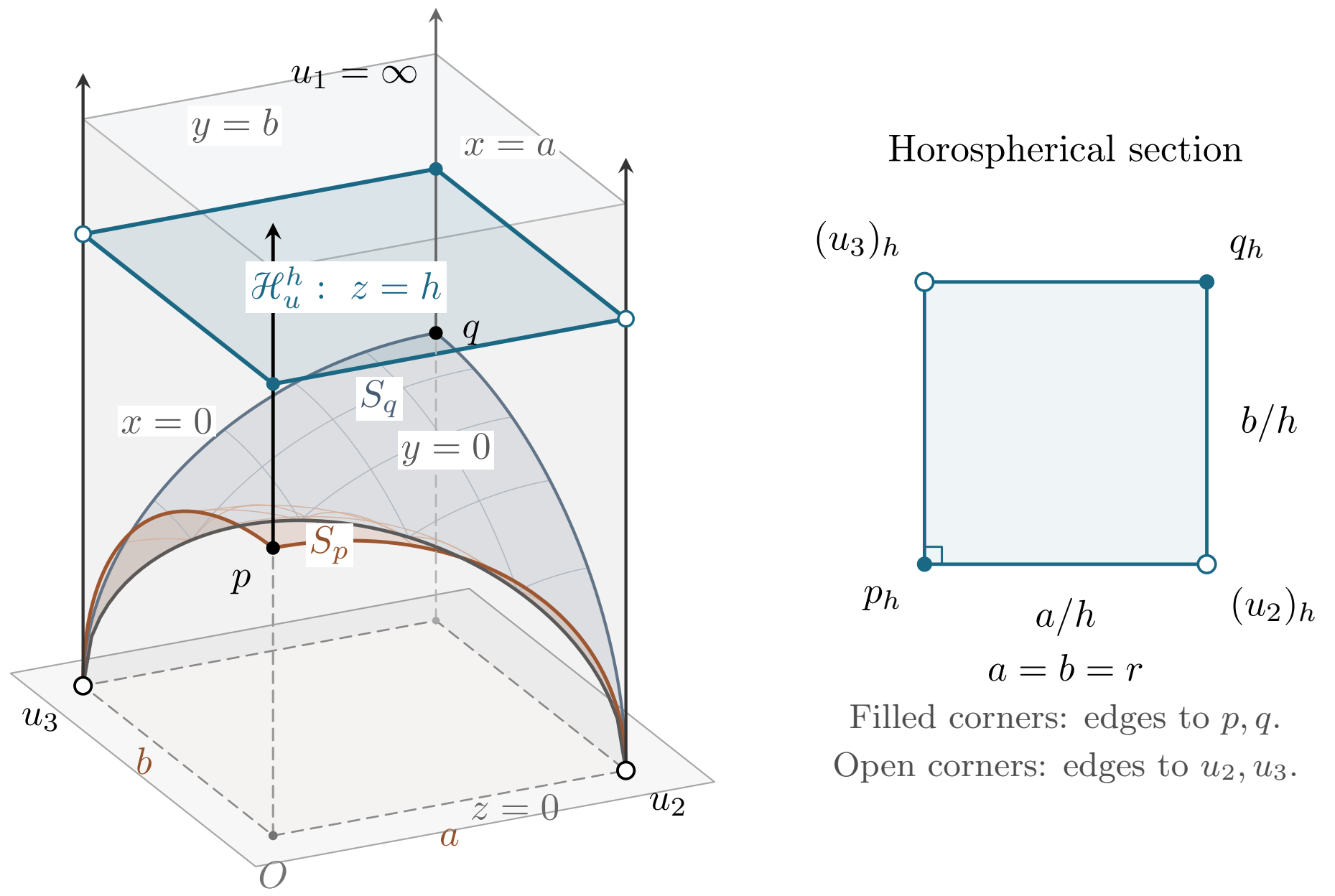

\centering

\caption[Square cusp sections of $\mathcal{P}$.]{Square cusp sections of $\mathcal{P}$,
with ideal vertices $u_1,u_2,u_3$ and finite apices $p,q$.
In the upper half-space model, put $u_1=\infty$, with four vertical
supporting planes $x=0$, $x=a$, $y=0$, $y=b$, where $a,b>0$ are
the horizontal separations of the two pairs. Let $S_p,S_q$ be the
hemispheres supporting $pu_2u_3,qu_2u_3$, respectively. Since $S_p$
is orthogonal to $x=0$ and $y=0$, both planes contain its center
$O=(0,0,0)$. Its boundary contains $u_2=(a,0,0)$ and $u_3=(0,b,0)$,
so $a=b=r$, where $r$ is its radius. Similarly, orthogonality of $S_q$
to $x=a$ and $y=b$ puts its center at $(a,b,0)$. Its boundary contains
the same two ideal vertices, so its radius is also $r$.
For sufficiently large $h$,
the horosphere $\mathcal H_u^h=\{z=h\}$ cuts out a square of side
$a/h=b/h=r/h$. In the right-hand panel, $v_h$ denotes the intersection
of the edge $u_1v$ with $\mathcal H_u^h$, for $v\in\{p,q,u_2,u_3\}$.}
\label{fig:bad-square-3d}
\end{figure}

For a 3-face $F\cong\mathcal{P}$ with ideal vertex $u$, its cusp section
$F\cap \mathcal H_u$ is the \emph{exceptional square section} shown in
Figure~\ref{fig:bad-square-3d}. It has the form
\begin{equation}\label{eq:bad-square-product}
F\cap \mathcal H_u\cong\Delta_i^1\times\Delta_j^1\times\prod_{k\ne i,j}\{v_k\},
\qquad |\Delta_i^1|=|\Delta_j^1|.
\end{equation}
The lengths are measured at the actual scales of the factors in $L_u$. Label each vertex
of $L_u$ as finite or ideal according to the other endpoint of its
corresponding edge of $P$. These labels alternate around every exceptional
square.
Counting all such squares therefore bounds the number of exceptional
3-faces through $u$ from above.

Following Prokhorov, call 3-faces isometric to $\mathcal{P}$ \emph{bad}
and all other 3-faces \emph{good}. Define the proportion of bad
3-faces by
\[
 b=\frac{\#\{F:\dim F=3,\ F\cong\mathcal{P}\}}{a_3},\qquad 0\leq b\leq1.
\]
Averaging \eqref{eq:euler-face} and using $v_4(\mathcal{P})=3$ gives
\begin{equation}\label{eq:euler-average}
 a_3^1=3a_3^2-6-\frac1{a_3}\sum_{\dim F=3}v_4(F)
 \leq 3a_3^2-6-3b.
\end{equation}

\FloatBarrier
\subsubsection{Equal-length edges in Euclidean Coxeter simplices}\label{ext:app:equal-length}
The square sections in \eqref{eq:bad-square-product} require edges of
equal length in two simplex factors. To bound the number of such
sections, we first bound the number of edges of each length in one factor.
For a Euclidean Coxeter simplex $\Delta$, let $m_\Delta(\rho)$ be the
number of edges with squared length $\rho>0$, and let
$m_\Delta(v;\rho)$ count those incident to a specified vertex $v$.
%We use squared lengths to simplify the calculations.

\begin{lem}\label{ext:lem:equal-length}
Let $\Delta$ be a Euclidean Coxeter simplex of dimension $d\geq1$.
For every prescribed squared length $\rho>0$,
$m_\Delta(\rho)\leq6(d+1)/5$. Unless $\Delta$ has type $\wt D_4$,
the stronger bound $d+1$ holds. In type $\wt D_4$, six of the ten
edges have the same length. In particular, $6/5$ is the optimal
constant in a uniform bound proportional to the number of vertices.
\end{lem}

\begin{proof}
The affine classification consists of four infinite families
$\wt A_{d},\wt B_{d},\wt C_{d},\wt D_{d}$, and five
exceptional types $\wt E_6,\wt E_7,\wt E_8,\wt F_4,\wt G_2$
\cite[Appendix~A1, p.~298]{BjornerBrenti2005}.
Each type determines a Euclidean Coxeter simplex up to similarity.
Explicit alcove vertex coordinates are tabulated in
\cite[Table~1, p.~151]{ChoudharyKachanovichWintraecken2020}.
For a general description in terms of fundamental coweights, see
\cite[Ch.~VI, \S\,2, no.~2, Corollary to Proposition~5]{Bourbaki1968}.
After rescaling, these coordinates take the following form. Write
$e_1,e_2,\ldots$ for the standard orthonormal basis in each ambient space.
Type $\wt A_1$ has only one edge.

\smallskip\noindent\emph{Type $\wt A_{d}$, $d\geq2$.}
In $\mathbb R^{d+1}$, take $\Delta$ to be the convex hull of the points
\[
v_i=e_1+\cdots+e_i-\frac{i}{d+1}(1,\ldots,1),\qquad 0\leq i\leq d,
\]
where $v_0=0$. For $0\leq i<j\leq d$,
\[
\begin{aligned}
v_j-v_i&=e_{i+1}+\cdots+e_j-\frac{j-i}{d+1}(1,\ldots,1),\\
|v_j-v_i|^2&=(j-i)-\frac{(j-i)^2}{d+1}
            =\frac{(j-i)(d+1-j+i)}{d+1}.
\end{aligned}
\]

The squared length depends only on the cyclic separation
$k=\min\{j-i,d+1-j+i\}$, and increases strictly for
$1\leq k\leq\lfloor(d+1)/2\rfloor$.
For $k<(d+1)/2$, the two index differences $k$ and $d+1-k$
give $(d+1-k)+k=d+1$ pairs. When $k=(d+1)/2$, they coincide
and give $(d+1)/2$ pairs. Thus

\begingroup
\small
\linespread{1}\selectfont
\begin{equation}\label{ext:eq:A-length-counts}
\begin{array}{c|c|c}
\toprule[1pt]
\text{Cyclic separation }k&\text{Squared length }\rho&m_\Delta(\rho)\\\hline
1\leq k<(d+1)/2&k(d+1-k)/(d+1)&d+1\\
k=(d+1)/2\quad(d\text{ odd})&(d+1)/4&(d+1)/2\\
\bottomrule[1pt]
\end{array}
\end{equation}
\endgroup
In every dimension $d\geq2$, the maximum multiplicity is $d+1$.

\smallskip\noindent\emph{Types $\wt B_{d},\wt C_{d},\wt D_{d}$.}
Write $p_i=e_1+\cdots+e_i$, with $p_0=0$. The vertex sets are
\begin{align}
\wt C_{d}:&\quad 0,p_1,\ldots,p_{d},\nonumber\\
\wt B_{d}:&\quad 0,2e_1,p_2,\ldots,p_{d},\label{ext:eq:classical-coordinate-models}\\
\wt D_{d}:&\quad 0,2e_1,p_2,\ldots,p_{d-2},
                  p_{d-1}-e_{d},p_{d-1}+e_{d}.\nonumber
\end{align}
These coordinates are twice those in the corresponding rows of
\cite[Table~1]{ChoudharyKachanovichWintraecken2020}.
Here $d\geq2$ for type $C$, $d\geq3$ for type $B$, and $d\geq4$ for type $D$.
Calculating the squared distances between all pairs of vertices and
grouping equal values gives the following multiplicities.
\begingroup
\small
\linespread{1}\selectfont
\renewcommand{\arraystretch}{0.92}
\setlength{\extrarowheight}{0pt}
\setlength{\arraycolsep}{3.5pt}
\renewcommand{\paddedmathcell}[1]{\vcenter{\kern2pt\hbox{$\displaystyle#1$}\kern2pt}}
\begin{equation}\label{ext:eq:classical-length-counts}
\begin{array}{c|c|c!{\vrule width 1.2pt}c|c|c}
\toprule[1pt]
\text{Type}&\text{Squared length }\rho&m_\Delta(\rho)
 &\text{Type}&\text{Squared length }\rho&m_\Delta(\rho)\\\hline
\paddedmathcell{\begin{array}{@{}c@{}}\wt C_d\\(d\geq2)\end{array}}
 &1\leq\rho\leq d&d+1-\rho
 &\wt D_4
 &\paddedmathcell{\begin{array}{@{}c@{}}2\\4\end{array}}
 &\paddedmathcell{\begin{array}{@{}c@{}}4\\6\end{array}}\\\hline
\wt B_3
 &\paddedmathcell{\begin{array}{@{}c@{}}1,\ 4\\2,\ 3\end{array}}
 &\paddedmathcell{\begin{array}{@{}c@{}}1\\2\end{array}}
 &\wt D_5
 &\paddedmathcell{\begin{array}{@{}c@{}}1\\2,\ 3,\ 5\\4\end{array}}
 &\paddedmathcell{\begin{array}{@{}c@{}}1\\4\\2\end{array}}\\\hline
\wt B_4
 &\paddedmathcell{\begin{array}{@{}c@{}}1,\ 3\\2,\ 4\end{array}}
 &\paddedmathcell{\begin{array}{@{}c@{}}2\\3\end{array}}
 &\wt D_6
 &\paddedmathcell{\begin{array}{@{}c@{}}1\\2\\3,\ 6\\4\end{array}}
 &\paddedmathcell{\begin{array}{@{}c@{}}2\\5\\4\\6\end{array}}\\\hline
\paddedmathcell{\begin{array}{@{}c@{}}\wt B_d\\(d\geq5)\end{array}}
 &\paddedmathcell{\begin{array}{@{}c@{}}
 1,\ 4\\2\leq\rho\leq d-1,\ \rho\ne4\\d
 \end{array}}
 &\paddedmathcell{\begin{array}{@{}c@{}}d-2\\d+1-\rho\\2\end{array}}
 &\paddedmathcell{\begin{array}{@{}c@{}}\wt D_d\\(d\geq7)\end{array}}
 &\paddedmathcell{\begin{array}{@{}c@{}}
 1\\4\\2\leq\rho\leq d-3,\ \rho\ne4\\d-2,\ d
 \end{array}}
 &\paddedmathcell{\begin{array}{@{}c@{}}d-4\\d-1\\d+1-\rho\\4\end{array}}\\
\bottomrule[1pt]
\end{array}
\end{equation}
\endgroup
Here $\rho$ ranges over positive integers. Each listed squared length
has the multiplicity shown in the corresponding $m_\Delta(\rho)$ column.
All unlisted squared lengths have multiplicity zero.
The maximum multiplicity is $d$ for type $\wt C_d$ and $d-1$ for
type $\wt B_d$. For type $\wt D_d$, it is $d-1$ except when $d=4$
or $6$, in which case it is $6$.
In particular, the four vertices $0,2e_1,p_3-e_4,p_3+e_4$ of the
$\wt D_4$ simplex have all six mutual squared distances equal to $4$.
The four remaining edges have squared length $2$.

\smallskip\noindent\emph{Exceptional types.}
The origin and the vectors in each row below are the vertices of the
indicated simplex:
\begingroup
\small
\linespread{1}\selectfont
\setlength{\LTpre}{6pt}
\setlength{\LTpost}{6pt}
\begin{longtable*}{>{$}c<{$}|>{$}l<{$}}
\toprule[1pt]
\text{Type}&\text{Nonzero vertices}\\\hline
\endfirsthead
\toprule[1pt]
\text{Type}&\text{Nonzero vertices (continued)}\\\hline
\endhead
\bottomrule[1pt]
\endfoot
\wt E_6&e_1,\ e_1+\sqrt3e_2,\ e_3,\ e_3+\sqrt3e_4,\ e_5,\ e_5+\sqrt3e_6\\\hline
\wt E_7&\paddedmathcell{\begin{aligned}
 &\sqrt3e_1,\ \sqrt2e_2,\ \sqrt2e_2+2e_3,\ \sqrt2e_2+2e_3+2\sqrt3e_4,\\
 &\sqrt2e_5,\ \sqrt2e_5+2e_6,\ \sqrt2e_5+2e_6+2\sqrt3e_7
\end{aligned}}\\\hline
\wt E_8&\paddedmathcell{\begin{aligned}
 &2\sqrt5e_1,\ \sqrt{15}e_2,\ \sqrt{15}e_2+3\sqrt5e_3,\ 2\sqrt3e_4,\ 2\sqrt3e_4+3\sqrt2e_5,\\
 &2\sqrt3e_4+3\sqrt2e_5+\sqrt{30}e_6,\ 2\sqrt3e_4+3\sqrt2e_5+\sqrt{30}e_6+2\sqrt{15}e_7,\\
 &2\sqrt3e_4+3\sqrt2e_5+\sqrt{30}e_6+2\sqrt{15}e_7+6\sqrt5e_8
\end{aligned}}\\\hline
\wt F_4&\sqrt3e_1,\ \sqrt3e_1+e_2,\ \sqrt3e_1+e_2+\sqrt2e_3,\ \sqrt3e_1+e_2+\sqrt2e_3+\sqrt6e_4\\\hline
\wt G_2&e_1,\ e_1+\sqrt3e_2\\
\end{longtable*}
\endgroup
Direct calculation of the squared distances between these vertices gives
the following maximum multiplicities:
\begingroup
\small
\linespread{1}\selectfont
\[
\begin{array}{c|rrrrr}
\toprule[1pt]
\Delta&\wt E_6&\wt E_7&\wt E_8&\wt F_4&\wt G_2\\\hline
\max_\rho m_\Delta(\rho)&6&3&3&2&1\\
\bottomrule[1pt]
\end{array}
\]
\endgroup
For every type other than $\wt D_4$, the maximum multiplicity is at
most $d+1$. Type $\wt D_4$ has five vertices and maximum multiplicity
$6$, which proves both the uniform bound and the optimality of $6/5$.
\end{proof}

\FloatBarrier
\subsection{Weighting by short distances}\label{subsec:distance-weights}
A weight is assigned to an incidence $(F,e)$, where $F$ is a
three-dimensional face and $e$ is one of its edges. Equivalently,
it is assigned to a dihedral angle of $F$. For each edge $e$, let
$T_e$ be the elliptic diagram of the $n-1$ facets containing it.
For a nonnegative rule $w$, write
\begin{equation}\label{eq:weight-sums}
\begin{aligned}
 \sigma(F)&=\sum_{e\subset F}w(F,e),\\W(T_e)&=\sum_{F\supset e}w(F,e),\\
 a_w&=\frac1{a_3}\sum_{\dim F=3}\sigma(F)
     =\frac1{a_3}\sum_eW(T_e).
\end{aligned}
\end{equation}
Here $\sigma(F)$ is the \emph{cellwise} sum over the edges of a
fixed 3-face $F$, and $W(T_e)$ is the \emph{edgewise} sum over the
3-faces containing a fixed edge $e$. The rules considered below determine
the edge sum from $T_e$, which accounts for the notation $W(T_e)$.
The quantity $a_w$ is the \emph{average total weight of a 3-face}.
The equality defining
$a_w$ follows by counting every incidence $(F,e)$ once in each order.

The eventual contradiction compares a lower and an upper bound for this
same average $a_w$. To obtain the upper bound, we seek an estimate
$W(T_e)\leq C(n-1)$ for every edge, with $C$ independent of $n$.
Thus $C$ bounds the total weight at an edge divided by $n-1$, the
number of nodes in its elliptic diagram.
Using \eqref{eq:incidences}, we have
\begin{equation}\label{eq:normalized-edge-bound}
 a_w\leq C(n-1)\frac{a_1}{a_3}
 =C(n-1)\frac{a_3^1}{\binom{n-1}{2}}
 =\frac{2C}{n-2}\,a_3^1.
\end{equation}
This is where the factor $n^{-1}$ enters: each edge belongs to
$\binom{n-1}{2}$ 3-faces, so the edge--face incidence count divides
the linear bound $C(n-1)$ by a quadratic factor. Since
\eqref{eq:A-E-bounds} gives $a_3^1<12+O(n^{-1})$, the resulting
upper bound is $24C/n+O(n^{-2})$, which tends to zero.

For the lower bound, we seek a rule giving $a_w\geq7-a_3^2-b$.
The choice of $7$ is motivated by the upper bound
$a_3^2<6+O(n^{-1})$ from \eqref{eq:A-E-bounds} and Prokhorov's
independent estimate $b\leq3/4+O(n^{-1})$, recalled in
Subsection~\ref{subsec:weights}. These upper bounds tend to $6$ and
$3/4$, respectively, leaving the positive gap $1/4$.
Thus the comparison has the form
\[
 \frac14-O(n^{-1})\leq a_w
 \leq\frac{24C}{n}+O(n^{-2}).
\]
The left-hand bound is positive for sufficiently large $n$, while the
right-hand bound tends to zero. The corresponding exact inequalities yield the dimension bound in
Subsection~\ref{subsec:weights}. This suggests looking for a weighting
with $\sigma(F)\geq7-f_2(F)$ for every good 3-face $F$.

For a fixed 3-face $F$, let $S_F$ be the elliptic diagram of the
$n-3$ ambient facets containing it. Each two-dimensional face of $F$
is cut out by one further ambient facet, whose node we call a
\emph{marked node}. A pair $\{a,b\}$ of such nodes is
\emph{admissible} if the corresponding two-dimensional faces meet
in an edge $e$ of $F$. In this case, $T_e=S_F\cup\{a,b\}$,
and $\{a,b\}$ is the \emph{marked pair} associated with $(F,e)$.

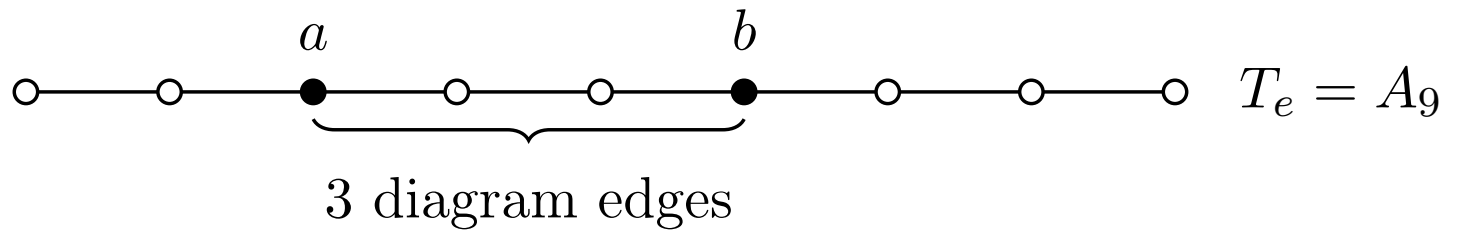
\begin{figure}[!htbp]
\centering
\begin{tikzpicture}[x=.95cm,y=1cm,line width=.65pt]
\draw(0,0)--(8,0);
\foreach\i in {0,1,2,3,4,5,6,7,8}{\draw[fill=white](\i,0)circle(2.2pt);}
\fill(2,0)circle(2.2pt);\fill(5,0)circle(2.2pt);
\node[above=4pt]at(2,0){$a$};\node[above=4pt]at(5,0){$b$};
\draw[decorate,decoration={brace,mirror,amplitude=4pt}](2,-.18)--(5,-.18)
node[midway,below=7pt]{$3$ diagram edges};
\node[right]at(8.3,0){$T_e=A_9$};
\end{tikzpicture}
\caption{Two marked nodes in an elliptic edge diagram. Their distance
is three, independently of the Coxeter labels. Diagram edges record
nonorthogonality of ambient facets, not edges of the polytope.}
\label{fig:marked-path}
\end{figure}

Conversely, for a fixed edge $e$, choosing two distinct nodes $a,b$
of $T_e$ leaves $n-3$ nodes cutting out a unique 3-face $F\supset e$.
Thus
\begin{equation}\label{eq:marked-bijection}
 \{F:\dim F=3,\ e\subset F\}
 \longleftrightarrow
 \{\{a,b\}\subset V(T_e):a\ne b\}.
\end{equation}
For nodes $a,b$ of a Coxeter diagram $T$, define $\dist_T(a,b)$
to be the number of edges in a shortest path from $a$ to $b$.
Set $\dist_T(a,b)=\infty$ if the nodes lie in different connected
components. Each edge counts as one step, regardless of its Coxeter label.
In Part~\ref{part:edge}, distances are measured in the elliptic edge diagram $T_e$
associated with the incidence $(F,e)$.
Figure~\ref{fig:marked-path} illustrates this definition.

The rule must give fixed positive totals on small 3-faces while
keeping $W(T_e)=O(n)$. Constant weights cannot do both: weight one
gives $W(T_e)=\binom{n-1}{2}$, whereas weight $1/n$ gives only
$6/n$ on a tetrahedron, whose required total is three.
Short diagram distance provides a criterion compatible with both estimates.

For the edgewise estimate, the number of node pairs within a fixed
distance in an elliptic diagram is bounded by a constant times its
number of nodes. More precisely, for an elliptic diagram $T$ and a
fixed positive integer $R$, the classification gives
\begin{equation}\label{eq:elliptic-pairs}
	N_R(T):=\#\{\{a,b\}\subset V(T):a\ne b,\ \dist_T(a,b)\leq R\}
	\leq R|T|.
\end{equation}
See \cite[Proposition~3]{Prokhorov1987}. For a chain of $r>R$ nodes,
there are $r-j$ unordered pairs at distance $j$, so
$N_R(T)=\sum_{j=1}^R(r-j)\leq Rr$. Thus a fixed distance cutoff
retains only linearly many of the quadratically many possible pairs.
In particular, any rule bounded by one and supported on distances at
most $R$ satisfies $W(T_e)\leq R(n-1)$. This controls the edgewise
sum uniformly in dimension. Prokhorov takes $R=15$ and uses the
following rule. The particular cutoffs $7$ and $15$, and the weights
$1$ and $1/3$, also account for the cellwise lower bounds proved in
Proposition~\ref{prop:facewise}.
A subscript $\mathrm P$ indicates Prokhorov's rule and the sums computed
with it.
\begin{equation}\label{eq:prokhorov-weight}
	w_{\mathrm P}(F,e)=
	\begin{cases}
		1,&1\leq \dist_{T_e}(a,b)\leq7,\\
		1/3,&8\leq \dist_{T_e}(a,b)\leq15,\\
		0,&\dist_{T_e}(a,b)>15\text{ or }\dist_{T_e}(a,b)=\infty.
	\end{cases}
\end{equation}
A pair at distance at most seven is counted in both $N_7$ and
$N_{15}$, whereas a pair at distance eight through fifteen is
counted only in $N_{15}$. Consequently the edge sum is
\begin{equation}\label{eq:old-edge-constant}
	W_{\mathrm P}(T_e)=\tfrac23N_7(T_e)+\tfrac13N_{15}(T_e)
	\leq C_{\mathrm P}(n-1),\qquad C_{\mathrm P}=\frac{29}{3}.
\end{equation}

For the cellwise estimates, a 3-face with four, five, or six facets
needs total weight at least $3,2,$ or $1$, respectively, to meet the
bound $7-f_2(F)$. Its small number of facets alone gives no distance
bound, since $S_F$ has $n-3$ nodes through which paths between marks
may pass. The Lann\'er and quasi-Lann\'er classification supplies the
additional restriction: these diagrams have uniformly bounded diameter.
For the good 3-face types in Figure~\ref{fig:small-faces}, the facet
incidences place selected marks in suitable subdiagrams of these types.
The short paths they supply yield sufficiently many pairs with
positive weight in the corresponding elliptic edge diagrams to meet
the required totals, independently of the ambient dimension.

For a bad 3-face $F\cong\mathcal{P}$, every dihedral angle is right.
By the face-angle criterion in
Subsection~\ref{subsec:hyperbolic-polytopes}, the corresponding marked nodes lie in
different connected components of $T_e$, so their distance is $\infty$
and the rule \eqref{eq:prokhorov-weight} assigns them weight zero. Hence
\begin{equation}\label{eq:zero-bipyramid}
	\sigma_{\mathrm P}(\mathcal{P})=0<7-f_2(\mathcal{P})=1.
\end{equation}
Thus $\mathcal{P}$ must be excluded from the desired cellwise bound.
The following proposition of Prokhorov \cite{Prokhorov1987} shows
that it is the only exception.

\begin{prop}\label{prop:facewise}
For the weights \eqref{eq:prokhorov-weight}, every three-dimensional
face satisfies
\begin{equation}\label{eq:local-target}
 \sigma_{\mathrm P}(F)\geq7-f_2(F)\quad\text{if }F\not\cong \mathcal{P},
 \qquad \sigma_{\mathrm P}(F)=0\quad\text{if }F\cong \mathcal{P}.
\end{equation}
In particular, $\mathcal{P}$ is the only combinatorial type for which
$\sigma_{\mathrm P}(F)\geq7-f_2(F)$ fails.
\end{prop}

\begin{proof}
If $f_2(F)\geq7$, the assertion follows from nonnegativity. For
$f_2(F)\leq6$, the vertex-valence restriction leaves exactly nine
types, shown in Figure~\ref{fig:small-faces}. Indeed, the dual is a
three-connected planar graph on at most six vertices, with triangular
or quadrilateral faces. There are respectively one, two, and seven
three-dimensional polytope types with four, five, and six facets.
The pentagonal pyramid is excluded by its five-valent apex. This
leaves eight good types and $\mathcal{P}$. The bound $f_2(F)\leq6$ is essential:
the valence restriction alone gives infinitely many types.

The required totals for the good types are $3,2,1$ when
$f_2(F)=4,5,6$, respectively. We recall the diagram argument that
establishes them \cite[Proposition~4 and Figure~4]{Prokhorov1987}.
Every diagram considered below is an induced subdiagram of
$\Sigma(P)$, so its Gram matrix has at most one negative eigenvalue.
Proposition~\ref{prop:prokhorov-minimal-indefinite} therefore supplies
a Lann\'er or quasi-Lann\'er subdiagram whenever such a diagram is
indefinite. For the path argument below, let
$\dist_L(a,b)$ denote the number of edges in a shortest path from
$a$ to $b$ in such a subdiagram $L$. The classification gives
\begin{equation}\label{eq:quasi-distance}
 \dist_L(a,b)\leq8,\qquad
 \text{with at most one unordered pair at distance }8.
\end{equation}
See \cite[Proposition~1 and the following discussion]{Prokhorov1987}.
This bound supplies paths for the proof. To infer a weight from such
a path, we must show that it lies in the relevant elliptic edge diagram.

Choose
a collection $M$ of marked nodes, write $H=S_F\cup M$ for the
induced diagram and put
\[
 \sigma(H)=\sum_{\substack{\{a,b\}\subset M\\\text{admissible}}}
              w(S_F\cup\{a,b\};a,b).
\]
In the following four cases all marked nodes belong to the same
indefinite connected component of $H$:
\begin{enumerate}[label=\textup{(\Roman*)},leftmargin=2.3em]
\item $|M|=4$, every indefinite induced subdiagram of $H$ contains
all four marked nodes, and every marked pair is admissible.
\item $M=\{a,b,c\}$, $H$ is indefinite, and all three pairs are admissible.
\item $M=\{a,b,c\}$, $H$ is indefinite, the pairs $\{a,b\}$ and
$\{a,c\}$ are admissible, and $S_F\cup\{b,c\}$ is singular positive
semidefinite. Geometrically the remaining pair meets only at an ideal
vertex of $F$.
\item $M=\{a,b,c,d\}$, the diagrams $S_F\cup\{a,b\}$ and
$S_F\cup\{c,d\}$ are indefinite, and the four cross-pairs between
$\{a,b\}$ and $\{c,d\}$ are admissible.
\end{enumerate}

The 2-faces represented by $M$ can be pictured as the four facets
of a tetrahedron in I; the three lateral facets of a triangular
prism in II; the base and two opposite triangular facets of a
quadrilateral pyramid in III; and the four lateral facets of a cube
in IV. In III the two triangular facets meet only at the ideal
apex; in IV the four facets form two opposite pairs.

For \eqref{eq:prokhorov-weight}, the corresponding lower bounds are
\begin{equation}\label{eq:four-local-bounds}
 \begin{array}{c|cccc}
 \text{case}&\mathrm I&\mathrm {II}&\mathrm {III}&\mathrm {IV}\\ \hline
 \sigma_{\mathrm P}(H)\geq&3&2&1&1/3.
 \end{array}
\end{equation}
In cases I--III, choose
$L\subset H$ minimal indefinite using
Proposition~\ref{prop:prokhorov-minimal-indefinite}.
It must contain all marked nodes: deleting any marked node leaves
a positive semidefinite diagram. Thus all the marks lie together
in this Lann\'er or quasi-Lann\'er diagram $L$.
The possible distance-eight pair in \eqref{eq:quasi-distance} has
only two endpoints. Since $|M|\geq3$, we can choose an initial mark
outside that pair. Join it to every other mark by a shortest path
in $L$, each of length at most seven. Split each path at every
internal marked node, as illustrated in
Figure~\ref{fig:split-marked-path}.
Since $L\subset S_F\cup M$, all internal nodes of each resulting
segment lie in $S_F$.

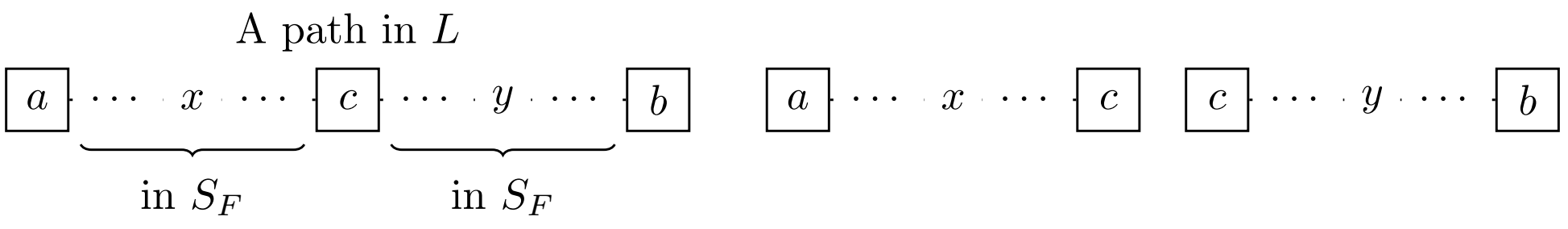
\begin{figure}[!htb]
\centering
\begin{tikzpicture}[x=1cm,y=1cm,scale=.85,
 font={\small\linespread{1}\selectfont},line width=.55pt,
 boxed mark/.style={draw,rectangle,minimum size=5.1mm,inner sep=1.3pt}]
\node[boxed mark](a)at(0,0){$a$};
\node(lx)at(.78,0){$\cdots$};\node(x)at(1.5,0){$x$};
\node(rx)at(2.22,0){$\cdots$};\node[boxed mark](c)at(3,0){$c$};
\node(ly)at(3.78,0){$\cdots$};\node(y)at(4.5,0){$y$};
\node(ry)at(5.22,0){$\cdots$};\node[boxed mark](b)at(6,0){$b$};
\draw(a)--(lx)--(x)--(rx)--(c)--(ly)--(y)--(ry)--(b);
\draw[decorate,decoration={brace,mirror,amplitude=2.5pt}](.42,-.43)--(2.58,-.43)
 node[midway,below=4.3pt]{in $S_F$};
\draw[decorate,decoration={brace,mirror,amplitude=2.5pt}](3.42,-.43)--(5.58,-.43)
 node[midway,below=4.3pt]{in $S_F$};
\node at(3,.65){A path in $L$};
\begin{scope}[xshift=7.35cm]
\node[boxed mark](aa)at(0,0){$a$};
\node(llx)at(.78,0){$\cdots$};\node(xx)at(1.5,0){$x$};
\node(rrx)at(2.22,0){$\cdots$};\node[boxed mark](cc)at(3,0){$c$};
\draw(aa)--(llx)--(xx)--(rrx)--(cc);
\node[boxed mark](ccc)at(4.05,0){$c$};
\node(lly)at(4.83,0){$\cdots$};\node(yy)at(5.55,0){$y$};
\node(rry)at(6.27,0){$\cdots$};\node[boxed mark](bb)at(7.05,0){$b$};
\draw(ccc)--(lly)--(yy)--(rry)--(bb);
\end{scope}
\end{tikzpicture}
\caption{Splitting a path in $L$ at the marked node $c$.
The internal nodes of each segment lie in $S_F$;
Coxeter labels are suppressed.}
\label{fig:split-marked-path}
\end{figure}

For each of the $|M|-1$ other marks, retain the final segment of
its chosen shortest path. The preceding mark on that path is strictly
closer to the initial mark in $L$. A pair cannot occur again with its
endpoints reversed, since this would require each endpoint to be
closer to the initial mark than the other. Thus these final segments
give $|M|-1$ distinct unordered pairs.

Figure~\ref{fig:distinct-final-segments} illustrates this choice with
$M=\{a,b,c,d\}$ and initial mark $c$. The retained segments
$c-x-a$, $a-y-b$, and $b-z-d$ give the three distinct pairs
$\{c,a\}$, $\{a,b\}$, and $\{b,d\}$.

\begin{figure}[!htb]
\centering
\begin{tikzpicture}[x=1cm,y=1cm,scale=.85,
 font={\small\linespread{1}\selectfont},line width=.55pt,
 boxed mark/.style={draw,rectangle,minimum size=5.1mm,inner sep=1.3pt},
 discarded/.style={black!42},
 retained/.style={line width=.9pt}]
\node at(-1.25,.8){Endpoint};
\node at(3,.8){Chosen shortest path from $c$};
\node at(8,.8){Retained pair};
\node at(-1.25,0){$a$};
\node[boxed mark](ca)at(0,0){$c$};\node(xa)at(1,0){$x$};
\node[boxed mark](aa)at(2,0){$a$};
\draw[retained](ca)--(xa)--(aa);
\node at(8,0){$\{c,a\}$};
\begin{scope}[yshift=-1.1cm]
\node at(-1.25,0){$b$};
\node[boxed mark,discarded](cb)at(0,0){$c$};
\node[discarded](xb)at(1,0){$x$};\node[boxed mark](ab)at(2,0){$a$};
\node(yb)at(3,0){$y$};\node[boxed mark](bb)at(4,0){$b$};
\draw[discarded](cb)--(xb)--(ab);
\draw[retained](ab)--(yb)--(bb);
\node at(8,0){$\{a,b\}$};
\end{scope}
\begin{scope}[yshift=-2.2cm]
\node at(-1.25,0){$d$};
\node[boxed mark,discarded](cd)at(0,0){$c$};
\node[discarded](xd)at(1,0){$x$};
\node[boxed mark,discarded](ad)at(2,0){$a$};\node[discarded](yd)at(3,0){$y$};
\node[boxed mark](bd)at(4,0){$b$};\node(zd)at(5,0){$z$};
\node[boxed mark](dd)at(6,0){$d$};
\draw[discarded](cd)--(xd)--(ad)--(yd)--(bd);
\draw[retained](bd)--(zd)--(dd);
\node at(8,0){$\{b,d\}$};
\end{scope}
\end{tikzpicture}
\caption{Three chosen shortest paths in $L$, starting at $c$.
Earlier portions are grey; retained final segments are black.
Boxes mark $a,b,c,d$, while $x,y,z\in S_F$.
Only the chosen paths are shown; Coxeter labels are suppressed.}
\label{fig:distinct-final-segments}
\end{figure}
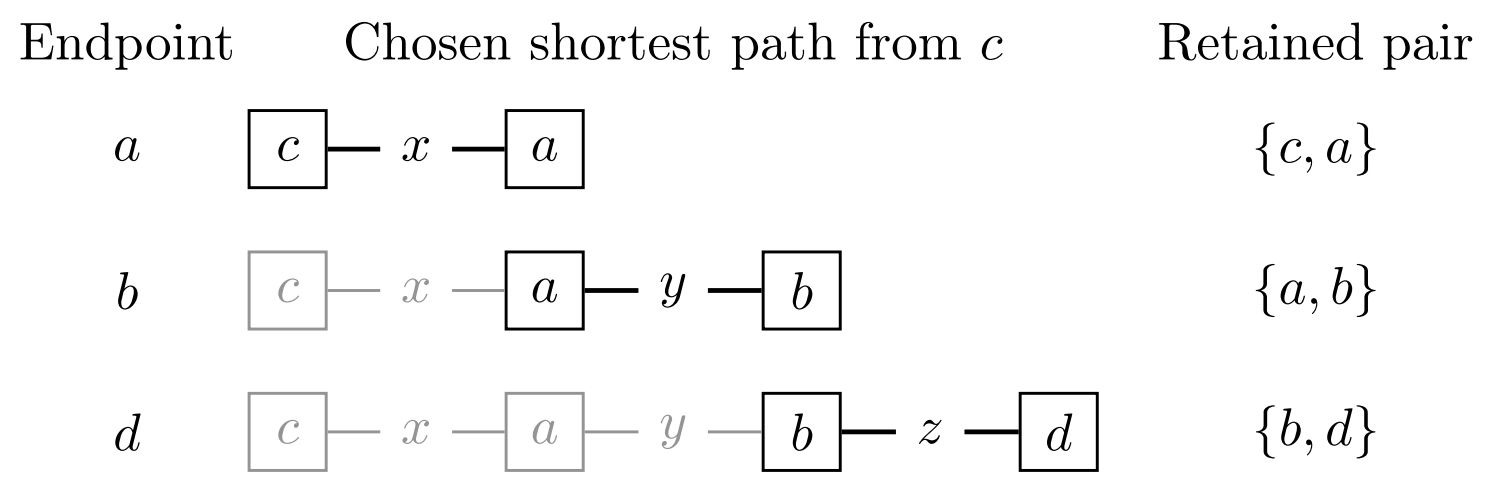

If a segment has admissible
endpoints $x,y$, it lies in the
elliptic diagram $S_F\cup\{x,y\}=T_e$, so $\dist_{T_e}(x,y)\leq7$
and the corresponding weight is one. In cases I and II all pairs
are admissible, giving at least three and two units, respectively.
In case III only $\{b,c\}$ is inadmissible, so discarding this pair
from the count of at least two leaves at least one unit.

In case IV put $X=\{a,b\}$ and $Y=\{c,d\}$. Applying
Proposition~\ref{prop:prokhorov-minimal-indefinite}, choose minimal
indefinite subdiagrams $L_X\subset S_F\cup X$ and
$L_Y\subset S_F\cup Y$. Each contains both of its own marks,
since any subdiagram omitting one of them lies in an elliptic
cross-pair diagram. Every node in either diagram is within distance
seven of at least one of its two marks. Otherwise both distances
would be eight, contradicting \eqref{eq:quasi-distance}.
A shortest path to a nearest mark cannot pass through the other
mark, which would then be closer. Thus its internal nodes lie in $S_F$.

If $L_X$ and $L_Y$ intersect, their common nodes lie in $S_F$.
Join one such node to a nearest mark in each diagram. Removing any
repetitions from the resulting walk gives a path between a mark in
$X$ and a mark in $Y$ of length at most $14$.
If the diagrams are disjoint, there must be a diagram edge joining
them. Otherwise the principal Gram matrix on their union would be
block diagonal with two indefinite blocks, giving two negative
directions in $G_H$, which has at most one. Choose endpoints
$u\in L_X$ and $v\in L_Y$ of such an edge. Joining $u$ and $v$
to nearest marks in their respective diagrams and using the single
edge $uv$ gives a path of length at most $7+1+7=15$.
The two constructions are illustrated in
Figure~\ref{fig:case-iv-short-paths}. In both cases the endpoints form
an admissible cross-pair and all internal nodes lie in $S_F$. The path therefore lies in the
corresponding elliptic diagram $T_e$, where the endpoint distance
is at most fifteen. Its weight is at least $1/3$, proving case IV.

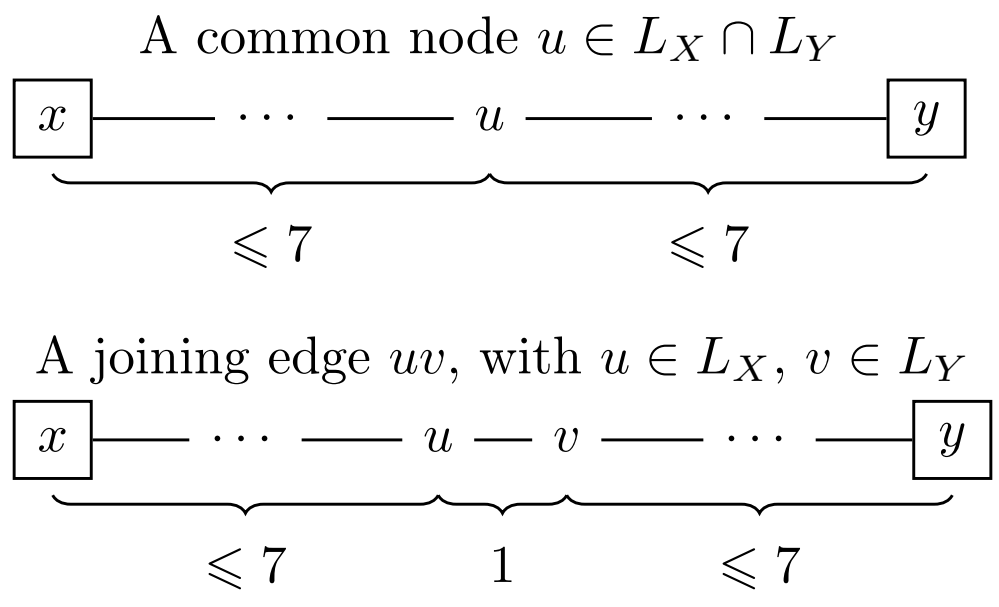
\begin{figure}[!htb]
\centering
\begin{tikzpicture}[x=1cm,y=1cm,scale=.85,
 font={\small\linespread{1}\selectfont},line width=.55pt,
 boxed mark/.style={draw,rectangle,minimum size=5.1mm,inner sep=1.3pt}]
\node[boxed mark](x)at(0,0){$x$};\node(l)at(1.7,0){$\cdots$};
\node(u)at(3.4,0){$u$};\node(r)at(5.1,0){$\cdots$};\node[boxed mark](y)at(6.8,0){$y$};
\draw(x)--(l);\draw(l)--(u)--(r);\draw(r)--(y);
\draw[decorate,decoration={brace,mirror,amplitude=3.4pt}](0,-.43)--(3.4,-.43)
 node[midway,below=6pt]{$\leq7$};
\draw[decorate,decoration={brace,mirror,amplitude=3.4pt}](3.4,-.43)--(6.8,-.43)
 node[midway,below=6pt]{$\leq7$};
\node at(3.4,.62){A common node $u\in L_X\cap L_Y$};
\begin{scope}[yshift=-2.5cm]
\node[boxed mark](xx)at(0,0){$x$};\node(ll)at(1.5,0){$\cdots$};
\node(uu)at(3,0){$u$};\node(vv)at(4,0){$v$};
\node(rr)at(5.5,0){$\cdots$};\node[boxed mark](yy)at(7,0){$y$};
\draw(xx)--(ll);\draw(ll)--(uu)--(vv)--(rr);\draw(rr)--(yy);
\draw[decorate,decoration={brace,mirror,amplitude=3.4pt}](0,-.43)--(3,-.43)
 node[midway,below=6pt]{$\leq7$};
\draw[decorate,decoration={brace,mirror,amplitude=3.4pt}](3,-.43)--(4,-.43)
 node[midway,below=6pt]{$1$};
\draw[decorate,decoration={brace,mirror,amplitude=3.4pt}](4,-.43)--(7,-.43)
 node[midway,below=6pt]{$\leq7$};
\node at(3.5,.62){A joining edge $uv$, with $u\in L_X$, $v\in L_Y$};
\end{scope}
\end{tikzpicture}
\caption{Paths in case IV, with $x\in X$ and $y\in Y$ chosen
as nearest marks. A common node gives a walk of length at most $14$;
a joining edge gives a path of length at most $15$. All internal
nodes lie in $S_F$. After removing repetitions, the path lies in
the elliptic diagram $S_F\cup\{x,y\}=T_e$. Coxeter labels are suppressed.}
\label{fig:case-iv-short-paths}
\end{figure}

For completeness, the reduction to \eqref{eq:four-local-bounds} can be
read directly from these types. The tetrahedron gives case I. The three
lateral facets of a triangular prism give case II. In a quadrilateral
pyramid, combine the base with either pair of opposite triangular
facets. This gives two case-III diagrams with disjoint admissible pairs,
hence total at least two. For the cube, choose two of its three pairs
of opposite facets. The three resulting case-IV diagrams partition its
twelve edges and give $3\cdot(1/3)=1$. In each truncated type (e) and (f),
the three facets adjacent to the new triangular facet give case II,
which provides at least two and hence exceeds the required one. The
remaining good types (g) and (h) each contain a case-III diagram.
Table~\ref{tab:small-face-contributions} summarizes these contributions,
which prove the lower bound for all eight good types.
Equation~\eqref{eq:zero-bipyramid} proves both the stated value for
$\mathcal{P}$ and the failure of the target there.\qedhere
\end{proof}

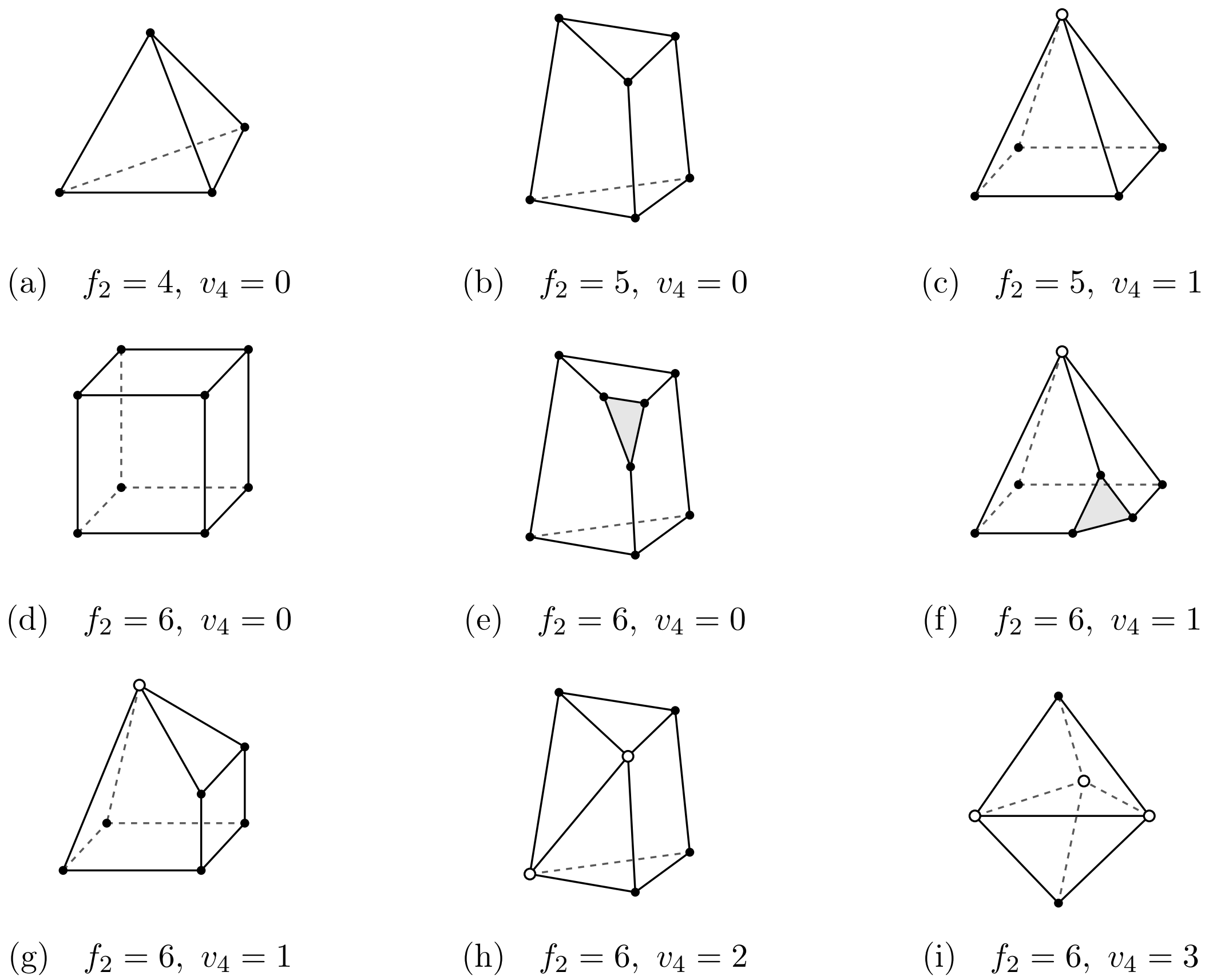
\begin{figure}[p]
\centering
\begin{minipage}[t]{.32\linewidth}\centering
\begin{tikzpicture}[scale=.84,line width=.65pt]
\path[use as bounding box] (-1.7,-1.45) rectangle (1.7,2.05);
\coordinate (A) at (0.0000,1.5000);
\coordinate (B) at (-1.2500,-0.7000);
\coordinate (C) at (0.8500,-0.7000);
\coordinate (D) at (1.3000,0.2000);
\draw (A)--(B);
\draw (A)--(C);
\draw (A)--(D);
\draw (B)--(C);
\draw[dash pattern=on 2pt off 2pt,black!65] (B)--(D);
\draw (C)--(D);
\draw[fill=black,draw=black] (A) circle (1.35pt);
\draw[fill=black,draw=black] (B) circle (1.35pt);
\draw[fill=black,draw=black] (C) circle (1.35pt);
\draw[fill=black,draw=black] (D) circle (1.35pt);
\end{tikzpicture}\par
(a)\quad $f_2=4,\ v_4=0$\end{minipage}
\hfill
\begin{minipage}[t]{.32\linewidth}\centering
\begin{tikzpicture}[scale=.84,line width=.65pt]
\path[use as bounding box] (-1.7,-1.45) rectangle (1.7,2.05);
\coordinate (A) at (-0.6500,1.7000);
\coordinate (B) at (0.9500,1.4500);
\coordinate (C) at (0.3000,0.8200);
\coordinate (D) at (-1.0500,-0.8000);
\coordinate (E) at (1.1500,-0.5000);
\coordinate (F) at (0.4000,-1.0500);
\draw (A)--(B);
\draw (A)--(C);
\draw (A)--(D);
\draw (B)--(C);
\draw (B)--(E);
\draw (C)--(F);
\draw[dash pattern=on 2pt off 2pt,black!65] (D)--(E);
\draw (D)--(F);
\draw (E)--(F);
\draw[fill=black,draw=black] (A) circle (1.35pt);
\draw[fill=black,draw=black] (B) circle (1.35pt);
\draw[fill=black,draw=black] (C) circle (1.35pt);
\draw[fill=black,draw=black] (D) circle (1.35pt);
\draw[fill=black,draw=black] (E) circle (1.35pt);
\draw[fill=black,draw=black] (F) circle (1.35pt);
\end{tikzpicture}\par
(b)\quad $f_2=5,\ v_4=0$\end{minipage}
\hfill
\begin{minipage}[t]{.32\linewidth}\centering
\begin{tikzpicture}[scale=.84,line width=.65pt]
\path[use as bounding box] (-1.7,-1.45) rectangle (1.7,2.05);
\coordinate (A) at (0.0000,1.7500);
\coordinate (B) at (-1.2000,-0.7500);
\coordinate (E) at (0.7800,-0.7500);
\coordinate (D) at (1.3800,-0.0800);
\coordinate (C) at (-0.6000,-0.0800);
\draw (A)--(B);
\draw[dash pattern=on 2pt off 2pt,black!65] (A)--(C);
\draw (A)--(D);
\draw (A)--(E);
\draw[dash pattern=on 2pt off 2pt,black!65] (B)--(C);
\draw (B)--(E);
\draw[dash pattern=on 2pt off 2pt,black!65] (C)--(D);
\draw (D)--(E);
\draw[fill=white,draw=black] (A) circle (2.1pt);
\draw[fill=black,draw=black] (B) circle (1.35pt);
\draw[fill=black,draw=black] (E) circle (1.35pt);
\draw[fill=black,draw=black] (D) circle (1.35pt);
\draw[fill=black,draw=black] (C) circle (1.35pt);
\end{tikzpicture}\par
(c)\quad $f_2=5,\ v_4=1$\end{minipage}
\par\vspace{3mm}
\begin{minipage}[t]{.32\linewidth}\centering
\begin{tikzpicture}[scale=.84,line width=.65pt]
\path[use as bounding box] (-1.7,-1.45) rectangle (1.7,2.05);
\coordinate (A) at (-1.0000,1.1500);
\coordinate (B) at (0.7500,1.1500);
\coordinate (C) at (0.7500,-0.7500);
\coordinate (D) at (-1.0000,-0.7500);
\coordinate (E) at (-0.4000,1.7800);
\coordinate (F) at (1.3500,1.7800);
\coordinate (G) at (1.3500,-0.1200);
\coordinate (H) at (-0.4000,-0.1200);
\draw (A)--(B);
\draw (A)--(D);
\draw (A)--(E);
\draw (B)--(C);
\draw (B)--(F);
\draw (C)--(D);
\draw (C)--(G);
\draw[dash pattern=on 2pt off 2pt,black!65] (D)--(H);
\draw (E)--(F);
\draw[dash pattern=on 2pt off 2pt,black!65] (E)--(H);
\draw (F)--(G);
\draw[dash pattern=on 2pt off 2pt,black!65] (G)--(H);
\draw[fill=black,draw=black] (A) circle (1.35pt);
\draw[fill=black,draw=black] (B) circle (1.35pt);
\draw[fill=black,draw=black] (C) circle (1.35pt);
\draw[fill=black,draw=black] (D) circle (1.35pt);
\draw[fill=black,draw=black] (E) circle (1.35pt);
\draw[fill=black,draw=black] (F) circle (1.35pt);
\draw[fill=black,draw=black] (G) circle (1.35pt);
\draw[fill=black,draw=black] (H) circle (1.35pt);
\end{tikzpicture}\par
(d)\quad $f_2=6,\ v_4=0$\end{minipage}
\hfill
\begin{minipage}[t]{.32\linewidth}\centering
\begin{tikzpicture}[scale=.84,line width=.65pt]
\path[use as bounding box] (-1.7,-1.45) rectangle (1.7,2.05);
\coordinate (A) at (-0.6500,1.7000);
\coordinate (B) at (0.9500,1.4500);
\coordinate (D) at (-1.0500,-0.8000);
\coordinate (E) at (1.1500,-0.5000);
\coordinate (F) at (0.4000,-1.0500);
\coordinate (CA) at (-0.0325,1.1280);
\coordinate (CB) at (0.5275,1.0405);
\coordinate (CF) at (0.3350,0.1655);
\fill[black!10] (CA)--(CB)--(CF)--cycle;
\draw (A)--(B);
\draw (A)--(CA);
\draw (A)--(D);
\draw (B)--(CB);
\draw (B)--(E);
\draw (CA)--(CB);
\draw (CA)--(CF);
\draw (CB)--(CF);
\draw (CF)--(F);
\draw[dash pattern=on 2pt off 2pt,black!65] (D)--(E);
\draw (D)--(F);
\draw (E)--(F);
\draw[fill=black,draw=black] (A) circle (1.35pt);
\draw[fill=black,draw=black] (B) circle (1.35pt);
\draw[fill=black,draw=black] (D) circle (1.35pt);
\draw[fill=black,draw=black] (E) circle (1.35pt);
\draw[fill=black,draw=black] (F) circle (1.35pt);
\draw[fill=black,draw=black] (CA) circle (1.35pt);
\draw[fill=black,draw=black] (CB) circle (1.35pt);
\draw[fill=black,draw=black] (CF) circle (1.35pt);
\end{tikzpicture}\par
(e)\quad $f_2=6,\ v_4=0$\end{minipage}
\hfill
\begin{minipage}[t]{.32\linewidth}\centering
\begin{tikzpicture}[scale=.84,line width=.65pt]
\path[use as bounding box] (-1.7,-1.45) rectangle (1.7,2.05);
\coordinate (A) at (0.0000,1.7500);
\coordinate (B) at (-1.2000,-0.7500);
\coordinate (D) at (1.3800,-0.0800);
\coordinate (C) at (-0.6000,-0.0800);
\coordinate (EA) at (0.5304,0.0500);
\coordinate (EB) at (0.1464,-0.7500);
\coordinate (ED) at (0.9720,-0.5356);
\fill[black!10] (EA)--(EB)--(ED)--cycle;
\draw (A)--(B);
\draw[dash pattern=on 2pt off 2pt,black!65] (A)--(C);
\draw (A)--(D);
\draw (A)--(EA);
\draw[dash pattern=on 2pt off 2pt,black!65] (B)--(C);
\draw (B)--(EB);
\draw[dash pattern=on 2pt off 2pt,black!65] (C)--(D);
\draw (D)--(ED);
\draw (EA)--(EB);
\draw (EA)--(ED);
\draw (EB)--(ED);
\draw[fill=white,draw=black] (A) circle (2.1pt);
\draw[fill=black,draw=black] (B) circle (1.35pt);
\draw[fill=black,draw=black] (D) circle (1.35pt);
\draw[fill=black,draw=black] (C) circle (1.35pt);
\draw[fill=black,draw=black] (EA) circle (1.35pt);
\draw[fill=black,draw=black] (EB) circle (1.35pt);
\draw[fill=black,draw=black] (ED) circle (1.35pt);
\end{tikzpicture}\par
(f)\quad $f_2=6,\ v_4=1$\end{minipage}
\par\vspace{3mm}
\begin{minipage}[t]{.32\linewidth}\centering
\begin{tikzpicture}[scale=.84,line width=.65pt]
\path[use as bounding box] (-1.7,-1.45) rectangle (1.7,2.05);
\coordinate (A) at (-0.1500,1.8000);
\coordinate (B) at (-1.2000,-0.7500);
\coordinate (C) at (-0.6000,-0.1000);
\coordinate (D) at (1.3000,-0.1000);
\coordinate (E) at (0.7000,-0.7500);
\coordinate (F) at (0.7000,0.3000);
\coordinate (G) at (1.3000,0.9500);
\draw (A)--(B);
\draw[dash pattern=on 2pt off 2pt,black!65] (A)--(C);
\draw (A)--(F);
\draw (A)--(G);
\draw[dash pattern=on 2pt off 2pt,black!65] (B)--(C);
\draw (B)--(E);
\draw[dash pattern=on 2pt off 2pt,black!65] (C)--(D);
\draw (D)--(E);
\draw (D)--(G);
\draw (E)--(F);
\draw (F)--(G);
\draw[fill=white,draw=black] (A) circle (2.1pt);
\draw[fill=black,draw=black] (B) circle (1.35pt);
\draw[fill=black,draw=black] (C) circle (1.35pt);
\draw[fill=black,draw=black] (D) circle (1.35pt);
\draw[fill=black,draw=black] (E) circle (1.35pt);
\draw[fill=black,draw=black] (F) circle (1.35pt);
\draw[fill=black,draw=black] (G) circle (1.35pt);
\end{tikzpicture}\par
(g)\quad $f_2=6,\ v_4=1$\end{minipage}
\hfill
\begin{minipage}[t]{.32\linewidth}\centering
\begin{tikzpicture}[scale=.84,line width=.65pt]
\path[use as bounding box] (-1.7,-1.45) rectangle (1.7,2.05);
\coordinate (A) at (-0.6500,1.7000);
\coordinate (B) at (0.9500,1.4500);
\coordinate (C) at (0.3000,0.8200);
\coordinate (D) at (-1.0500,-0.8000);
\coordinate (E) at (1.1500,-0.5000);
\coordinate (F) at (0.4000,-1.0500);
\draw (A)--(B);
\draw (A)--(C);
\draw (A)--(D);
\draw (B)--(C);
\draw (B)--(E);
\draw (C)--(D);
\draw (C)--(F);
\draw[dash pattern=on 2pt off 2pt,black!65] (D)--(E);
\draw (D)--(F);
\draw (E)--(F);
\draw[fill=black,draw=black] (A) circle (1.35pt);
\draw[fill=black,draw=black] (B) circle (1.35pt);
\draw[fill=white,draw=black] (C) circle (2.1pt);
\draw[fill=white,draw=black] (D) circle (2.1pt);
\draw[fill=black,draw=black] (E) circle (1.35pt);
\draw[fill=black,draw=black] (F) circle (1.35pt);
\end{tikzpicture}\par
(h)\quad $f_2=6,\ v_4=2$\end{minipage}
\hfill
\begin{minipage}[t]{.32\linewidth}\centering
\begin{tikzpicture}[scale=.84,line width=.65pt]
\path[use as bounding box] (-1.7,-1.45) rectangle (1.7,2.05);
\coordinate (A) at (-0.0500,1.6500);
\coordinate (B) at (-1.2000,0.0000);
\coordinate (C) at (1.2000,0.0000);
\coordinate (D) at (0.3000,0.4800);
\coordinate (E) at (-0.0500,-1.2000);
\draw (A)--(B);
\draw (A)--(C);
\draw[dash pattern=on 2pt off 2pt,black!65] (A)--(D);
\draw (B)--(C);
\draw[dash pattern=on 2pt off 2pt,black!65] (B)--(D);
\draw (B)--(E);
\draw[dash pattern=on 2pt off 2pt,black!65] (C)--(D);
\draw (C)--(E);
\draw[dash pattern=on 2pt off 2pt,black!65] (D)--(E);
\draw[fill=black,draw=black] (A) circle (1.35pt);
\draw[fill=white,draw=black] (B) circle (2.1pt);
\draw[fill=white,draw=black] (C) circle (2.1pt);
\draw[fill=white,draw=black] (D) circle (2.1pt);
\draw[fill=black,draw=black] (E) circle (1.35pt);
\end{tikzpicture}\par
(i)\quad $f_2=6,\ v_4=3$\end{minipage}
\caption{The eight good $3$-face types (a)--(h) and the bad
polyhedron $\mathcal{P}$ in (i). Here $v_4=v_4(F)$ counts four-valent
vertices. Panels (a)--(h) redraw Prokhorov's Figure~4; panel (i) is
also shown in his Figure~1. Open vertices have valence four and are ideal.
Filled vertices have valence three and need not all be finite.}
\label{fig:small-faces}
\end{figure}

\begin{table}[p]
\centering
\caption{Local contributions for the good $3$-faces in
Figure~\ref{fig:small-faces}. Multiple configurations have disjoint
scored edge sets, so their contributions add.}
\label{tab:small-face-contributions}
\begingroup
\small
\linespread{1}\selectfont
\begin{tabular}{c|l|c|c|c}
\toprule[1pt]
Panels & Combinatorial type & \shortstack{Local\\configurations}
 & \shortstack{Guaranteed\\contribution} & $7-f_2(F)$\\\hline
(a) & Tetrahedron & I & $3$ & $3$\\\hline
(b) & Triangular prism & II & $2$ & $2$\\\hline
(c) & Quadrilateral pyramid & III (twice) & $2\cdot1$ & $2$\\\hline
(d) & Cube & IV (three times) & $3\cdot(1/3)$ & $1$\\\hline
(e), (f) & Six-facet types & II & $2$ & $1$\\\hline
(g), (h) & Six-facet types & III & $1$ & $1$\\
\bottomrule[1pt]
\end{tabular}
\endgroup
\end{table}

\subsection{Double counting and Prokhorov's bound}\label{subsec:weights}
Suppose more generally that a nonnegative rule $w$ satisfies
$\sigma(F)\geq7-f_2(F)$ for every good 3-face and
\begin{equation}\label{eq:edge-target}
 W(T_e)\leq C(n-1)\qquad\text{for every edge }e,
\end{equation}
where $C>0$ is independent of $n$ and depends on the chosen weight rule.
For Prokhorov's rule, \eqref{eq:old-edge-constant} gives $C=C_{\mathrm P}$.
Each bad 3-face has six facets, so
replacing the target $7-f_2(F)=1$ by nonnegativity loses one unit.
Summing over the 3-faces, and then over the edges, gives
\begin{equation}\label{eq:global-comparison}
 7-a_3^2-b\leq a_w\leq C(n-1)\frac{a_1}{a_3}
 =\frac{2C}{n-2}a_3^1.
\end{equation}
The lower bound on the left is not asserted for every 3-face separately.
It is the average after accounting for the exceptional 3-faces.

Prokhorov controlled their proportion by
\begin{equation}\label{eq:bad-proportion}
 b\leq
 \begin{cases}
 \displaystyle\frac34+\frac{15}{4(n-3)},&n\text{ even},\\[5pt]
 \displaystyle\frac34+\frac{15}{4(n-2)},&n\text{ odd},
 \end{cases}\qquad n\geq8.
\end{equation}
He first proved that a four-dimensional face with $k$ facets has at
least $10-k$ good facets, and then applied
\eqref{eq:nikulin-khovanskii} with $(i,k)=(3,4)$
\cite[Lemma~1]{Prokhorov1987}. The upper bound tends to $3/4$.
The proportion $b$ itself need not do so.
The explicit bound for $a_3^2$ in \eqref{eq:A-E-bounds} controls
the left side of \eqref{eq:global-comparison}, together with
\eqref{eq:bad-proportion}. The bound for $a_3^1$ controls its right
side. Substitution gives
\begin{equation}\label{eq:old-comparison}
\begin{aligned}
 \frac14-\frac{12}{n-4}-\frac{15}{4(n-3)}
 &<a_w<\frac{24C(n-1)}{(n-2)(n-4)}&& (n\text{ even}),\\[4pt]
 \frac14-\frac{12}{n-3}-\frac{15}{4(n-2)}
 &<a_w<\frac{24Cn}{(n-2)(n-3)}&& (n\text{ odd}).
\end{aligned}
\end{equation}
This explains the choice of $7-f_2(F)$: the average lower bound tends
to $7-6-3/4=1/4$, whereas the edgewise upper bound is
$24C/n+O(n^{-2})$ and tends to zero. Thus the two bounds become
incompatible for large $n$. The exact comparison gives
\begin{equation}\label{eq:prokhorov-lemma}
 n<96C+68
\end{equation}
by \cite[Lemma~2]{Prokhorov1987}. Applying
Proposition~\ref{prop:facewise} and \eqref{eq:old-edge-constant} yields
$n<96C_{\mathrm P}+68=996$.

\subsection{The estimates to be sharpened}\label{subsec:motivation-162}
We retain the cellwise lower bound and improve both upper estimates
in \eqref{eq:global-comparison}. Section~\ref{ext:sec:weights}
reduces the leading edgewise coefficient from $29/3$ to $73/12$.
The larger improvement comes from the exceptional 3-faces: their
square cusp sections force an edgewise deficit, giving
$a_w\leq\lambda_n(a_3^1-9b)$ in Section~\ref{ext:sec:cusp}.
Section~\ref{ext:sec:h} uses the same square sections and Timorin's
truncation identities to sharpen the bound for $a_3^2$.
After eliminating the average face numbers, the comparison is
affine in $b$ and fails at both $b=0$ and $b=1$ for $n\geq163$;
Section~\ref{ext:sec:conclusion} therefore gives $n\leq162$.
Prokhorov's independent estimate \eqref{eq:bad-proportion} is needed
only for the comparison with his original argument.

\FloatBarrier

\section{Refined weights and local estimates}\label{ext:sec:weights}

We reduce the weight of certain marked pairs while preserving
Prokhorov's four local lower bounds $3,2,1,1/3$. The resulting estimates are
\[
W(T_e)\leq\frac{73}{12}(n-1)-\frac{23}{4},
\qquad
\sigma(F)\geq7-f_2(F)\quad(F\text{ good}),
\]
as proved in Lemma~\ref{ext:lem:budget} and
Proposition~\ref{ext:prop:face}.

\subsection{Two-sided extension}
Every connected elliptic diagram with more than eight nodes has
type $A_r,B_r$, or $D_r$. Away from its endpoints, branches, and
heavy edges, it consists of ordinary chains. For each fixed
distance $k$, these contain $r+O(1)$ marked pairs, whose weights
determine the leading coefficient of $W(T)$.

To reduce that coefficient, we distinguish a path that continues
far in both directions from one near an endpoint, branch, or heavy
edge. The extension is measured in the full elliptic edge diagram
$T_e$ and may leave a chosen minimal indefinite subdiagram.
The local estimates below verify which reductions are permissible.

\subsection{The rule and its effect on long chains}
\label{ext:subsec:path-rule}
We fix the path terminology used in both parts. A path in a tree
component of a Coxeter diagram $T$ is \emph{special} if it contains
a node of degree at least three in $T$, including a path endpoint,
or a heavy edge. As in Subsection~\ref{subsec:hyperbolic-polytopes},
a bold edge is treated as heavy. A path that is not special is
called \emph{nonspecial}.

A marked pair $(a,b)$ in a tree component has a
\emph{two-sided extension of length $s$} if its path extends by
$s$ edges at each end, every edge of the extended path has label
three, and every internal node has degree two in $T$. Only the
two new outer endpoints are exempt from the degree condition.
In particular, the original marks $a,b$ must have degree two.

For the elliptic diagrams used in this part, put
$k=\dist_T(a,b)$ and define
\begin{equation}\label{ext:eq:preliminary-weight}
w_{\mathrm{pre}}(T;a,b)=
\begin{cases}
1,&1\leq k\leq6,\text{ or }k=7\text{ and the path is special},\\
1/2,&k=7\text{ and the path is nonspecial},\\
1/3,&8\leq k\leq12,\text{ or }k=13\text{ and the path is special},\\
1/6,&k=13\text{ and the path is nonspecial},\\
1/6,&k=14\text{ and the path is special},\\
1/12,&k=14\text{ and the path is nonspecial},\\
0,&\text{otherwise}.
\end{cases}
\end{equation}
The refined weight is
\begin{equation}\label{ext:eq:weight}
w(T;a,b)=
\begin{cases}
1/2,&3\leq k\leq5\text{ and an extension of length }9-k\text{ exists},\\
1/3,&k\in\{6,7\}\text{ and an extension of length }9-k\text{ exists},\\
w_{\mathrm{pre}}(T;a,b),&\text{otherwise}.
\end{cases}
\end{equation}
Here and below, ``extension'' in a weight rule means a two-sided
extension in the full diagram $T$.
For any elliptic diagram $T$, extend the edge-sum notation by setting
\[
W(T)=\sum_{\substack{\{a,b\}\subset V(T)\\a\ne b}}w(T;a,b).
\]
This total is additive over connected components. For an incidence
$e\subset F$, use the marked pair from \eqref{eq:marked-bijection}
and set $w(F,e)=w(T_e;a,b)$. The sum $W(T_e)$ then agrees with the
edgewise sum in \eqref{eq:weight-sums}, and $a_w$ retains its definition
there.

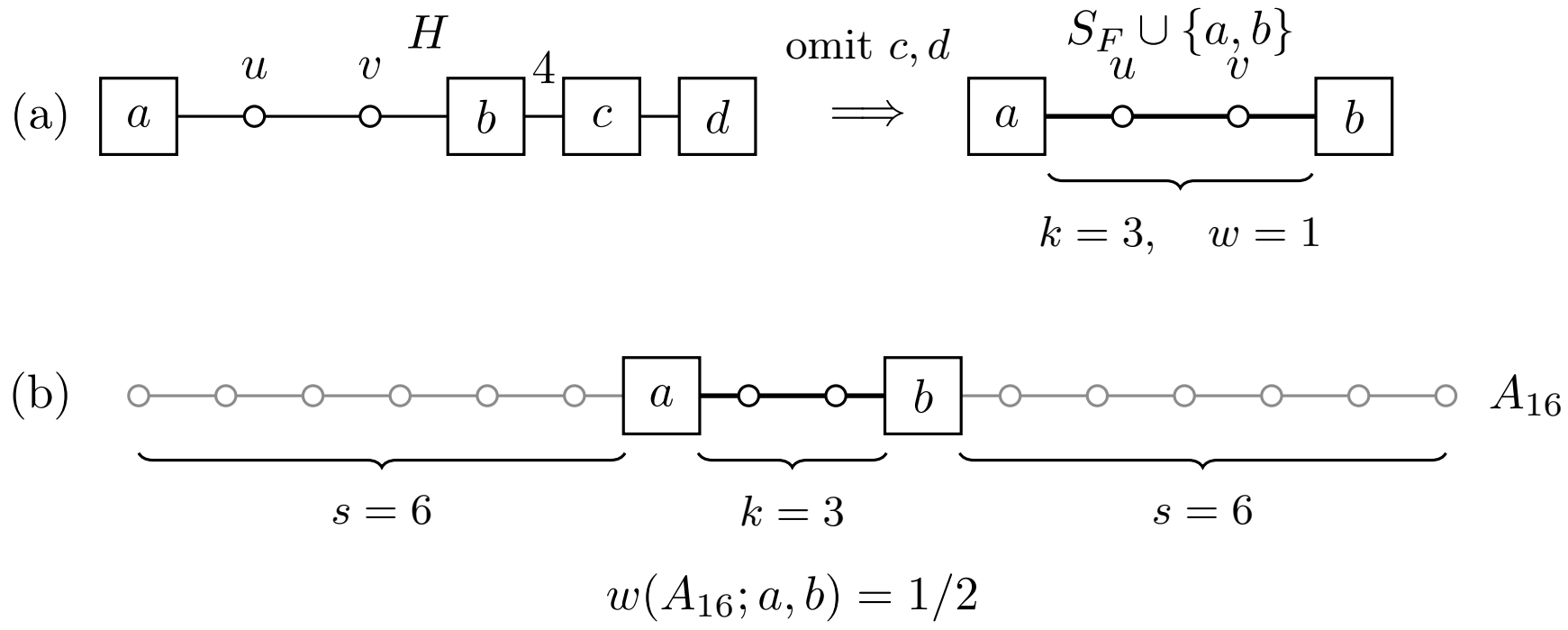
\begin{figure}[!htb]
\centering
\begin{tikzpicture}[x=.85cm,y=1cm,line width=.6pt,font=\small,
 boxed mark/.style={draw,rectangle,fill=white,minimum size=5.6mm,inner sep=1.3pt}]
\node at (-.85,0){(a)};
\node at (2.5,.62){$H$};
\draw (0,0)--(3,0);
\draw (3,0)--node[midway,above=3pt]{$4$}(4,0);
\draw (4,0)--(5,0);
\node[boxed mark] at (0,0){$a$};
\draw[fill=white](1,0)circle(2.1pt);\node[above=4pt] at(1,0){$u$};
\draw[fill=white](2,0)circle(2.1pt);\node[above=4pt] at(2,0){$v$};
\node[boxed mark] at (3,0){$b$};
\node[boxed mark] at (4,0){$c$};
\node[boxed mark] at (5,0){$d$};
\node at (6.3,0){$\Longrightarrow$};
\node[font=\footnotesize] at (6.3,.50){omit $c,d$};
\begin{scope}[xshift=6.375cm]
\node at (1.5,.62){$S_F\cup\{a,b\}$};
\draw[line width=1pt](0,0)--(3,0);
\node[boxed mark] at (0,0){$a$};
\draw[fill=white](1,0)circle(2.1pt);\node[above=4pt] at(1,0){$u$};
\draw[fill=white](2,0)circle(2.1pt);\node[above=4pt] at(2,0){$v$};
\node[boxed mark] at (3,0){$b$};
\draw[decorate,decoration={brace,mirror,amplitude=3pt}](.36,-.42)--(2.64,-.42)
 node[midway,below=5pt,font=\footnotesize]{$k=3,\quad w=1$};
\end{scope}
\begin{scope}[x=.64cm,yshift=-2.05cm]
\node at (-1.13,0){(b)};
\draw[black!45](0,0)--(6,0);\draw[black!45](9,0)--(15,0);
\draw[line width=1pt](6,0)--(9,0);
\foreach\i in {0,...,5}{\draw[black!45,fill=white](\i,0)circle(2.1pt);}
\foreach\i in {10,...,15}{\draw[black!45,fill=white](\i,0)circle(2.1pt);}
\foreach\i in {7,8}{\draw[fill=white](\i,0)circle(2.1pt);}
\node[boxed mark] at(6,0){$a$};\node[boxed mark] at(9,0){$b$};
\node[right=5pt] at(15,0){$A_{16}$};
\draw[decorate,decoration={brace,mirror,amplitude=3pt}](0,-.42)--(5.58,-.42)
 node[midway,below=5pt,font=\footnotesize]{$s=6$};
\draw[decorate,decoration={brace,mirror,amplitude=3pt}](6.42,-.42)--(8.58,-.42)
 node[midway,below=5pt,font=\footnotesize]{$k=3$};
\draw[decorate,decoration={brace,mirror,amplitude=3pt}](9.42,-.42)--(15,-.42)
 node[midway,below=5pt,font=\footnotesize]{$s=6$};
\node at(7.5,-1.48){$w(A_{16};a,b)=1/2$};
\end{scope}
\end{tikzpicture}
\caption{Distance-three pairs without (a) and with (b) a six-edge
extension at each end. Boxes mark selected nodes; unlabeled edges
have Coxeter label three.}
\label{ext:fig:extension}
\end{figure}

Figure~\ref{ext:fig:extension} compares two pairs at distance three:
one whose unit weight is needed for a cellwise lower bound, and one
to which the extension rule assigns weight $1/2$.

In Figure~\ref{ext:fig:extension}(a), $H$ is the diagram of a
finite-volume hyperbolic Coxeter $5$-simplex. The nodes $u,v$ define
a three-dimensional face $F$, which is a tetrahedron; thus
$S_F=\{u,v\}$ and its four facets correspond to the marks
$a,b,c,d$. For each scored pair, retain $S_F$ and that pair,
deleting the other two marks. The resulting weights are:
\begin{center}
\begingroup
\small
\linespread{1}\selectfont
\begin{tabular}{l|l|c}
\toprule[1pt]
Marked pair & In its scored diagram & Weight\\\hline
$\{a,b\}$ & The path $a-u-v-b$, of distance three & $1$\\\hline
$\{b,c\},\ \{c,d\}$ & Distance one & $1$ each\\\hline
$\{a,c\},\ \{a,d\},\ \{b,d\}$ & Disconnected & $0$\\
\bottomrule[1pt]
\end{tabular}
\endgroup
\end{center}
The diagram for $\{a,b\}$ is $A_4$, so neither end admits the
required six-edge extension. Consequently
\[
\sigma(F)=1+1+1=3=7-f_2(F).
\]
Lowering every nonspecial distance-three pair to $1/2$, without the
extension condition, would instead give $\sigma(F)=5/2<3$.
Thus distance alone does not identify the pairs that can be lowered.

Part~(b) illustrates the extension rule in the elliptic diagram
$A_{16}$. Number its nodes from $0$ to $15$, with $a=6$ and $b=9$.
Their path has length three and extends by six edges at each
end, so its weight is $1/2$. The pair $\{0,3\}$ has the same
distance but no extension on the left, and retains weight $1$.
The extension condition therefore distinguishes these pairs while
preserving the necessary unit contribution in part~(a). The local
lower bounds for the complete rule are proved below.

The cutoffs in \eqref{ext:eq:weight} are part of an explicit feasible choice
of weights. We make no assertion that this choice is optimal among all
possible rules. One can regard the design problem as minimizing the
coefficient of $|T|$ in $W(T)$ subject to the four local lower bounds. The
preliminary rule $w_{\mathrm{pre}}$ first reduces contributions at distances seven and
thirteen through fifteen, distinguishing paths that contain a branch or
a larger label. Its values at distances thirteen and fourteen retain a total of $1/3$
in four long configurations of case IV. The available pairs and
their assigned weights give
\begin{center}
\begingroup
\small
\linespread{1}\selectfont
\begin{tabular}{l|c|c|c}
\toprule[1pt]
Path type & Distance & Number of pairs & Total weight\\\hline
Special & $13$ & $1$ & $1\cdot\frac13=\frac13$\\\hline
Special & $14$ & $2$ & $2\cdot\frac16=\frac13$\\\hline
Nonspecial & $13$ & $2$ & $2\cdot\frac16=\frac13$\\\hline
Nonspecial & $14$ & $4$ & $4\cdot\frac1{12}=\frac13$\\
\bottomrule[1pt]
\end{tabular}
\endgroup
\end{center}
The configurations imposing these requirements are checked in
Subsection~\ref{ext:app:long}. Thus the weights $1/6$ and $1/12$
reflect the number of available pairs in these cases.

The extension condition permits further reductions at distances three
through seven. Its length $9-k$ is constrained by the marked deletion
conditions, not just by the diameter of a minimal indefinite diagram.
The type-III calculation following Table~\ref{ext:tab:marking-counts}
shows why replacing it by $8-k$ would lose the local lower bound.
Lemma~\ref{ext:lem:local} verifies the length $9-k$ for all configurations.

Number the nodes of $A_r$ from $0$ to $r-1$. A pair $\{i,i+k\}$
has a two-sided extension of length $s$ precisely when
$s\leq i\leq r-1-k-s$. Thus exactly $\max\{r-k-2s,0\}$ pairs
admit it, and at most $2s$ fail. For $s=9-k$ and $3\leq k\leq7$,
there are at most $12+10+8+6+4=40$ exceptions altogether.
They affect only the constant term in $W(A_r)$.

Far from the ends of a long chain, the weights at each fixed distance are
\begingroup
\small
\linespread{1}\selectfont
\[
\begin{array}{c|c|c|c|c|c|c}
\toprule[1pt]
k&1,2&3,4,5&6,7&8,9,10,11,12&13&14\\\hline
w&1&1/2&1/3&1/3&1/6&1/12\\
\bottomrule[1pt]
\end{array}
\]
\endgroup
At each fixed distance $k$, the count $r-k$ has leading coefficient
one in $r$, and the bounded end corrections do not change this
coefficient. Thus the coefficient for such a chain is
\begin{equation}\label{ext:eq:chain-coefficient}
2+\frac32+\frac73+\frac16+\frac1{12}=\frac{73}{12}.
\end{equation}
For comparison, the preliminary rule has coefficient $101/12$. The reductions enabled by extensions subtract
\[
3\cdot\frac12+\frac23+\frac16=\frac73,
\]
giving $101/12-7/3=73/12$. This explains how the proposed modification
improves the leading constant. The next two steps establish that it also
preserves the lower bounds and controls the constant terms uniformly over
\emph{all} elliptic diagrams. For this fixed rule the coefficient $73/12$
cannot be improved uniformly, since $W(A_r)/r\to73/12$.

\subsection{Local lower bounds}\label{ext:sec:local-algorithm-input}
We use the four marked situations I--IV introduced in the proof of
Proposition~\ref{prop:facewise}. In case I the scored diagrams are the six
double marked deletions. In cases II and III they are, respectively, all
three pairs and the two pairs containing the distinguished mark $a$. In
case IV they are the four cross-pairs. In case I every single marked
deletion is positive semidefinite and every double marked deletion is
positive definite. In case II every single marked deletion is positive
definite. In case III deletion of $a$ is positive semidefinite, while
deletions of $b,c$ are positive definite. Let $\mathcal C$ denote the
complete list of these tests in the applicable case.

Two elementary matrix facts justify the finite reduction and the
hand calculation below.

\begin{lem}\label{ext:lem:ordinary-arms}
Suppose $\Sigma$ consists of paths of lengths
$\ell_1,\ldots,\ell_p\geq1$ meeting only at one common endpoint,
with every edge labeled three and no additional edges. Set
\[
\eta=2-\sum_{i=1}^p\frac{\ell_i}{\ell_i+1}.
\]
Then $G(\Sigma)$ is positive definite for $\eta>0$, positive
semidefinite of corank one for $\eta=0$, and nonsingular with
exactly one negative eigenvalue for $\eta<0$.
\end{lem}
\begin{proof}
In $2G(\Sigma)$, the block on an arm of length $\ell$ is the
positive definite tridiagonal matrix $K_\ell$ with diagonal two
and adjacent entries minus one. Its endpoint inverse entry is
$(K_\ell^{-1})_{11}=\ell/(\ell+1)$. Eliminating the arm blocks
therefore leaves the scalar Schur complement $\eta$ at the common
endpoint. Sylvester's law of inertia proves the assertions.
\end{proof}

\begin{lem}\label{ext:lem:label-reduction}
Let $G=I-A$, where $A$ is symmetric, entrywise nonnegative, and has
zero diagonal. If $\overline A$ is symmetric and
$0\leq\overline A\leq A$ entrywise, each principal
submatrix of $I-\overline A$ corresponding to a positive definite
or positive semidefinite principal submatrix of $G$ has the same
respective property. In particular, lowering Coxeter labels or
deleting edges preserves these conditions. Lowering labels
greater than four to four also preserves every weight
\eqref{ext:eq:weight}.
\end{lem}
\begin{proof}
For a retained node set $U$, monotonicity of the spectral radius gives
$\rho(\overline A_U)\leq\rho(A_U)$; the required condition is
$\rho(A_U)<1$ or $\rho(A_U)\leq1$, respectively.
Lowering a label decreases $\cos(\pi/m)$, and deleting an edge
sets it to zero. Lowering labels greater than four to four leaves the
underlying graph, special paths, and ordinary extensions unchanged,
so it also leaves the weights unchanged.
\end{proof}

\begin{lem}\label{ext:lem:local}
The scores of configurations \textup{I--IV} are at least
$3,2,1,1/3$, respectively.
\end{lem}
\begin{proof}
For types I--III, choose a minimal indefinite principal subdiagram $L$ of
$H$ using Proposition~\ref{prop:prokhorov-minimal-indefinite}.
If a marked node were absent from $L$, the corresponding positive
semidefinite marked deletion would contain the indefinite principal
subdiagram $L$, a contradiction. Thus $M\subset V(L)$. Since
$|M|\geq3$, the diagram $L$ has order at least three, and the
classification of Lann\'er and quasi-Lann\'er diagrams listed in
Appendix~\ref{app:lanner} gives $|L|\leq10$.

Replace every label at least four in $L$ by four, obtaining $\overline L$ on the
same nodes, and make the same replacements in the full diagram $H$,
obtaining $\overline H$. All marks remain in $\overline L$. By
Lemma~\ref{ext:lem:label-reduction}, every prescribed marked deletion
retains its positive definite or semidefinite condition. For every
scored pair, the full diagrams $H$ and $\overline H$ also give identical
weights: their paths, branches, heavy edges, and ordinary extensions
are unchanged. It therefore suffices to consider full extensions of
$\overline L$ satisfying these necessary deletion conditions; $\overline L$ itself need
not remain indefinite.

For a fixed representative $\overline L$ and marking $M$, the quantity to control
is the total weight of the scored pairs in every compatible full
extension $J$ satisfying these deletion conditions. Let $\mathcal E(M)$ be the set of scored
pairs and put $D_{xy}=M\setminus\{x,y\}$. For each
$\{x,y\}\in\mathcal E(M)$, we calculate a uniform lower bound $\omega_{xy}$ satisfying
\[
\omega_{xy}\leq w(J\setminus D_{xy};x,y)
\]
for every such full extension $J$. The weight in the core alone cannot
replace this quantity, since nodes outside $\overline L$ may provide an extension. The required calculation is then
\[
\sum_{\{x,y\}\in\mathcal E(M)}\omega_{xy}\geq
\begin{cases}
3,&\text{in type I},\\
2,&\text{in type II},\\
1,&\text{in type III}.
\end{cases}
\]
To calculate $\omega_{xy}$, first inspect the pair in
$\overline L\setminus D_{xy}$. A disconnected pair receives the bound zero.
Otherwise its path remains the unique path in every actual elliptic
extension, and an existing branch or heavy edge cannot disappear.
Special paths and paths outside the extension range retain at least
their preliminary weights. For a nonspecial path of distance
$k\in\{3,4,5,6,7\}$, the remaining question is whether it can acquire
an ordinary extension of length $9-k$ at both ends.

Follow each end for at most six steps. At every step retain the whole
candidate, including the marks omitted from the scored diagram, and
check all prescribed positive definite and semidefinite deletions.
One either follows an existing edge, adds a new node, or, from a node
outside the candidate, joins a different existing component. If one
end cannot be extended, the reduction is excluded and the preliminary
weight is a valid lower bound. If both ends pass these separate
necessary tests, use the smaller weight $1/2$ or $1/3$ specified in
\eqref{ext:eq:weight}. The two extensions need not coexist; allowing
the smaller weight nevertheless gives a valid lower bound. Thus the
deletion tests determine which reductions can be excluded. Every actual
extension gives a branch of the search: retain the whole $\overline L$ and the
portion of the extension being followed, remove other external nodes
and edges, and reduce heavy labels to four. These operations preserve
all required deletion conditions by Lemma~\ref{ext:lem:label-reduction}.
If the path meets an existing node, use that node rather than copying
it. Hence the failure of a necessary test excludes the reduction.
Algorithm~\ref{ext:alg:arm} specifies the finite search.

The classification, with every label at least four replaced by four
and isomorphic labeled diagrams identified, gives $61$ representatives.
All orders are treated by the same reduction. The numbers at orders
three through ten are
\[
5,\quad17,\quad13,\quad12,\quad3,\quad4,\quad4,\quad3.
\]
These counts follow directly from the classification in Appendix~\ref{app:lanner}.
At order three the reduced labels are $(2,3,4)$, $(2,4,4)$, $(3,3,4)$,
$(3,4,4)$ and $(4,4,4)$. At order four, the shapes are three stars,
five paths, two triangles with a leaf, five four-cycles, and the complete
diagram with either one edge removed or none: $3+5+2+5+1+1=17$.
At order five there are one four-leaf star, four trees with arm lengths
$(1,1,2)$, three paths, two four-cycles with a leaf, two five-cycles,
and one diagram formed by three two-edge paths sharing their endpoints:
$1+4+3+2+2+1=13$.
A reduced representative need not remain indefinite; all nodes and marks
are retained.

For a representative $\overline L$, put
\[
 r=|\overline L|,\qquad
 p=p(\overline L)=\#\{v\in V(\overline L):G_{\overline L\setminus\{v\}}\text{ is positive definite}\}.
\]
Every single-node deletion is positive semidefinite, and every
double-node deletion is positive definite. Consequently the numbers
of eligible markings are obtained by the following choices:
\begin{enumerate}[label=\textup{(\Roman*)},leftmargin=2.3em]
\item Choose any four nodes, giving $\binom r4$ markings.
\item Choose three nodes with positive definite single-node deletions,
giving $\binom p3$ markings.
\item Choose the unordered pair $\{b,c\}$ among those $p$ nodes,
then choose $a$ among the remaining $r-2$ nodes, giving
$(r-2)\binom p2$ markings.
\end{enumerate}
Here $\binom rk=0$ for $r<k$. In III the actual deletion at $a$ is
singular, but its restriction to the core, or its label reduction,
can be positive definite. Thus the necessary condition on this deletion
is positive semidefiniteness. Distinct choices of marked nodes are
retained even when a diagram automorphism relates them.

Table~\ref{ext:tab:marking-counts} is a direct count of these choices.
In its third column, $p:m$ means that $m$ representatives have $p$
positive definite single-node deletions. The three marking columns
are the sums of $\binom r4$, $\binom p3$, and $(r-2)\binom p2$
with these multiplicities.
\begin{table}[!htb]
\centering
\setlength{\tabcolsep}{4pt}
\caption{Reduced classification diagrams and their eligible markings.}
\label{ext:tab:marking-counts}
\begingroup
\small
\linespread{1}\selectfont
\begin{tabular}{r|r|l|r|r|r}
\toprule[1pt]
$r$ & Diagrams & $p:\text{multiplicity}$ & I & II & III\\\hline
3 & 5 & $3:5$ & 0 & 5 & 15\\\hline
4 & 17 & $0:2,\ 1:1,\ 2:3,\ 3:5,\ 4:6$ & 17 & 29 & 108\\\hline
5 & 13 & $2:2,\ 3:2,\ 4:5,\ 5:4$ & 65 & 62 & 234\\\hline
6 & 12 & $0:1,\ 1:1,\ 2:1,\ 3:2,\ 4:3,\ 5:4$ & 180 & 54 & 260\\\hline
7 & 3 & $6:3$ & 105 & 60 & 225\\\hline
8 & 4 & $7:4$ & 280 & 140 & 504\\\hline
9 & 4 & $7:1,\ 8:3$ & 504 & 203 & 735\\\hline
10 & 3 & $7:1,\ 8:1,\ 9:1$ & 630 & 175 & 680\\\hline
Total & 61 & & 1781 & 728 & 2761\\
\bottomrule[1pt]
\end{tabular}
\endgroup
\end{table}
The totals are $1781$, $728$, and $2761$ markings in types I, II,
and III, respectively.

Computer verification using exact arithmetic applies the deletion tests
$\mathcal C$ and Algorithm~\ref{ext:alg:arm} to every eligible marking,
then sums its pairwise lower bounds. The minimum totals are $3,2,1$
in I, II, III, respectively. Thus the finite checks establish the
required inequalities; the marking counts alone do not establish them. The calculation following this proof illustrates a
critical marking and the obstruction to shortening the extension threshold.

For type IV, a cross pair at distance at most twelve has weight at least
$1/3$. When all cross distances exceed twelve, the two minimal indefinite
subdiagrams admit only finitely many intersections. The reduction in Subsection~\ref{ext:app:long} leaves $92$ gluing checks;
exact computer verification leaves thirteen possibilities, each of score
at least $1/3$.
\end{proof}

\noindent\emph{A critical type-III configuration.}
Take the ordinary chain $v_0-v_1-\cdots-v_8$, attach a leaf $c$ only
to $v_2$, and call the resulting diagram $L$. Set $a=v_8$ and $b=v_1$.
All edge labels are three. With central node $v_2$, its arm lengths are
$(1,2,6)$, so Lemma~\ref{ext:lem:ordinary-arms} gives
\[
\eta_L=2-\left(\frac12+\frac23+\frac67\right)=-\frac1{42}<0.
\]
Deleting $a$ leaves the three arms $(1,2,5)$ and gives $\eta=0$.
Every other single-node deletion leaves positive definite components:
either paths with all labels three or a central component with shorter arms and
$\eta>0$. Thus $L$ is minimal indefinite and $p(L)=9$. Its contributions
to the three marking counts in Table~\ref{ext:tab:marking-counts} are
\[
\binom{10}{4}=210,\qquad
\binom93=84,\qquad
8\binom92=288.
\]

\begin{figure}[!htb]
\centering
\begin{tikzpicture}[x=.86cm,y=1cm,
 font={\small\linespread{1}\selectfont},line width=.55pt,
 boxed mark/.style={draw,rectangle,fill=white,minimum size=5.1mm,inner sep=1.3pt}]
\draw(0,0)--(9,0);
\draw(2,0)--(2,1.05);
\foreach\i in {0,2,3,4,5,6,7,9}{
 \draw[fill=white](\i,0)circle(1.7pt);
 \node at(\i,-.42){$v_{\i}$};}
\node[boxed mark](b)at(1,0){$b$};
\node at(1,-.42){$v_1$};
\node[boxed mark](a)at(8,0){$a$};
\node at(8,-.42){$v_8$};
\node[boxed mark]at(2,1.05){$c$};
\draw[densely dashed,black!55](-1.05,0)--(0,0);
\draw[fill=white,black!55](-1.05,0)circle(1.7pt);
\node[text=black!55]at(-1.05,-.42){$q$};
\draw[decorate,decoration={brace,mirror,amplitude=3pt}](0,-.82)--(8,-.82)
 node[midway,below=5pt]{$L$};
\end{tikzpicture}
\caption{The type-III diagram $H$, with $a=v_8$ and $b=v_1$.
All solid edges have label three. Removing $v_9$ gives the minimal
indefinite core $L$. The dashed ordinary edge $v_0q$ is not part of $H$;
it represents a proposed second outward step at $b$.}
\label{ext:fig:threshold}
\end{figure}
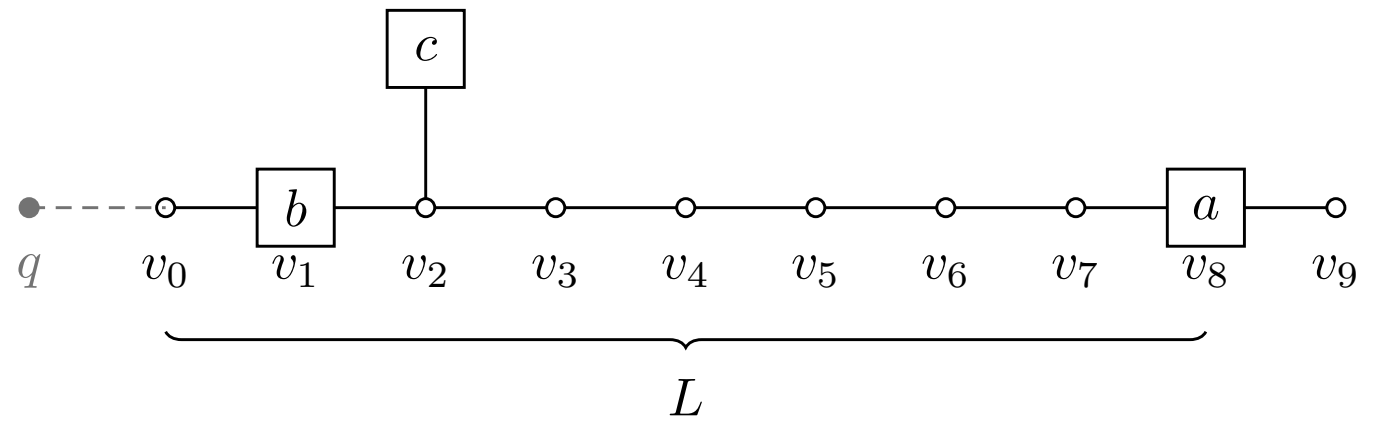

Now attach one more node $v_9$ to $a$, obtaining the full diagram $H$
in Figure~\ref{ext:fig:threshold}. Its arm lengths are $(1,2,7)$,
with $\eta_H=-1/24$, so it has exactly one negative eigenvalue.
For the marking $M=\{a,b,c\}$, the three deletions are
\[
\begin{aligned}
H\setminus\{a\}&\cong\widetilde E_8\sqcup A_1,
&&\text{singular positive semidefinite},\\
H\setminus\{b\}&\cong A_9\sqcup A_1,
&&\text{positive definite},\\
H\setminus\{c\}&\cong A_{10},
&&\text{positive definite}.
\end{aligned}
\]
The isolated nodes in the first two lines are $v_9$ and $v_0$,
respectively. Hence this marked diagram satisfies the local conditions
of III.

Both scored paths have distance seven. In $H\setminus\{c\}$, the
path from $a$ to $b$ has one outward edge at each end, namely $av_9$
and $bv_0$. In $H\setminus\{b\}$, the path from $a$ to $c$ has one
outward edge at $a$ and none at $c$. Thus the two-edge requirement
$9-7=2$ leaves both preliminary weights unchanged:
\[
w(H\setminus\{c\};a,b)+w(H\setminus\{b\};a,c)
=\frac12+\frac12=1.
\]
If the threshold were replaced by $8-k$, the first pair would satisfy
the one-edge requirement $8-7=1$ and receive weight $1/3$, while the
second would still receive $1/2$. The total would become
\begin{equation}\label{ext:eq:threshold-eight}
\frac13+\frac12=\frac56<1,
\end{equation}
violating the type-III lower bound.

The matrix obstruction can also be read directly. A second outward
step at $b$ would require an ordinary edge $v_0q$. After deleting $a$,
even the diagram with only this required new edge has a nontrivial
component with arm lengths $(1,3,5)$, giving
\[
\eta=2-\left(\frac12+\frac34+\frac56\right)=-\frac1{12}<0.
\]
Additional edges cannot restore positive semidefiniteness: deleting
them would have to preserve it by Lemma~\ref{ext:lem:label-reduction},
but the diagram with only $v_0q$ already fails the test. This explains
how the marked deletion at $a$ protects the half-unit contribution.
The example excludes $8-k$ within the stated family of rules; the
validity of $9-k$ for every configuration is the separate assertion
proved in Lemma~\ref{ext:lem:local}.

\subsubsection{Long configurations of type IV}\label{ext:app:long}
The remaining type-IV configurations require a separate reduction:
their two minimal indefinite cores need not lie in a single
classification diagram. Let $S$ be the elliptic unmarked diagram,
$X=\{x_1,x_2\}$, and $Y=\{y_1,y_2\}$.
Choose $L_X\subseteq S\cup X$ and $L_Y\subseteq S\cup Y$, minimal
among the indefinite induced subdiagrams of these respective diagrams.
Each contains its own marked pair, and $L_X\setminus\{x_i\}$ and
$L_Y\setminus\{y_i\}$ are positive definite. We may assume all cross-pair
distances exceed twelve, since otherwise one pair already contributes $1/3$.

If $L_X$ has order two, it is the dotted edge $X$ and contains no
node of $S$. The cores are then disjoint, and a joining edge $uv$
with $u\in X$ gives a cross pair at distance at most $1+7=8$.
The same holds for $L_Y$. We may therefore assume that both cores
have order at least three.

Put $d_X(u)=\min_{x\in X}\dist_{L_X}(u,x)$, and define $d_Y$ similarly.
A path to a nearest mark avoids the other mark, so joining such
paths produces a walk in a scored cross deletion.

\smallskip\noindent\emph{Classification data.}
Appendix~\ref{app:lanner} gives $d_X\leq7$. There are sixteen
marked classification diagrams having a position with $d_X\geq6$:
together they have seventeen positions at distance six and four
at distance seven. The four distance-seven positions belong to
different marked diagrams, each with an adjacent distance-six
position. For a marked classification diagram $(L,X)$, exact elimination
gives $\det G_{L\setminus\{u\}}=0$ at distance seven, and
\[
\max_{x\in X}(G_{L\setminus\{x\}}^{-1})_{uu}\geq\frac{16}{5}
\]
at distance six whenever $\det G_{L\setminus\{u\}}\ne0$.
The supplementary records contain the individual values.

\smallskip\noindent\emph{Reduction to overlapping cores.}
Disjoint cores must have a joining edge, or their Gram matrix would
have two negative directions. For every joining edge $uv$,
$d_X(u)+d_Y(v)+1\geq13$, so both distances are at least five.
There is at most one such edge: paths between its endpoints and
those of a second joining edge have length at most eight within
each core, whereas meeting a mark would require length at least
ten. They would therefore form a cycle in the elliptic forest $S$.

For this unique edge, write its Gram entry as $-t$, with $t\geq1/2$.
The determinant formula for the two blocks joined by this edge is
\[
\det G_{L_X\cup L_Y}
=\det G_{L_X}\det G_{L_Y}
-t^2\det G_{L_X\setminus\{u\}}\det G_{L_Y\setminus\{v\}}.
\]
Both core determinants are negative. If either cofactor vanishes,
the displayed determinant is positive; the existing negative
direction therefore forces a second one. Otherwise neither
attachment distance is seven, so both are six. Choose $x\in X$
and $y\in Y$ realizing the inverse-entry bounds above, and put
$A=G_{L_X\setminus\{x\}}$, $B=G_{L_Y\setminus\{y\}}$.
Positive definiteness of the corresponding cross deletion requires
$t^2(A^{-1})_{uu}(B^{-1})_{vv}<1$, whereas
\[
t^2(A^{-1})_{uu}(B^{-1})_{vv}
\geq\frac14\left(\frac{16}{5}\right)^2=\frac{64}{25}>1.
\]
Thus the cores overlap.

At a shared node $u$, $d_X(u)+d_Y(u)\geq13$, so at least one
distance is seven. The uniqueness of that position on each side
allows at most two shared nodes. Two shared nodes identify the
adjacent six--seven positions in reverse order. An additional
edge between the cores would again form an unmarked cycle.
Hence the overlap is a single node or a single edge.

\smallskip\noindent\emph{The finite check.}
Pair distinct marked diagrams without ordering them. For two copies
of the same diagram, retain both orientations of a one-node gluing
at different positions and identify inverse edge gluings.
There are
\[
17\cdot4+4+\binom52=82
\quad\text{node gluings},\qquad
\binom52=10\quad\text{edge gluings}.
\]
The extra four node gluings are the reversed six--seven
identifications within the same marked diagram. These are
$92$ gluing checks, not necessarily distinct isomorphism classes.
Exact arithmetic excludes $76$ by a non-positive-definite cross
deletion and three by inconsistent shared-edge labels; all
thirteen remaining gluings have score at least $1/3$.
For each gluing, the check deletes the two marks outside each
cross pair, tests positive definiteness, and evaluates
\eqref{ext:eq:weight}.

\smallskip\noindent\emph{A hand calculation at the lower bound.}
Take two copies of the ten-node diagram consisting of the path
$v_0-v_1-\cdots-v_8$ and a node $z$ joined only to $v_2$, with all
labels three. Identify $u=v_8^L=v_8^R$ and call the union $J$. Choose
\[
X=\{v_1^L,z^L\},\qquad Y=\{v_1^R,z^R\}.
\]
Put $S_J=J\setminus(X\cup Y)$. It has type $A_{13}\sqcup2A_1$, and
$S_J\cup X$, $S_J\cup Y$ each contain one of the minimal indefinite
cores. Deleting $u$ leaves $\wt E_8\sqcup\wt E_8$, which is
positive semidefinite. Interlacing therefore gives at most one negative
direction in $J$; the contained indefinite core gives exactly one.
For each cross pair $\{x,y\}$, delete the other two marks. The resulting scored diagram is
$A_{17}$ if neither endpoint is $z$, $A_{16}\sqcup A_1$ if exactly
one is $z$, and $A_{15}\sqcup2A_1$ if both are $z$. All four are
positive definite. Deleting the other mark on each side removes its
branch, and each retained endpoint lies seven edges from $u$. Thus
every scored path is nonspecial and has distance fourteen, giving
\[
\sigma(J)=\sum_{\substack{x\in X\\y\in Y}}
w\bigl(J\setminus((X\cup Y)\setminus\{x,y\});x,y\bigr)
=4\cdot\frac1{12}=\frac13.
\]
At these long distances, the extension rule makes no further reduction.
The scored paths in every passing gluing remain unchanged in an actual
elliptic extension, and special paths remain special. Their scores
therefore give the required lower bound for the full type-IV configuration.
This completes Lemma~\ref{ext:lem:local}. The supplementary
material supplies the exact computations and complete records.

\begin{prop}\label{ext:prop:face}
For the weight \eqref{ext:eq:weight}, every good 3-face satisfies
$\sigma(F)\geq7-f_2(F)$, while $\sigma(\mathcal{P})=0$. In particular,
\begin{equation}\label{ext:eq:extension-face-average}
a_w\geq7-a_3^2-b.
\end{equation}
\end{prop}
\begin{proof}
The reduction to the nine types in Figure~\ref{fig:small-faces} is unchanged
from Proposition~\ref{prop:facewise}. The tetrahedron uses case I, the
triangular prism case II, and the quadrilateral pyramid two instances of case III
on disjoint pairs of base edges. The other good six-facet types use case
II or III, except for the cube, whose twelve edges split into three sets
of four, each giving case IV. Lemma~\ref{ext:lem:local} therefore gives the
totals $3,2,1$ required for good 3-faces with four, five, and six facets.
For 3-faces with at least seven facets, nonnegativity suffices.
As in the argument for \eqref{eq:zero-bipyramid}, every scored pair
of $\mathcal{P}$ is disconnected, so all weights vanish. Summing, with one unit
subtracted for each bad 3-face, gives the displayed inequality.
\end{proof}

\subsection{The available weight at an edge}
We now specialize the coefficient $C$ in \eqref{eq:edge-target} to the
refined weight rule. Throughout Part~\ref{part:edge} set
\begin{equation}\label{ext:eq:edge-constants}
C=\frac{73}{12},\qquad \delta=\frac{23}{4},
\qquad \ell(T)=C|T|-W(T).
\end{equation}
\begin{lem}\label{ext:lem:budget}
Every nonempty elliptic diagram satisfies $\ell(T)\geq\delta$. For connected
diagrams of order $r\geq18$,
\begin{equation}\label{ext:eq:finite-tail}
\ell(A_r)=\frac{41}{3},\qquad
\ell(B_r)=\frac{127}{12},\qquad
\ell(D_r)=\frac{23}{4}.
\end{equation}
\end{lem}
\begin{proof}
A connected elliptic diagram of order greater than eight is of type
$A_r,B_r$, or $D_r$. For $r\geq18$, summing the preliminary weight by distance
gives, respectively,
\[
\frac{101r}{12}-\frac{89}{2},\qquad
\frac{101r}{12}-\frac{175}{4},\qquad
\frac{101r}{12}-\frac{165}{4}.
\]
For each distance $k\in\{3,4,5,6,7\}$, the numbers of pairs with the required
extension are $r-18+k$, $r-19+k$, and $r-20+k$ in these three families.
Subtract $1/2$ at distances three to five, $2/3$ at distance six, and $1/6$
at distance seven. This gives \eqref{ext:eq:finite-tail}.

For a connected diagram of order $r\leq8$, each pair has weight at most one,
so $\ell(T)\geq r(C-(r-1)/2)\geq\delta$. For $9\leq r\leq17$,
directly summing the fixed weight rule gives
\[
\ell(A_r)\geq\frac{41}{3},\qquad
\ell(B_r)\geq\frac{127}{12},\qquad
\ell(D_r)\geq\frac{23}{4}.
\]
The individual finite values are supplied with the supplementary
verification material. Additivity proves the result for disconnected diagrams.
\end{proof}

For later use, write the resulting edge coefficient as
\begin{equation}\label{ext:eq:lambda}
\lambda_n=\frac{2}{n-2}\left(C-\frac{\delta}{n-1}\right)
=\frac{2}{n-2}\left(\frac{73}{12}-\frac{23}{4(n-1)}\right).
\end{equation}
The incidence identities \eqref{eq:incidences} now give
\begin{equation}\label{ext:eq:basic-edge}
a_w=\frac1{a_3}\sum_e W(T_e)\leq\lambda_n a_3^1.
\end{equation}

\subsection{The weight bound and the two cusp improvements}
The reduction in the linear coefficient is $C_{\mathrm P}-C=43/12$.
Without using either new cusp estimate, we may discard the constant
improvement $\delta$ and apply \eqref{eq:prokhorov-lemma} with this value of $C$.
This immediately gives
\begin{equation}\label{ext:eq:weight-only-bound}
n<96C+68=96\cdot\frac{73}{12}+68=652,\qquad\text{hence }n\leq651.
\end{equation}
The bound $162$ also uses the two cusp estimates below.

To see what the next two sections must improve, place the lower and
upper bounds for the same average $a_w$ side by side. Equations
\eqref{ext:eq:extension-face-average}, \eqref{ext:eq:basic-edge}, and
\eqref{eq:euler-average}, together with the average face-number bound,
give
\begin{equation}\label{ext:eq:comparison-start}
\begin{aligned}
a_w&\geq7-a_3^2-b,\\
a_w&\leq\lambda_n a_3^1,\\
a_3^1&\leq3a_3^2-6-3b,\qquad a_3^2\leq u_{n,3},\qquad 0\leq b\leq1,
\end{aligned}
\end{equation}
where we write the universal upper bound as
\begin{equation}\label{ext:eq:un}
u_{n,3}=\begin{cases}
6+12/(n-4),&n\text{ even},\\[2pt]
6+12/(n-3),&n\text{ odd}.
\end{cases}
\end{equation}
The next two sections sharpen this comparison for $n\geq163$.

\smallskip\noindent\emph{Section~\ref{ext:sec:cusp}: lower the weight upper bound directly.}
Square cusp sections force an additional deficit in the edgewise
sum. Proposition~\ref{ext:prop:charge} gives
$a_w\leq\lambda_n(a_3^1-9b)$; with Euler's estimate this becomes

\begin{equation}\label{ext:eq:comparison-after-cusps}
7-a_3^2-b\leq a_w\leq\lambda_n(3a_3^2-6-12b).
\end{equation}
The coefficient $12b$ combines the $3b$ already present in Euler's
estimate with the new $9b$ from the weight deficit.

\smallskip\noindent\emph{Section~\ref{ext:sec:h}: improve the average face-number bound, and hence both sides.}
A smaller upper bound for $a_3^2$ improves both sides of
\eqref{ext:eq:comparison-after-cusps}. Writing the gap explicitly,

\begin{equation}\label{ext:eq:comparison-gap}
\begin{aligned}
a_w&\geq7-u_{n,3}-b+(u_{n,3}-a_3^2),\\
a_w&\leq\lambda_n(3u_{n,3}-6-12b)-3\lambda_n(u_{n,3}-a_3^2).
\end{aligned}
\end{equation}
Section~\ref{ext:sec:h} estimates this gap using the cusp terms
in Timorin's identities and the square sections of bad 3-faces.
Its bound \eqref{ext:eq:average} involves both $b$ and $a_3^1$;
Euler's relation eliminates the latter.
Section~\ref{ext:sec:conclusion} then proves incompatibility for
every $n\geq163$.

\section{Square sections and edgewise deficits}\label{ext:sec:cusp}
Square cusp sections allow us to sharpen the edgewise estimate
\eqref{ext:eq:comparison-start} to
$a_w\leq\lambda_n(a_3^1-9b)$. For an edge $e$, define its
\emph{weight deficit} by
\begin{equation}\label{ext:eq:edge-loss}
\ell_e:=\ell(T_e)-\delta=C(n-1)-\delta-W(T_e)\geq0.
\end{equation}
Here $C$, $\delta$, and $\ell$ are as in \eqref{ext:eq:edge-constants}.
By the incidence identity
\eqref{eq:incidences},
\begin{equation}\label{ext:eq:deficit-identity}
a_w=\lambda_n a_3^1-\frac1{a_3}\sum_e\ell_e.
\end{equation}
Consequently, to obtain the desired upper bound, it suffices to prove
\begin{equation}\label{le:deficit-identity}
\sum_e\ell_e\geq9\lambda_n B,
\end{equation}
where $B=ba_3$ is the number of bad 3-faces.

\subsection{Factor estimates}\label{ext:sec:metrics}
Throughout this section, $L_u$ denotes the cusp link at an ideal
vertex $u$, with Euclidean simplex factors $\Delta_i$. Write
\[
L_u=\prod_{i=1}^r\Delta_i,\qquad d_i=\dim\Delta_i,
\qquad\sum_i d_i=n-1.
\]
A vertex $v=(v_1,\ldots,v_r)$ of $L_u$ corresponds to an edge of $P$
incident to $u$. A 2-face of $L_u$ corresponds to a 3-face of $P$
containing $u$. In particular, each bad $3$-face through $u$
gives the square
\[
F\cap \mathcal H_u=\Delta_i^1\times\Delta_j^1\times\prod_{k\ne i,j}\{v_k\},
\qquad |\Delta_i^1|=|\Delta_j^1|,
\]
as shown in Figure~\ref{fig:bad-square-3d}. We bound the number of
\emph{squares arising from bad 3-faces} by the total number of square
2-faces of $L_u$. Let $\Sigma_i$ be the parabolic diagram of $\Delta_i$.
It has $d_i+1$ nodes, one for each facet. A vertex $v_i$ lies on exactly
$d_i$ facets. Write $\Sigma_i\setminus v_i$ for the possibly disconnected
elliptic diagram obtained by deleting the node corresponding to the
facet opposite $v_i$. Define the deficit at $v_i$ and its vertex average by
\begin{equation}\label{ext:eq:factor-loss}
\ell_{\Delta_i}(v_i)=Cd_i-W(\Sigma_i\setminus v_i),
\qquad
\bar\ell_{\Delta_i}=\frac1{d_i+1}
\sum_{v_i\in\operatorname{Vert}(\Delta_i)}\ell_{\Delta_i}(v_i).
\end{equation}
Thus $\ell_{\Delta_i}(v_i)=\ell(\Sigma_i\setminus v_i)$ in the
notation of \eqref{ext:eq:edge-constants}.

For the edge $e$ corresponding to a product vertex $v=(v_i)_{i=1}^r$,
orthogonality of the factors gives
\[
T_e=\bigsqcup_{i=1}^r(\Sigma_i\setminus v_i),
\qquad W(T_e)=\sum_{i=1}^rW(\Sigma_i\setminus v_i).
\]
Consequently,
\begin{equation}\label{ext:eq:factor-additivity}
\ell_e=\sum_{i=1}^r\ell_{\Delta_i}(v_i)-\delta.
\end{equation}
Each vertex of $\Delta_i$ occurs in $f_0(L_u)/(d_i+1)$ product
vertices, so
\begin{equation}\label{ext:eq:product-average}
\sum_{v\in\operatorname{Vert}(L_u)}\ell_{\Delta_i}(v_i)
=\frac{f_0(L_u)}{d_i+1}
  \sum_{v_i\in\operatorname{Vert}(\Delta_i)}\ell_{\Delta_i}(v_i)
=f_0(L_u)\bar\ell_{\Delta_i}.
\end{equation}

Let $B_u$ count bad 3-faces through $u$, and let $B_{u,ij}$ count
those whose square sections use the factors $i,j$. Each such square
uses exactly two factors, so $B_u=\sum_{i<j}B_{u,ij}$. Every bad 3-face has three ideal vertices, so
$\sum_{u\text{ ideal}}B_u=3B$.
We also define two sums of edge deficits:
\[
\begin{aligned}
I_u&:=\sum_{e\ni u}\ell_e
     =f_0(L_u)\left(\sum_{i=1}^r\bar\ell_{\Delta_i}-\delta\right),\\
J_u&:=\sum_{\substack{e\ni u\\\text{the other endpoint is finite}}}\ell_e.
\end{aligned}
\]
Their sums over all ideal vertices satisfy
\[
\sum_{u\text{ ideal}}(I_u+J_u)
=2\sum_{\substack{e\\\text{$e$ has at least one ideal endpoint}}}\ell_e
\leq2\sum_e\ell_e.
\]
We shall bound $I_u$ and $J_u$ from below in terms of $B_u$ to obtain
\eqref{le:deficit-identity}.

The sum $I_u$ uses all product vertices, so
\eqref{ext:eq:product-average} applies. The vertices selected by
$J_u$ need not use each factor vertex equally often; that sum
requires a pointwise estimate.

Use the edge-length counts of
Subsection~\ref{ext:app:equal-length}, writing
$m_i(\rho)=m_{\Delta_i}(\rho)$ and
$m_i(v_i;\rho)=m_{\Delta_i}(v_i;\rho)$. An equal-length pair of edges
in factors $i,j$ determines one square for each choice of vertices
in the other factors. Hence
\[
B_{u,ij}\leq\sum_\rho m_i(\rho)m_j(\rho)
\prod_{k\ne i,j}(d_k+1).
\]
Define
\begin{equation}\label{ext:eq:lengths}
M_i:=\frac{\sum_\rho m_i(\rho)^2}{(d_i+1)^2},
\qquad S_i(v_i):=\sum_\rho m_i(v_i;\rho)^2.
\end{equation}
Dividing the preceding count by $f_0(L_u)$ and applying
Cauchy--Schwarz gives
\begin{equation}\label{ext:eq:cauchy}
\frac{B_{u,ij}}{f_0(L_u)}\leq
\sum_\rho\frac{m_i(\rho)m_j(\rho)}{(d_i+1)(d_j+1)}
\leq\sqrt{M_iM_j}.
\end{equation}
Similarly, if $b_u(v)$ counts the squares arising from bad 3-faces
through a fixed product vertex $v=(v_i)$, then
\[
b_u(v)\leq\sum_{i<j}\sum_\rho m_i(v_i;\rho)m_j(v_j;\rho)
\leq\sum_{i<j}\sqrt{S_i(v_i)S_j(v_j)}.
\]
All length matches use the actual relative scales of the factors;
$M_i$ and $S_i(v_i)$ themselves are invariant under rescaling.

For a simplex $\Delta$ of dimension $d$, use the analogous notation
$M_\Delta,S_\Delta(v),\ell_\Delta(v),\bar\ell_\Delta$.
We compare the pointwise and mean deficits with $S_\Delta(v)/d$
and $M_\Delta/d$, respectively. A small adjustment is necessary:
for an equilateral triangle, $S_\Delta(v)=4$, $M_\Delta=1$, and
\[
\ell_\Delta(v)=2C-1=\frac C2S_\Delta(v)-1.
\]
The following constants allow this adjustment in the pointwise
estimate while recovering it in the mean estimate:
\begin{equation}\label{ext:eq:correction}
c=
\begin{cases}
1,&\Delta\text{ has type }\wt A_2,\\
2,&\Delta\text{ has type }\wt B_{d},\ d\geq9,\\
3,&\Delta\text{ has type }\wt D_{d},\ d\geq9,\\
0,&\text{otherwise}.
\end{cases}
\end{equation}
Each $c$ depends only on the affine type and dimension and does not
change the weight rule.

\begin{lem}\label{ext:lem:factor}
Let $\Delta$ be an irreducible Euclidean Coxeter simplex of dimension
$d$, and let $c$ be its constant in \eqref{ext:eq:correction}.
\begin{enumerate}[label=\textnormal{(\roman*)},leftmargin=2.2em]
\item At every vertex $v$,
\begin{equation}\label{ext:eq:point-factor}
\ell_\Delta(v)\geq\frac{C}{d}S_\Delta(v)-c.
\end{equation}
\item The mean deficit satisfies
\begin{equation}\label{ext:eq:mean-floor}
\bar\ell_\Delta\geq C+c.
\end{equation}
If $d\geq9$, it also satisfies
\begin{equation}\label{ext:eq:mean-factor}
\bar\ell_\Delta\geq\frac{4C}{d}M_\Delta+c.
\end{equation}
\end{enumerate}
\end{lem}
The two estimates apply to the selected vertices and all vertices,
respectively. The restriction $d\geq9$ in \eqref{ext:eq:mean-factor}
cannot simply be dropped: for $\wt A_2$,
$\bar\ell_\Delta-c=2C-2<4CM_\Delta/2$. Small factors will be
handled separately in Lemma~\ref{ext:lem:pair}.

\begin{proof}
The length multiplicities in
\eqref{ext:eq:A-length-counts} give
\begin{equation}\label{ext:eq:A-metric}
M_\Delta=\begin{cases}
d/2,&d\text{ even},\\
d/2-1/4,&d\text{ odd}
\end{cases}
\qquad(\Delta\text{ of type }\wt A_d).
\end{equation}
For the other classical families, the multiplicities in
\eqref{ext:eq:classical-length-counts} give, for $d\geq9$,
\begin{equation}\label{ext:eq:BCD-metric}
\begin{array}{c|c}
\text{Type of }\Delta&\sum_\rho m_\Delta(\rho)^2\\\hline
\wt C_{d}&d(d+1)(2d+1)/6\\
\wt B_{d}&d(d+1)(2d+1)/6-2d+2\\
\wt D_{d}&d(d+1)(2d+1)/6-4d+26
\end{array}
\end{equation}
After division by $(d+1)^2$, these formulas imply
$M_\Delta<d/3$ in types $\wt B_d,\wt C_d,\wt D_d$;
in type $\wt A_d$ we have $M_\Delta\leq d/2$.

For the coordinate vertices $v_i$ of Subsection~\ref{ext:app:equal-length},
the finite components obtained by deleting the opposite facet give
Table~\ref{ext:tab:deletion}. Assign deficit zero to the empty diagram, and use
$D_2=A_1\sqcup A_1$ and $D_3=A_3$.
\begin{table}[!htb]
\centering
\begingroup
\small
\linespread{1}\selectfont
\setlength{\tabcolsep}{6pt}
\caption{Deletion deficits at the coordinate vertices of a classical affine simplex.}
\label{ext:tab:deletion}
\begin{tabular}{c|l|l}
\toprule[1pt]
Type & Vertex indices & $\ell_\Delta(v_i)$\\\hline
$\wt A_d$ & all $i$ & $\ell(A_d)$\\\hline
$\wt C_d$ & $0\leq i\leq d$ & $\ell(B_i)+\ell(B_{d-i})$\\\hline
$\wt B_d$ & $i=0,1$ & $\ell(B_d)$\\\hline
& $2\leq i\leq d$ & $\ell(D_i)+\ell(B_{d-i})$\\\hline
$\wt D_d$ & $i=0,1,d-1,d$ & $\ell(D_d)$\\\hline
& $2\leq i\leq d-2$ & $\ell(D_i)+\ell(D_{d-i})$\\
\bottomrule[1pt]
\end{tabular}
\endgroup
\end{table}

\smallskip\noindent\emph{Mean deficits.}
Deleting a node of an affine diagram leaves a union of finite diagrams.
Using Table~\ref{ext:tab:deletion}, the finite deficits computed from the
weight rule, and \eqref{ext:eq:finite-tail} gives, for $d\geq20$,
\begin{equation}\label{ext:eq:means}
\begin{array}{c|c}
\text{Type}&\bar\ell_\Delta\\\hline
\wt A_{d}&41/3\\
\wt B_{d}&(49d+673)/[3(d+1)]\\
\wt C_{d}&(127d+1021)/[6(d+1)]\\
\wt D_{d}&(23d+557)/[2(d+1)]
\end{array}
\end{equation}
For instance, deleting node $i$ of $\wt C_{d}$ gives $B_i\sqcup B_{d-i}$,
so $\bar\ell_\Delta=2\sum_{j=1}^{d}\ell(B_j)/(d+1)$.
The finite values for $9\leq d\leq19$ satisfy the lower estimates
needed here. After subtraction of $c$, uniform bounds for
$d\geq9$ are
\[
\begin{array}{c|cccc}
\text{Type}&\wt A_{d}&\wt B_{d}&\wt C_{d}&\wt D_{d}\\\hline
\text{Lower bound for }\bar\ell_\Delta-c&41/3&19/2&23/2&17/2.
\end{array}
\]
Since $41/3>2C$ and $17/2>4C/3$, these prove
\eqref{ext:eq:mean-floor} and \eqref{ext:eq:mean-factor} in this range.

\smallskip\noindent\emph{Pointwise deficits.}
For type $A$, the coordinates give
$S_\Delta(v_i)=2d-(d\bmod2)$ at every vertex. Its deletion deficit is
$41/3>2C$ in the stable range. For type $C$,
$S_\Delta(v_i)=d+2\min(i,d-i)$. An internal deletion has two nonempty
$B$ components, each of deficit at least $C$. An end deletion has one
such component and $S_\Delta(v_i)=d$.
The type $\wt A_2$ was checked in the example preceding the lemma.

In type $\wt B_{d}$, the coordinate vertices $v_0,v_1,v_{d}$ have
$S_\Delta(v_i)\leq d+2$ and a single nonempty deletion component, with deficit
at least $\delta$. In particular, deleting the node opposite $v_{d}$
gives $D_{d}$. Every other vertex has $S_\Delta(v_i)\leq2d+2$ and at least
two nonempty components, with total deficit at least $2\delta$.
In type $\wt D_{d}$, the four end vertices have $S_\Delta(v_i)\leq d+4$
and deficit at least $\delta$. Internal vertices have
$S_\Delta(v_i)\leq2d+8$ and deficit at least $2\delta$.
Consequently, for $d\geq21$ it is enough to check
\begin{align*}
\delta+2&\geq C(1+2/d),&2\delta+2&\geq C(2+2/d),\\
\delta+3&\geq C(1+4/d),&2\delta+3&\geq C(2+8/d).
\end{align*}
All hold already at $d=21$. The smallest difference is $1/63$ in the
last inequality. In particular, the endpoint $v_{d}$ in type $B$ has
difference $(10d-73)/(6d)>0$.

\smallskip\noindent\emph{Small dimensions and exceptional types.}
For each classical case $d\leq20$, compute the squared edge lengths
from the coordinates in Subsection~\ref{ext:app:equal-length}, group equal
lengths to obtain $S_\Delta(v_i)$, and use Table~\ref{ext:tab:deletion}
to calculate $\ell_\Delta(v_i)$. Substitution gives
$\ell_\Delta(v_i)-C S_\Delta(v_i)/d+c\geq0$ at every vertex.
The same calculation for the five exceptional simplices uses their
explicit coordinates and the finite diagram obtained at each deletion;
it gives \eqref{ext:eq:point-factor} with $c=0$. All comparisons are exact;
the individual values are included in the supplementary material.
For $d\leq8$, every node pair has weight at most one, so
$W\leq\binom{d}{2}$. The correction is at most one in dimension two
and zero in the other dimensions in this range. Therefore
\begin{equation}\label{ext:eq:small-mean}
\bar\ell_\Delta-c
\geq d\left(C-\frac{d-1}{2}\right)-c
\geq C.
\end{equation}
The one-dimensional interval has $S_\Delta(v)=1$, $\ell_\Delta(v)=C$, and $c=0$.
This completes all cases.
\end{proof}

\subsection{Summing the edgewise deficits}\label{ext:sec:charge}
Retain the factor notation of Subsection~\ref{ext:sec:metrics}, and
write $c_i$ for the constant of $\Delta_i$ in
\eqref{ext:eq:correction}. We shall prove
$I_u+J_u\geq6\lambda_n B_u$. Since each bad 3-face has three ideal
vertices and $\sum_u(I_u+J_u)\leq2\sum_e\ell_e$, this yields
\eqref{le:deficit-identity}.

The first step compares the square count for two factors with their
mean deficits.
\begin{lem}\label{ext:lem:pair}
For every pair of factors of $L_u$,
\begin{equation}\label{ext:eq:pair}
8C\frac{B_{u,ij}}{f_0(L_u)}
\leq d_i(\bar\ell_{\Delta_j}-c_j)+d_j(\bar\ell_{\Delta_i}-c_i).
\end{equation}
\end{lem}
\begin{proof}
If $d_i,d_j\geq9$, Lemma~\ref{ext:lem:factor} and the arithmetic--geometric
mean inequality give
\[
\begin{aligned}
d_i(\bar\ell_{\Delta_j}-c_j)+d_j(\bar\ell_{\Delta_i}-c_i)
&\geq4C\left(\frac{d_i}{d_j}M_j+\frac{d_j}{d_i}M_i\right)\\
&\geq8C\sqrt{M_iM_j}
\geq8C\frac{B_{u,ij}}{f_0(L_u)}.
\end{aligned}
\]
It remains to treat pairs with a factor of dimension at most eight.
We use the length multiplicity bound and the alternating labels to
bound their square count, and then compare with the mean deficits.

\smallskip\noindent\emph{Counting when a factor is small.}
For any pair $i,j$, counting all edge pairs
and then using the length multiplicity bound in Lemma~\ref{ext:lem:equal-length}
gives
\begin{equation}\label{ext:eq:qgeneral}
\frac{B_{u,ij}}{f_0(L_u)}
\leq\min\left\{\frac{d_i d_j}{4},\frac35\min(d_i,d_j)\right\}.
\end{equation}
The first bound counts all edge pairs. For the second, each edge in
factor $i$ matches at most $6(d_j+1)/5$ edges in factor $j$;
interchanging the factors gives the minimum. If at least one dimension
is at least nine, the multiplicity bound for that factor improves to
its number of vertices, and
\[
\frac{B_{u,ij}}{f_0(L_u)}\leq\frac12\min(d_i,d_j).
\]

The smallest pairs admit a stronger bound because the finite and ideal
labels alternate around each square arising from a bad 3-face. Color
finite-labeled vertices zero and ideal-labeled vertices one, and fix
all other factor coordinates.
For dimensions $1,2$, the product vertices form a $2\times3$ array.
An alternating rectangle pairs a column $(0,1)$ with a column $(1,0)$,
so there are at most two such rectangles among its six vertices.
For dimensions $2,2$, a $3\times3$ array has at most four alternating
rectangles. Indeed, if two rows agree, they give none together and
at most two each with the third row. If all rows differ, at most one
pair is complementary and contributes two. Each other pair contributes
at most one. Applying these bounds in each fixed-coordinate array yields
\begin{equation}\label{ext:eq:small-color}
\frac{B_{u,ij}}{f_0(L_u)}\leq
\begin{cases}
1/3,&\{d_i,d_j\}=\{1,2\},\\
4/9,&d_i=d_j=2.
\end{cases}
\end{equation}

\smallskip\noindent\emph{Comparing with the mean deficits.}
Assume $d_i\leq d_j$. For $d_i\leq d_j\leq8$ with $d_j\geq3$,
the explicit lower bound \eqref{ext:eq:small-mean} gives
\[
\frac{\bar\ell_{\Delta_i}-c_i}{d_i}\geq C-\frac{d_j-1}{2},
\qquad
\bar\ell_{\Delta_j}-c_j\geq d_j\left(C-\frac{d_j-1}{2}\right)\geq\frac{61}{4}.
\]
The possible correction at $d_i=2$ gives a lower bound $C-1$
per dimension, sufficient since $d_j\geq3$. The second bound follows
by taking the minimum of the concave quadratic
$x(C-(x-1)/2)$ on $[3,8]$.
Consequently,
\[
\begin{aligned}
d_i(\bar\ell_{\Delta_j}-c_j)+d_j(\bar\ell_{\Delta_i}-c_i)
&\geq2d_i d_j\left(C-\frac{d_j-1}{2}\right)\\
&\geq\frac{61}{2}d_i>\frac{24C}{5}d_i.
\end{aligned}
\]
Together with $B_{u,ij}/f_0(L_u)\leq3d_i/5$ from
\eqref{ext:eq:qgeneral}, this proves \eqref{ext:eq:pair}.

For $d_j\leq2$, the remaining pairs are $(1,1),(1,2),(2,2)$.
The upper bounds for $B_{u,ij}/f_0(L_u)$ are respectively
$1/4,1/3,4/9$, by \eqref{ext:eq:qgeneral} and
\eqref{ext:eq:small-color}. The corresponding lower bounds for
$d_i(\bar\ell_{\Delta_j}-c_j)+d_j(\bar\ell_{\Delta_i}-c_i)$ from
\eqref{ext:eq:small-mean} are
\[
2C=8C\cdot\frac14,\qquad
4C-2\geq8C\cdot\frac13,\qquad
8C-8\geq8C\cdot\frac49.
\]
Thus \eqref{ext:eq:pair} also holds for these three pairs.

Finally, if $d_i\leq8<d_j$, the large-factor mean bounds proved in
Lemma~\ref{ext:lem:factor} and
\eqref{ext:eq:small-mean} give
\[
\bar\ell_{\Delta_j}-c_j\geq\frac{17}{2},\qquad
\frac{\bar\ell_{\Delta_i}-c_i}{d_i}\geq\frac{31}{12}.
\]
Hence
\[
\frac{d_i(\bar\ell_{\Delta_j}-c_j)+d_j(\bar\ell_{\Delta_i}-c_i)}{d_i}
\geq\frac{17}{2}+9\cdot\frac{31}{12}
=\frac{127}{4}>4C.
\]
This proves \eqref{ext:eq:pair} using the stronger bound
$B_{u,ij}/f_0(L_u)\leq d_i/2$, and completes all remaining cases.
\end{proof}

\begin{prop}\label{ext:prop:charge}
For $n\geq3$, $a_w\leq\lambda_n(a_3^1-9b)$.
\end{prop}
\begin{proof}
We estimate $I_u$ and $J_u$ separately at a fixed ideal vertex $u$.

\smallskip\noindent\emph{All edges at the cusp.}
Equations~\eqref{ext:eq:factor-additivity} and
\eqref{ext:eq:product-average} give the exact identity
\[
I_u=f_0(L_u)\left(\sum_i\bar\ell_{\Delta_i}-\delta\right).
\]
Summing \eqref{ext:eq:pair} over $i<j$ gives
\[
\begin{aligned}
8C\frac{B_u}{f_0(L_u)}
&\leq\sum_i(n-1-d_i)(\bar\ell_{\Delta_i}-c_i)\\
&=(n-1)\sum_i(\bar\ell_{\Delta_i}-c_i)-\sum_i d_i(\bar\ell_{\Delta_i}-c_i)\\
&\leq(n-1)\left(\sum_i(\bar\ell_{\Delta_i}-c_i)-C\right).
\end{aligned}
\]
Here $\bar\ell_{\Delta_i}-c_i\geq C$ and $\sum_i d_i=n-1$.
Substituting into the identity for $I_u$ gives
\begin{equation}\label{ext:eq:I}
I_u\geq\frac{8C}{n-1}B_u
       +f_0(L_u)(C-\delta)+f_0(L_u)\sum_i c_i.
\end{equation}

\smallskip\noindent\emph{Edges whose other endpoint is finite.}
Recall that $b_u(v)$ counts the squares arising from bad 3-faces
through the product vertex $v=(v_i)$. The matching count from
Subsection~\ref{ext:sec:metrics} gives
\[
b_u(v)\leq\sum_{i<j}\sum_\rho m_i(v_i;\rho)m_j(v_j;\rho)
\leq\sum_{i<j}\sqrt{S_i(v_i)S_j(v_j)}.
\]
For each pair of factors, the pointwise estimate gives
\[
\begin{aligned}
2C\sqrt{S_i(v_i)S_j(v_j)}
&\leq C\left(\frac{d_i}{d_j}S_j(v_j)
                    +\frac{d_j}{d_i}S_i(v_i)\right)\\
&\leq d_i\bigl(\ell_{\Delta_j}(v_j)+c_j\bigr)
      +d_j\bigl(\ell_{\Delta_i}(v_i)+c_i\bigr).
\end{aligned}
\]
Since the length multiplicities are nonnegative integers,
$S_i(v_i)\geq\sum_\rho m_i(v_i;\rho)=d_i$; hence
$\ell_{\Delta_i}(v_i)+c_i\geq C$. Summing over factor pairs gives
\[
2C b_u(v)
\leq\sum_i(n-1-d_i)\bigl(\ell_{\Delta_i}(v_i)+c_i\bigr)
\leq(n-1)\left(\sum_i\bigl(\ell_{\Delta_i}(v_i)+c_i\bigr)-C\right).
\]
By \eqref{ext:eq:factor-additivity}, the edge $e$ corresponding to $v$
therefore satisfies
\[
\ell_e\geq\frac{2C}{n-1}b_u(v)+C-\delta-\sum_i c_i.
\]

Let $f_u$ count the edges from $u$ to finite vertices. Each square
arising from a bad 3-face has exactly two finite-labeled corners, so
\[
\sum_{v\text{ finite-labeled}}b_u(v)=2B_u.
\]
Summing the pointwise inequality over these $f_u$ vertices yields
\begin{equation}\label{ext:eq:J}
J_u\geq\frac{4C}{n-1}B_u
       +f_u(C-\delta)-f_u\sum_i c_i.
\end{equation}
\smallskip\noindent\emph{Combining the two sums.}
Adding \eqref{ext:eq:I} and \eqref{ext:eq:J} gives
\begin{equation}\label{ext:eq:transfer}
\begin{aligned}
I_u+J_u\geq{}&\frac{12C}{n-1}B_u
 +(f_0(L_u)+f_u)(C-\delta)\\
&+(f_0(L_u)-f_u)\sum_i c_i.
\end{aligned}
\end{equation}
The last term is nonnegative because $f_u\leq f_0(L_u)$ and
$c_i\geq0$. We retain the middle term, since $C-\delta=1/3$.
To bound it in terms of $B_u$, count all rectangular 2-faces of
$L_u$, without requiring equal side lengths or a bad ambient face.
At a product vertex, choosing one incident edge in factor $i$ and
one in factor $j$ gives $d_i d_j$ such rectangles. Their total is
\[
\sum_{i<j}d_i d_j
=\frac{(n-1)^2-\sum_i d_i^2}{2}
\leq\frac{(n-1)(n-2)}2,
\]
since $d_i^2\geq d_i$ and $\sum_i d_i=n-1$.
Apply this incidence count to the $B_u$ squares arising from bad 3-faces.
Counting all four corners, and then only their two finite-labeled
corners, gives respectively
\[
4B_u\leq f_0(L_u)\frac{(n-1)(n-2)}2,
\qquad
2B_u\leq f_u\frac{(n-1)(n-2)}2.
\]
Consequently,
\[
f_0(L_u)+f_u\geq\frac{12B_u}{(n-1)(n-2)}.
\]
Substituting in \eqref{ext:eq:transfer} and discarding its nonnegative
last term gives
\begin{equation}\label{ext:eq:cusp-loss}
I_u+J_u\geq
\left[\frac{12C}{n-1}
 +\frac{12(C-\delta)}{(n-1)(n-2)}\right]B_u
=6\lambda_n B_u.
\end{equation}
The last equality uses \eqref{ext:eq:lambda}.

\smallskip\noindent\emph{Summing over all cusps.}
Using $\sum_u(I_u+J_u)\leq2\sum_e\ell_e$ and
$\sum_uB_u=3B$ in \eqref{ext:eq:cusp-loss}, we obtain
\begin{equation}\label{ext:eq:total-loss}
2\sum_e\ell_e\geq6\lambda_n\sum_uB_u=18\lambda_n B,
\qquad\text{hence}\qquad
\sum_e\ell_e\geq9\lambda_n B.
\end{equation}
Finally, \eqref{ext:eq:deficit-identity} gives
\[
a_w\leq\lambda_n a_3^1-9\lambda_n\frac{B}{a_3}
=\lambda_n(a_3^1-9b).
\]
\end{proof}

Together with Euler's relation \eqref{eq:euler-average}, this proves
the strengthened upper bound in \eqref{ext:eq:comparison-after-cusps}.
The next section constrains the remaining average $a_3^2$, thereby
improving both sides of that comparison.

\section{Cusp geometry and average face numbers}\label{ext:sec:h}

We next bound the average number $a_3^2$ of facets of a 3-face
in terms of $b$. Retaining the cusp contributions in the $h$-vector
identities gives a lower bound for $u_{n,3}-a_3^2$; the square counts
then relate this gap to the number of bad 3-faces. All local estimates
in this section hold for $n\geq6$.

\subsection{Deriving the cusp contribution from truncation}

The incidence identity $a_3a_3^2=(n-2)a_2$ gives
\begin{equation}\label{ext:eq:normalized-gap}
u_{n,3}-a_3^2=\frac{n-2}{a_3}\left(\frac{u_{n,3}}{n-2}a_3-a_2\right).
\end{equation}
Expanding $a_3$ and $a_2$ using \eqref{ext:eq:h-inverse}, we obtain
\begin{equation}\label{ext:eq:gap-expansion}
u_{n,3}-a_3^2
=\frac{n-2}{a_3}\sum_{j=0}^n
 \left[\frac{u_{n,3}}{n-2}\binom j3-\binom j2\right]h_j(P).
\end{equation}
For $1\leq j<n/2$, substitute
$h_j(P)=h_{n-j}(P)-\sum_u g_j(L_u)$ from
\eqref{ext:eq:asymmetry} and pair degrees $j$ and $n-j$.

The degree-zero term vanishes, and the degree-$n$ term inside the sum
is $n(n-1)(u_{n,3}-3)/6$, since $h_n(P)=1$.
For even $n$ the unpaired middle term has coefficient
\[
\binom{n/2}{2}
\left[\frac{u_{n,3}(n/2-2)}{3(n-2)}-1\right]=0.
\]
Thus the gap takes the form
\begin{equation}\label{ext:eq:gap-pair}
\begin{aligned}
u_{n,3}-a_3^2={}&\frac{n-2}{a_3}\Biggl(
 \sum_u\sum_{j=2}^{\lfloor(n-1)/2\rfloor}
 g_j(L_u)\left[\binom j2-\frac{u_{n,3}}{n-2}\binom j3\right]\\
 &\qquad+R_n(P)\Biggr).
\end{aligned}
\end{equation}
The first term in parentheses is a sum over cusp links. For a product
link $L$ of dimension $n-1$, define
\begin{equation}\label{ext:eq:theta}
\Theta_n^{(3)}(L):=\sum_{j=2}^{\lfloor(n-1)/2\rfloor}
 g_j(L)\left[\binom j2-\frac{u_{n,3}}{n-2}\binom j3\right].
\end{equation}
Here the remainder is
\[
\begin{aligned}
R_n(P)={}&\frac{n(n-1)(u_{n,3}-3)}6\\
&+\sum_{j=1}^{\lfloor(n-1)/2\rfloor}
\left[\frac{u_{n,3}}{n-2}\left(\binom j3+\binom{n-j}3\right)
-\binom j2-\binom{n-j}2\right]h_{n-j}(P).
\end{aligned}
\]

\begin{lem}\label{ext:lem:retained}
For $n\geq6$, $R_n(P)\geq0$.
\end{lem}
\begin{proof}
The degree-$n$ contribution $n(n-1)(u_{n,3}-3)/6$ is positive because
$u_{n,3}>3$. For $1\leq j\leq\lfloor(n-1)/2\rfloor$, the coefficient of
$h_{n-j}(P)$ in the paired terms is, by \eqref{ext:eq:un},
\begin{equation}\label{ext:eq:paired-coeff}
\begin{aligned}
&\frac{u_{n,3}}{n-2}\left(\binom j3+\binom{n-j}3\right)
       -\binom j2-\binom{n-j}2\\
&\qquad=\begin{cases}
\displaystyle\frac{(2j-n)^2(n-1)}{2(n-4)},&n\text{ even},\\[4pt]
\displaystyle\frac{n\bigl((2j-n)^2-1\bigr)}{2(n-3)},&n\text{ odd}.
\end{cases}
\end{aligned}
\end{equation}
These expressions are nonnegative. In the odd case this follows
because $2j-n$ is an odd integer. Since $n-j>n/2$,
\eqref{ext:eq:h-upper-nonnegative} gives $h_{n-j}(P)\geq0$,
completing the proof.
\end{proof}

The retained cusp contribution is also nonnegative.
For every integer $0\leq j\leq\lfloor(n-1)/2\rfloor$, direct factoring gives
\begin{equation}\label{ext:eq:theta-coefficient}
\begin{aligned}
\binom j2-\frac{u_{n,3}}{n-2}\binom j3
&=\frac{j(j-1)(n-2j)}{2(n-4)}&& (n\text{ even}),\\
&=\frac{j(j-1)(n-2j)}{2(n-4)}
 +\frac{2j(j-1)(j-2)}{(n-4)(n-3)(n-2)}&& (n\text{ odd}).
\end{aligned}
\end{equation}
These coefficients are nonnegative throughout the summation range,
as are $g_j(L)$ by Lemma~\ref{ext:lem:link-polynomial}. Thus
Lemma~\ref{ext:lem:retained} yields
\begin{equation}\label{ext:eq:retained}
u_{n,3}-a_3^2\geq\frac{n-2}{a_3}\sum_u\Theta_n^{(3)}(L_u)\geq0.
\end{equation}

For a simplex link all $g_j$ with $j\geq1$ vanish, so $\Theta_n^{(3)}(L)=0$.
For $L=\Delta^2\times\Delta^{n-3}$, $n\geq6$, the lower-half differences
are $g_0=g_1=g_2=1$ and all later ones are zero. Hence $\Theta_n^{(3)}(L)=1$.
It therefore contributes $(n-2)/a_3$ to this lower bound for the gap.

\subsection{Normalizing the cusp term by the number of vertices}
The normalized coefficients in \eqref{ext:eq:normalized-g}
are $\pi_j=(n-2j)g_j(L)/f_0(L)$. Their sum is one. We apply these
common weights to the 3-face cusp term.

Since $g_j(L)/f_0(L)=\pi_j/(n-2j)$ and the coefficients at $j=0,1$
are zero, \eqref{ext:eq:theta} becomes the exact weighted average
\begin{equation}\label{ext:eq:theta-weighted}
\frac{\Theta_n^{(3)}(L)}{f_0(L)}
=\sum_j\pi_j\frac{\binom j2-\dfrac{u_{n,3}}{n-2}\binom j3}{n-2j}.
\end{equation}
Since $n-2j>0$, dividing \eqref{ext:eq:theta-coefficient} by $n-2j$
and using $N=\sum_jj(j-1)\pi_j$ gives
\begin{equation}\label{ext:eq:theta-moment}
\frac{\Theta_n^{(3)}(L)}{f_0(L)}\geq\frac{N}{2(n-4)},
\end{equation}
with equality for even $n$.

For the common quantities $m,N,e_2$, Lemma~\ref{ext:lem:moments} gives
\[
\sum_j j(n-j)\pi_j=e_2/2,\qquad m\geq e_2/(2n),\qquad
N\geq m^2-m,\quad m\leq N+1,\quad m^2\leq2N+1.
\]
These estimates relate the retained cusp term to the factor dimensions.

\subsection{Counting bad 3-faces at a cusp}\label{ext:sec:average}

Retain $B_u=\sum_{i<j}B_{u,ij}$ and $L_u=\prod_i\Delta_i$ from
Subsection~\ref{ext:sec:metrics}. For each edge of $\Delta_i$,
Lemma~\ref{ext:lem:equal-length} gives at most $6(d_j+1)/5$
equally long edges of $\Delta_j$. As in \eqref{ext:eq:qgeneral},
\begin{equation}\label{ext:eq:normalized-square-count}
\frac{B_{u,ij}}{f_0(L_u)}
\leq\frac{\binom{d_i+1}{2}\,\frac65(d_j+1)}{(d_i+1)(d_j+1)}
=\frac35d_i.
\end{equation}
Interchanging the factors and summing yields
\begin{equation}\label{ext:eq:q-square-sum}
\frac{B_u}{f_0(L_u)}\leq\frac35\sum_{i<j}\min(d_i,d_j).
\end{equation}
This ratio is a number of square 2-faces divided by a number of
vertices, not a proportion; it need not be at most one.

\subsection{From the square count to the retained cusp term}

We now relate the dimension sum in \eqref{ext:eq:q-square-sum}
to the moment $N$ and hence to the retained cusp term.

\begin{lem}\label{ext:lem:local-defect}
For $n\geq6$, every cusp satisfies
\begin{equation}\label{ext:eq:local-theta}
B_u\leq\mathcal A_n\Theta_n^{(3)}(L_u)+\mathcal B_nf_0(L_u),
\end{equation}
where
\begin{equation}\label{ext:eq:local-defect-constants}
\mathcal A_n=\frac{12n(9n-1)(n-4)}{5(n-1)^2},\qquad
\mathcal B_n=\frac{6n(5n-1)}{5(n-1)^2}.
\end{equation}
Equivalently,
\begin{equation}\label{ext:eq:local-theta-lower}
\Theta_n^{(3)}(L_u)\geq
\frac{B_u-\mathcal B_nf_0(L_u)}{\mathcal A_n}.
\end{equation}
\end{lem}
\begin{proof}
Reorder the factor dimensions so that $d_1$ is largest, and put
$d_{\mathrm{rest}}=n-1-d_1$. Since the dimensions are positive integers,
\begin{align}
\sum_{i<j}\min(d_i,d_j)
&\leq d_{\mathrm{rest}}+\sum_{2\leq i<j\leq r}d_i d_j\nonumber\\
&=d_{\mathrm{rest}}+\frac12\left(d_{\mathrm{rest}}^2-\sum_{i=2}^r d_i^2\right)
\leq d_{\mathrm{rest}}+\frac{d_{\mathrm{rest}}^2-d_{\mathrm{rest}}}{2}=\frac{d_{\mathrm{rest}}(d_{\mathrm{rest}}+1)}2.
\label{ext:eq:dimension-outside-largest}
\end{align}
Also $\sum_i d_i^2\leq d_1(n-1)$, so
\begin{equation}\label{ext:eq:e2-controls-t}
e_2=\frac{(n-1)^2-\sum_i d_i^2}{2}
\geq\frac{(n-1)(n-1-d_1)}2=\frac{(n-1)d_{\mathrm{rest}}}{2}.
\end{equation}
Lemma~\ref{ext:lem:moments} also gives $e_2\leq2nm$, so
\[
0\leq d_{\mathrm{rest}}\leq\frac{2e_2}{n-1}\leq\frac{4n}{n-1}m.
\]
Since $d_{\mathrm{rest}}(d_{\mathrm{rest}}+1)$ increases for $d_{\mathrm{rest}}\geq0$, equations
\eqref{ext:eq:q-square-sum} and \eqref{ext:eq:dimension-outside-largest}
yield
\begin{equation}\label{ext:eq:q-in-m}
\begin{aligned}
\frac{B_u}{f_0(L_u)}&\leq\frac3{10}d_{\mathrm{rest}}(d_{\mathrm{rest}}+1)\\
 &\leq\frac3{10}\left[\left(\frac{4n}{n-1}m\right)^2
                              +\frac{4n}{n-1}m\right]\\
 &=\frac{24n^2}{5(n-1)^2}m^2+\frac{6n}{5(n-1)}m.
\end{aligned}
\end{equation}
Apply \eqref{ext:eq:first-second-control}, retaining the constant terms:
\begin{align}
\frac{B_u}{f_0(L_u)}&\leq\frac{24n^2}{5(n-1)^2}(2N+1)
       +\frac{6n}{5(n-1)}(N+1)\nonumber\\
 &=\frac{48n^2+6n(n-1)}{5(n-1)^2}N
   +\frac{24n^2+6n(n-1)}{5(n-1)^2}\nonumber\\
 &=\frac{6n(9n-1)}{5(n-1)^2}N
   +\frac{6n(5n-1)}{5(n-1)^2}.
\label{ext:eq:q-in-N}
\end{align}
Finally, \eqref{ext:eq:theta-moment} gives
$N\leq2(n-4)\Theta_n^{(3)}(L_u)/f_0(L_u)$.
Substitution into \eqref{ext:eq:q-in-N} proves the stated bounds.
\end{proof}

The additive term is needed: for $L=\Delta^4\times\Delta^1$,
one has $m=2/5$ but $N=\Theta_6^{(3)}(L)=0$. Thus the constants in
$m\leq N+1$ and $m^2\leq2N+1$ cannot be omitted. In particular,
\eqref{ext:eq:local-theta-lower} gives a positive lower bound only if
$B_u>\mathcal B_nf_0(L_u)$.

\subsection{Summing the local estimates over the cusps}

\begin{prop}\label{ext:prop:average}
For $n\geq6$,
\begin{equation}\label{ext:eq:average}
\begin{aligned}
a_3^2\leq{}&u_{n,3}
-\frac{5(n-1)^2(n-2)}{4n(n-4)(9n-1)}\,b\\
&+\frac{2(5n-1)}{(n-1)(n-4)(9n-1)}\,a_3^1.
\end{aligned}
\end{equation}
\end{prop}
\begin{proof}
Every bad 3-face has exactly three ideal vertices. Counting its
incidences with cusps gives $\sum_u B_u=3B$, where $B=ba_3$ is the
total number of bad 3-faces. A vertex of $L_u$ corresponds to an
edge of $P$ incident to $u$. An ambient edge has at most two ideal
endpoints, so $\sum_u f_0(L_u)\leq2a_1$. The edge--3-face incidence
identity is
\[
\binom{n-1}{2}a_1=a_3a_3^1,
\quad\text{hence}\quad
\sum_u f_0(L_u)\leq\frac{4a_3a_3^1}{(n-1)(n-2)}.
\]
Summing \eqref{ext:eq:local-theta-lower} therefore gives
\[
\sum_u\Theta_n^{(3)}(L_u)
\geq\frac{3ba_3}{\mathcal A_n}
 -\frac{4\mathcal B_na_3a_3^1}
        {\mathcal A_n(n-1)(n-2)}.
\]
Substitution in the retained-gap inequality \eqref{ext:eq:retained}
yields
\[
u_{n,3}-a_3^2\geq\frac{3(n-2)}{\mathcal A_n}b
 -\frac{4\mathcal B_n}{\mathcal A_n(n-1)}a_3^1.
\]
Substituting \eqref{ext:eq:local-defect-constants} proves
\eqref{ext:eq:average}.
\end{proof}

The term involving $a_3^1$ comes from the additive vertex term in
\eqref{ext:eq:local-theta}; it remains in the final comparison.

\section{The bound \texorpdfstring{$n\leq162$}{n <= 162}}\label{ext:sec:conclusion}
Combining Proposition~\ref{ext:prop:charge},
Proposition~\ref{ext:prop:average}, the cellwise lower bound, and
Euler's relation gives
\begin{equation}\label{ext:eq:coupled}
\begin{gathered}
7-a_3^2-b\leq a_w\leq\lambda_n(3a_3^2-6-12b),\\
a_3^2\leq u_{n,3}
-\frac{5(n-1)^2(n-2)}{4n(n-4)(9n-1)}\,b
+\frac{2(5n-1)}{(n-1)(n-4)(9n-1)}\,a_3^1,\\
a_3^1\leq3a_3^2-6-3b,\qquad 0\leq b\leq1.
\end{gathered}
\end{equation}
At leading order, the weight improvement changes the comparison
from $1-b\lesssim146/n$ to
$1-b\lesssim(146/n)(1-b)$. This explains the gain but is not the dimension
argument: the finite-dimensional terms remain significant, particularly
near $b=1$.

For the elimination, write
\begin{equation}\label{ext:eq:global-coefficients}
\begin{aligned}
\kappa_n&:=\frac{3(n-2)}{\mathcal A_n}
 =\frac{5(n-1)^2(n-2)}{4n(n-4)(9n-1)},\\
\varepsilon_n&:=\frac{4\mathcal B_n}{\mathcal A_n(n-1)}
 =\frac{2(5n-1)}{(n-1)(n-4)(9n-1)}.
\end{aligned}
\end{equation}
Since $\varepsilon_n>0$, the last two lines of
\eqref{ext:eq:coupled} give
\[
(1-3\varepsilon_n)a_3^2
\leq u_{n,3}-6\varepsilon_n-(\kappa_n+3\varepsilon_n)b.
\]
For $n\geq163$, $1-3\varepsilon_n>0$. The coefficient of $a_3^2$
in the lower weight bound minus the upper weight bound is
$-(1+3\lambda_n)<0$. Substituting its upper bound therefore gives
the necessary inequality
\begin{equation}\label{ext:eq:necessary}
0\geq\Phi_n(b):=7+6\lambda_n+(12\lambda_n-1)b
-\frac{1+3\lambda_n}{1-3\varepsilon_n}
\bigl(u_{n,3}-6\varepsilon_n-(\kappa_n+3\varepsilon_n)b\bigr).
\end{equation}
The endpoint values of this affine function can be checked exactly.

\begin{lem}\label{ext:lem:endpoints}
For every integer $n\geq163$, $\Phi_n(0)>0$ and $\Phi_n(1)>0$.
\end{lem}
\begin{proof}
Let $H(n)=9n^3-46n^2+11n+2$. Then
\[
1-3\varepsilon_n=\frac{H(n)}{(n-1)(n-4)(9n-1)}.
\]
Exact simplification of \eqref{ext:eq:necessary} gives, for even $n$,
\begin{equation}\label{ext:eq:even}
\Phi_n(0)=\frac{P_e}{(n-2)H},\qquad
\Phi_n(1)=\frac{(2n^2+67n-138)Q_e}{8n(n-2)(n-1)H},
\end{equation}
where
\begin{align*}
P_e&=9n^4-1486n^3+4335n^2-2994n+256,\\
Q_e&=5n^4-457n^3+165n^2-11n+10.
\end{align*}
For odd $n$,
\begin{equation}\label{ext:eq:odd}
\Phi_n(0)=\frac{P_o}{(n-3)(n-2)H},\qquad
\Phi_n(1)=\frac{(2n^2+67n-138)Q_o}{8n(n-3)(n-2)(n-1)H},
\end{equation}
where
\begin{align*}
P_o&=9n^5-1513n^4+8901n^3-12393n^2+1384n+60,\\
Q_o&=5n^5-472n^4+1968n^3-986n^2+91n-30.
\end{align*}
Writing $n=163+t$ gives
\begin{align*}
H&=9t^3+4355t^2+702378t+37756344,\\
P_e&=9t^4+4382t^3+712407t^2+38872506t+32404656,\\
P_o&=9t^5+5822t^4+1413635t^3+152914044t^2\\
&\qquad+6261781674t+5746994076,\\
Q_e&=5t^4+2803t^3+573762t^2+50242620t+1554789528,\\
Q_o&=5t^5+3603t^4+1022674t^3+142255308t^2\\
&\qquad+9627885720t+250624461888.
\end{align*}
Every displayed coefficient is positive. All other denominator factors and
$2n^2+67n-138$ are positive for $n\geq163$. This proves the lemma.
\end{proof}

By Lemma~\ref{ext:lem:endpoints}, both endpoints are positive for
$n\geq163$. Since $\Phi_n$ is affine in $b$, we obtain
\[
\Phi_n(b)=(1-b)\Phi_n(0)+b\Phi_n(1)>0\qquad(0\leq b\leq1).
\]
This contradicts \eqref{ext:eq:necessary} for every possible $b$,
including $b=0$, when no exceptional 3-faces occur. Hence $n\leq162$,
proving Theorem~\ref{edge:thm:dimension-bound}.
Only the weight \eqref{ext:eq:weight} is used. No separate upper bound on
$b$, such as the four-dimensional estimate \eqref{eq:bad-proportion},
is needed in this final argument.

\part{Extending Vinberg's planar-angle method}\label{part:vertex}
We extend Vinberg's facewise--vertexwise comparison to finite-volume
polytopes and prove Theorem~\ref{thm:dimension-bound}. We now assign
weights to incidences of a 2-face with one of its vertices. The notation
for these weights is local to this part; the common face numbers and
product-link identities retain their meanings. The quantities
$u_{n,2}$ and $\Theta_n^{(2)}$ are the planar counterparts of
$u_{n,3}$ and $\Theta_n^{(3)}$ in Part~\ref{part:edge}.
% Confidential local draft: the two-face/vertex branch of the combined paper.
% This fragment neither starts a document nor changes appendix numbering.
\section{The facewise--vertexwise comparison}\label{two:sec:comparison}

The same total weight can be summed over 2-faces or over vertices.
At an ideal vertex, a marked pair in one simplex factor of the link
represents a distinct planar angle for every choice of a vertex in
the other factors. These multiplicities can make the weight per
incident edge large. The truncation identity supplies a cusp term
that we subtract from the vertex total before dividing by the number
of incident edges. The comparison below explains why it is this
difference that must be bounded.

Throughout this part, $F$ denotes a 2-face and the notation for weights
is local to incidences $(F,v)$, with $v$ a finite or ideal vertex of $F$.
We may assume $n\geq4$. The numbers $a_i$, $a_k^i$, and in particular
$a_2^0$, retain their definitions from the common preliminaries.

\subsection{The planar ideal-vertex term}\label{two:subsec:retained-cusp}
Recall that $V_u=f_0(L_u)$ counts the edges incident to an ideal
vertex, as in \eqref{two:eq:link-vertices}. We use the endpoint and corner
incidence identities \eqref{two:eq:endpoint-incidence} and
\eqref{two:eq:corner-incidences} from the common preliminaries.

The truncation identity \eqref{ext:eq:asymmetry} expresses
$h_{n-j}(P)-h_j(P)$ as $\sum_u g_j(L_u)$. We retain these
link contributions when deriving the planar face-number inequality;
no nonnegativity of the lower-half $h$-numbers of $P$ is needed.

The planar-angle coefficient is
\begin{equation}\label{two:eq:cn}
c_{n,2}:=\frac{u_{n,2}}{n-1}
=\begin{cases}
4/(n-2),&n\text{ even},\\[2pt]
4n/(n-1)^2,&n\text{ odd}.
\end{cases}
\end{equation}
For a link of dimension $n-1$, put
\begin{equation}\label{two:eq:Xi}
\Theta_n^{(2)}(L):=
\sum_{j=0}^{\lfloor(n-1)/2\rfloor}
g_j(L)\left(j-c_{n,2}\binom j2\right),
\qquad h_{-1}(L)=0,\quad g_0(L)=1.
\end{equation}
Its degree-zero term is zero.
The contribution also vanishes for a simplex link.
Subsection~\ref{two:sec:cusp-h} will express this quantity
directly in terms of the simplex-factor dimensions.

\begin{lem}\label{two:lem:retained}
The face numbers satisfy
\begin{equation}\label{two:eq:retained}
c_{n,2}a_2-a_1\geq\sum_{u\text{ ideal}}\Theta_n^{(2)}(L_u).
\end{equation}
Equivalently,
\begin{equation}\label{two:eq:corrected-average}
a_2^0\leq u_{n,2}-\frac{n-1}{a_2}
\sum_{u\text{ ideal}}\Theta_n^{(2)}(L_u).
\end{equation}
\end{lem}
\begin{proof}
Set $A(j)=c_{n,2}\binom j2-j$.
The inverse transformation \eqref{ext:eq:h-inverse} gives
$c_{n,2}a_2-a_1=\sum_{j=0}^nA(j)h_j(P)$.
Pair the terms of degrees $j$ and $n-j$, with
$1\leq j\leq\lfloor(n-1)/2\rfloor$.
Using \eqref{ext:eq:asymmetry}, their sum becomes
\[
[A(j)+A(n-j)]h_{n-j}(P)
-A(j)\sum_{u\text{ ideal}}g_j(L_u).
\]
For even and odd $n$, respectively, the first coefficient is
\[
\frac{(2j-n)^2}{n-2},\qquad
\frac{n((2j-n)^2-1)}{(n-1)^2}.
\]
Both are nonnegative in the paired range, and
$h_{n-j}(P)\geq0$ by \eqref{ext:eq:h-upper-nonnegative}.
If $n$ is even, the unpaired middle coefficient
$A(n/2)$ is zero.
The degree-zero coefficient is zero, while the degree-$n$
term is positive because $h_n(P)=1$ and
\[
A(n)=
\begin{cases}
n^2/(n-2),&n\text{ even},\\
n(n+1)/(n-1),&n\text{ odd}.
\end{cases}
\]
Discarding these nonnegative terms leaves
$\sum_u\Theta_n^{(2)}(L_u)$, proving \eqref{two:eq:retained}.
Finally substitute $a_1/a_2=a_2^0/(n-1)$ from
\eqref{two:eq:edge-face-incidence}.
\end{proof}

Every coefficient defining $\Theta_n^{(2)}$ is nonnegative in the stated
range, and $g_j(L_u)\geq0$ by the product-link decomposition.
Thus \eqref{two:eq:corrected-average} strengthens the usual
Nikulin--Khovanskii bound. The positive degree-$n$ term makes
the inequality strict, although the weak form suffices below.

\subsection{The comparison criterion}
The weights are combinatorial: even when $v$ is ideal and the
geometric angle is zero, its diagram may receive positive weight.
For a nonnegative rule $w$, write
\[
\sigma(F):=\sum_{v\in F}w(F,v).
\]
We shall construct $w$ and an auxiliary nonnegative rule $q$
for which
\begin{equation}\label{two:eq:face-targets}
\sigma(F)\geq5-f_0(F),\qquad
\sum_{v\in F}q(F,v)\geq1\quad\text{if }f_0(F)=3.
\end{equation}
Thus triangles require at least two units of $w$,
quadrilaterals require one, and polygons with at least
five vertices require only nonnegativity.

Define the global totals and the averages per 2-face by
\begin{equation}\label{two:eq:total-weights}
\begin{aligned}
W&:=\sum_{\dim F=2}\sigma(F),
&a_w&:=\frac{W}{a_2},\\
Z&:=\sum_{\dim F=2}\sum_{v\in F}q(F,v),
&a_q&:=\frac{Z}{a_2},\\
J&:=W+\tfrac14Z,
&\frac{J}{a_2}&=a_w+\tfrac14a_q.
\end{aligned}
\end{equation}
Here $a_w$ is the \emph{average total weight of a 2-face},
and $a_q$ is the corresponding average for the auxiliary rule.
The first target in \eqref{two:eq:face-targets} gives
\[
a_w\geq5-a_2^0>
\begin{cases}
1-\dfrac4{n-2},&n\text{ even},\\[4pt]
1-\dfrac4{n-1},&n\text{ odd},
\end{cases}
\]
by \eqref{two:eq:average-two-face-bound}. Thus the constant five
leaves a lower bound tending to one, since the upper bound for the
average number of vertices of a 2-face tends to four.

The weight of a fixed vertex is
\[
J_v:=\sum_{\substack{\dim F=2\\v\in F}}
\left(w(F,v)+\tfrac14q(F,v)\right).
\]
Counting each incidence $(F,v)$ once in either order gives
$J=\sum_vJ_v$.
We use \emph{facewise} for sums on a fixed 2-face and
\emph{vertexwise} for sums at a fixed vertex.

The auxiliary rule improves the triangle contribution without changing
the required quadrilateral contribution. The triangle estimate and
nonnegativity give
$\sum_{v\in F}q(F,v)\geq4-f_0(F)$ for every $F$. Thus, for any
$\beta\geq0$, the local lower bound for $w+\beta q$ is
\[
(5-f_0(F))+\beta(4-f_0(F))
=(5+4\beta)-(1+\beta)f_0(F).
\]
The choice $\beta=1/4$ gives the constants six and $5/4$ below.
Its role in the vertex estimates is explained after the comparison
criterion.

\begin{lem}\label{two:lem:global-lower}
The facewise estimates \eqref{two:eq:face-targets} imply
\begin{equation}\label{two:eq:global-lower}
a_w+\tfrac14a_q\geq6-\tfrac54a_2^0,
\qquad
J\geq6a_2-\tfrac54(n-1)a_1.
\end{equation}
\end{lem}
\begin{proof}
Sum the preceding local inequality with $\beta=1/4$ over all
2-faces. Dividing by $a_2$ gives the first assertion; the identity
$\sum_F f_0(F)=(n-1)a_1$ gives the second.
\end{proof}

A finite vertex has $n$ incident edges, whereas an ideal vertex $u$
has $V_u$. We use these numbers to normalize the vertex contributions.

\begin{prop}\label{two:prop:criterion}
Suppose \eqref{two:eq:face-targets} holds.
If a constant $C_n$ satisfies
\begin{equation}\label{two:eq:vertex-targets}
J_v\leq nC_n\quad\text{at every finite vertex},\qquad
J_u-\frac6{c_{n,2}}\Theta_n^{(2)}(L_u)\leq V_uC_n
\quad\text{at every ideal vertex},
\end{equation}
then the existence of $P$ requires
\begin{equation}\label{two:eq:threshold}
C_n\geq\tau_n:=
\frac12\left(\frac6{c_{n,2}}-\frac54(n-1)\right)
=
\begin{cases}
(n-7)/8,&n\text{ even},\\[2pt]
(n-7)/8+3/(4n),&n\text{ odd}.
\end{cases}
\end{equation}
\end{prop}
\begin{proof}
Put $D=\sum_{u\text{ ideal}}\Theta_n^{(2)}(L_u)$ and
$K=J-6D/c_{n,2}$. The facewise lower bound and
$D\leq c_{n,2}a_2-a_1$ give
\[
K\geq6a_2-\frac54(n-1)a_1
       -\frac6{c_{n,2}}(c_{n,2}a_2-a_1)
=\left(\frac6{c_{n,2}}-\frac54(n-1)\right)a_1.
\]
On the other hand, summing \eqref{two:eq:vertex-targets} and using
the edge-endpoint identity gives
\[
K\leq C_n\left(nV_{\mathrm f}+\sum_{u\text{ ideal}}V_u\right)=2C_na_1.
\]
Comparison of these two bounds proves \eqref{two:eq:threshold}.
\end{proof}

The subtraction thus comes from the face-number inequality, not from
a change to the angle weights. It is applied once to the combined
total $J_u$, and the resulting difference is divided by $V_u$.

The coefficient $1/4$ balances the two vertex estimates. With
$W_u=\sum_{F\ni u}w(F,u)$ and $Q_u=\sum_{F\ni u}q(F,u)$,
the corresponding ideal contribution for a general coefficient
$\beta$ would be
\[
\frac{W_u-5\Theta_n^{(2)}(L_u)/c_{n,2}}{V_u}
+\beta\frac{Q_u-4\Theta_n^{(2)}(L_u)/c_{n,2}}{V_u}.
\]
Increasing $\beta$ can only increase the finite-vertex totals, but may lower this
ideal-vertex contribution. The comparison threshold also depends on
$\beta$: it becomes
$\frac12((5+4\beta)/c_{n,2}-(1+\beta)(n-1))$.
The estimates below show that $\beta=1/4$ suffices for the bound
$n\leq81$; no optimality of this choice is claimed.

\section{Planar-angle weights and facewise estimates}\label{two:sec:weights}

The rules below preserve the required totals on triangles and
quadrilaterals while reducing the sums over marked pairs in elliptic
and affine components. Starting with a distance rule $w_0$, we assign
zero to two classes of long paths in chains and affine cycles.
The justification uses the four angle diagrams of a quadrilateral
together.

\subsection{The diagram of an angle}
Let $S_F$ be the elliptic diagram of the $n-2$ ambient facets
containing a fixed 2-face $F$. Each side of $F$ is cut out by one
further ambient facet; its node will be called a \emph{marked node}.
If the two sides meeting at a vertex $v$ have marked nodes $a,b$,
put
\[
T_{F,v}=S_F\cup\{a,b\}.
\]
By the face criterion, each side lies in exactly one ambient facet
in addition to those represented by $S_F$; hence $T_{F,v}$ has $n$ nodes. Its Gram matrix is positive
definite at a finite vertex and positive semidefinite of corank
one at an ideal vertex. In both cases,
$S_F$, $S_F\cup\{a\}$, and $S_F\cup\{b\}$ are positive definite.
We call $T_{F,v}$ the \emph{angle diagram}. Here and below, unions
of node sets always carry the induced diagram.
At a finite vertex it is the complete vertex diagram $T_v$.
At an ideal vertex it can be a proper subdiagram of the full
parabolic vertex diagram: it consists of one complete affine
factor and an elliptic node deletion in every other factor.

For a quadrilateral, write $a,b,c,d$ for the marked nodes in cyclic
order. Thus the adjacent pairs are
$\{a,b\},\{b,c\},\{c,d\},\{d,a\}$, and the opposite pairs are
$\{a,c\},\{b,d\}$. The diagrams $S_F\cup\{a,c\}$ and
$S_F\cup\{b,d\}$ are indefinite. Indeed, positive definiteness would
give a common finite vertex of the opposite sides by
Theorem~\ref{Vinberg:thm3.1}. Singular positive
semidefiniteness would give a common ideal point by Vinberg's
face correspondence \cite{Vinberg1985}. Both possibilities
contradict the disjointness of the closures of the opposite sides.

\subsection{The weight rules}\label{two:subsec:weight-rules}
We use the graph distance of Subsection~\ref{subsec:distance-weights}
and the path terminology of Subsection~\ref{ext:subsec:path-rule}.
Let $T$ be a positive semidefinite Coxeter diagram with distinct
marked nodes $a,b$, and put $k=\dist_T(a,b)$. For an angle $(F,v)$,
we take $T=T_{F,v}$: distances, degrees, and extensions are all
measured in this actual angle diagram.
Every connected component of $T$ is a tree or an affine $A$-cycle.
For a cycle, we use the shorter cyclic distance and the nonspecial
cases of the rules below.

The initial weight $w_0$ is chosen to preserve the quadrilateral
lower bound, including the long configurations in
Subsection~\ref{two:app:minimal}. The auxiliary weight $q$ is
supported on shorter paths and strengthens the triangle estimate.
Define
\begin{equation}\label{two:eq:weight-table}
\begin{aligned}
w_0(T;a,b)&=
\begin{cases}
1,&1\leq k\leq12,\\
1,&k\in\{13,14\}\text{ and the path is special},\\
1/2,&k=13\text{ and the path is nonspecial},\\
1/4,&k=14\text{ and the path is nonspecial},\\
0,&\text{otherwise},
\end{cases}\\[6pt]
q(T;a,b)&=
\begin{cases}
1,&1\leq k\leq3\text{ and the path is special},\\
1/2,&1\leq k\leq4\text{ and the path is nonspecial},\\
1/2,&4\leq k\leq7\text{ and the path is special},\\
0,&\text{otherwise}.
\end{cases}
\end{aligned}
\end{equation}
In particular, both weights vanish for disconnected pairs.

Both $w_0$ and $q$ are nondecreasing under induced positive
semidefinite extension. In a tree the old path remains, and a special
path stays special. If the extended component is an affine $A$-cycle,
the old path was nonspecial and its distance can only decrease.

For a connected pair, let $\Omega$ be its component in $T$. The final
main weight is
\begin{equation}\label{two:eq:final-weight}
w(T;a,b)=
\begin{cases}
0,&\Omega\text{ is an affine }A\text{-cycle and }k\geq9,\\
0,&\begin{aligned}[t]
   &\Omega\text{ is a tree},\ k\geq8,\text{ and }(a,b)\text{ has a}\\
   &\text{two-sided extension of length }5,
   \end{aligned}\\
w_0(T;a,b),&\text{otherwise}.
\end{cases}
\end{equation}
For disconnected pairs set $w=0$. The auxiliary weight $q$ is
unchanged. In the tree case the extended path has at least nineteen
nodes, so the classification excludes the exceptional finite and
affine diagrams. In the remaining classical types, the degree
condition and the five added simple edges at each end also force
the marked path to have only simple edges; see the proof of
Lemma~\ref{two:lem:tree-delete}.

For an incidence $(F,v)$ with marked pair $\{a,b\}$, set
\[
\begin{aligned}
w_0(F,v)&=w_0(T_{F,v};a,b),\\
w(F,v)&=w(T_{F,v};a,b),\\
q(F,v)&=q(T_{F,v};a,b).
\end{aligned}
\]
Unlike $w_0$ and $q$, the final weight $w$ need not increase under
induced extension. The two quadrilateral lemmas in
Section~\ref{two:sec:deletions} show that the zero assignments
preserve the required facewise lower bounds.

For each fixed distance between eight and fourteen, a sufficiently
long chain of simple edges has $r+O(1)$ pairs with two-sided
extensions of length five. Assigning them zero changes the coefficient
of $r$ in the total weight. Away from the ends, only distances one
through seven remain, each with weight one, giving $7r+O(1)$.
The auxiliary rule contributes $2r+O(1)$, so the combined rule
has total $15r/2+O(1)$. The constants at the ends are computed in
Section~\ref{two:sec:budgets}.

\subsection{Triangles}
For a triangle, the three marked nodes lie together in one
minimal indefinite subdiagram. The classification supplies two
angles of distance at most seven, so neither zero assignment
affects them. In addition to the two units of $w$, we need a
total of one for $q$. The following lemma establishes both
statements from the same three-marked diagrams.

\begin{lem}\label{two:lem:triangle}
Every triangular 2-face has at least two angles whose marked
nodes have distance at most seven in their angle diagrams.
Moreover,
\[
\sum_{v\in F}q(F,v)\geq1.
\]
\end{lem}
\begin{proof}
Let $a,b,c$ be the three marked nodes and put
$H=S_F\cup\{a,b,c\}$. This diagram has $n+1$ nodes, so its
Gram matrix cannot be positive definite
in $\mathbb E^{n,1}$. If it were positive semidefinite, it would
be singular, giving a common ideal point of all three sides
by Vinberg's face correspondence \cite{Vinberg1985}. This is impossible for a triangle.
Thus $H$ is indefinite. Choose an induced subdiagram $L\subseteq H$
minimal among its indefinite induced subdiagrams. It contains all
three marks, since a subdiagram omitting
one of them lies in a positive semidefinite angle diagram.
Such a minimal diagram has order between three and ten.

If $L$ has order $r\leq5$, take a shortest path from $a$ to $b$.
If it passes through
$c$, its $a$--$c$ and $c$--$b$ segments avoid the remaining
mark. Otherwise, join $c$ to this path, stopping at the first
node met. The $a$--$b$ path and the route from $c$ to one
of its ends give two distinct marked pairs whose paths avoid
the third mark. Each path lies in the corresponding angle
diagram and omits at least that third node, so it has at most
$r-2\leq3$ edges. Both pairs therefore receive auxiliary
weight at least $1/2$.

In orders $6,\ldots,10$, the classification has respectively
$12,3,4,4,3$ minimal indefinite diagrams. For each such diagram
and each three-element marked subset, delete each mark in turn
and compute the distance and the special-path condition for
the other two. All
\[
12\binom63+3\binom73+4\binom83+4\binom93+3\binom{10}3=1265
\]
marked subsets satisfy both assertions. The diagram list and
the exact finite procedure are given in
Subsection~\ref{two:app:minimal}.

The actual angle diagrams are induced positive semidefinite
extensions of these deletion diagrams. Monotonicity gives the
assertion for $q$, and the two paths of length at most seven
remain available. Neither of the zero assignments applies
at these distances, so both angles have final weight one.
\end{proof}

\subsection{\texorpdfstring{The quadrilateral estimate for $w_0$}{The quadrilateral estimate for w0}}
The following classification data concern a minimal indefinite
Coxeter diagram with two marked nodes. In the quadrilateral proof
below, they are applied to subdiagrams of the opposite-pair diagrams
$S_F\cup\{a,c\}$ and $S_F\cup\{b,d\}$, with the marks in cyclic order.
\begin{lem}\label{two:lem:tails}
Let $L$ be a Coxeter diagram whose Gram matrix is indefinite and
whose proper induced subdiagrams have positive semidefinite Gram
matrices. Let $X=\{x_1,x_2\}\subseteq V(L)$ consist of two distinct
nodes such that $G_{L\setminus\{x_i\}}$ is positive definite for
$i=1,2$. For $u\in V(L)$, put
\[
d_X(u)=\min_{x\in X}\dist_L(u,x).
\]
Then:
\begin{enumerate}[label=\textup{(\roman*)},leftmargin=2.3em]
\item $d_X(u)\leq7$, and for a fixed marked pair there is at
most one node at distance seven.
\item If $d_X(u)=7$, then $\det G_{L\setminus\{u\}}=0$.
\item If $d_X(u)\geq6$ and $\det G_{L\setminus\{u\}}\ne0$, then
\[
\max_{x\in X}(G_{L\setminus\{x\}}^{-1})_{uu}\geq16/5.
\]
\item Nodes at nearest-mark distances six and seven, when both
occur, are adjacent.
\end{enumerate}
All Gram matrices are normalized to have diagonal entries one.
The data for $d_X(u)\geq6$ comprise twenty-one positions in sixteen
marked-diagram records; see Subsection~\ref{two:app:minimal}.
\end{lem}

\begin{lem}\label{two:lem:quad-preliminary}
For every quadrilateral 2-face $F$,
\[
\sum_{v\in F}w_0(F,v)\geq1.
\]
\end{lem}
\begin{proof}
Put $S=S_F$, $X=\{a,c\}$, and $Y=\{b,d\}$.
The opposite-pair diagrams $S\cup X$ and $S\cup Y$ are indefinite.
Choose $L_X\subseteq S\cup X$ and $L_Y\subseteq S\cup Y$, minimal
among the indefinite induced subdiagrams of these respective diagrams.
They contain both of their own marks and none of the other pair.
Every single-mark deletion in either diagram is positive definite,
because it lies in $S$ together with just one marked node.
If any adjacent pair has distance at most twelve, its weight
alone proves the assertion. We therefore assume that all finite
adjacent-pair distances are at least thirteen.

If $L_X$ has order two, it is the dotted edge $X$ and contains no
node of $S$. The cores are then disjoint, and a joining edge from
$X$ gives an adjacent pair at distance at most $1+7=8$. The same
applies to $L_Y$, so both cores have order at least three, and
$d_X$, $d_Y$ below are computed in $L_X$ and $L_Y$, respectively.

Suppose first that $L_X$ and $L_Y$ are disjoint. An edge must join
them; otherwise their principal Gram matrix would be block diagonal
with two indefinite blocks and would have at least two negative
directions. For a joining edge $uv$, its endpoints must satisfy
\[
d_X(u)+d_Y(v)+1\geq13.
\]
Each distance is at least five. Two joining edges would, together
with shortest paths in $L_X$ and $L_Y$, form a cycle entirely in $S$.
Indeed, a path between two joining endpoints cannot pass through
a mark: it has length at most eight, while passing through a
mark would require length at least ten. The same reasoning
applies if the joining edges share an endpoint. Since $S$ is
elliptic and hence a forest, there is exactly one joining edge.

Write its Gram entry as $-t$, where $t\geq1/2$. Then
\[
\det G_{L_X\cup L_Y}
=\det G_{L_X}\det G_{L_Y}
-t^2\det G_{L_X\setminus\{u\}}\det G_{L_Y\setminus\{v\}}.
\]
Both core determinants are negative. If either cofactor vanishes,
the displayed determinant is positive. The matrix is then
nonsingular and already has a negative direction, so it has
at least two negative eigenvalues, a contradiction.
Otherwise Lemma~\ref{two:lem:tails}\textup{(ii)} excludes distance
seven on both sides. Both nearest-mark distances are therefore six.
By part~\textup{(iii)}, choose $x\in X$ and $y\in Y$ such that
$A=G_{L_X\setminus\{x\}}$ and $B=G_{L_Y\setminus\{y\}}$ satisfy
$(A^{-1})_{uu},(B^{-1})_{vv}\geq16/5$.
Positive semidefiniteness of the corresponding adjacent-pair
deletion requires $t^2(A^{-1})_{uu}(B^{-1})_{vv}\leq1$, whereas
\[
t^2(A^{-1})_{uu}(B^{-1})_{vv}
\geq\frac14\left(\frac{16}{5}\right)^2=\frac{64}{25}>1.
\]
Thus $L_X$ and $L_Y$ intersect.

For a common node $u$,
$d_X(u)+d_Y(u)\geq13$. Each distance is at most seven, and each
side has at most one distance-seven node. Consequently the
intersection has at most two nodes. If there are two, their
distance pairs are $(7,6)$ and $(6,7)$, since each side has
only one distance-seven node. Lemma~\ref{two:lem:tails}\textup{(iv)}
therefore makes them adjacent in both subdiagrams. There are no additional
edges between $L_X\setminus L_Y$ and $L_Y\setminus L_X$.
The endpoints of such an edge have nearest-mark distance at
least five, while a common node has nearest-mark distance at
least six in each subdiagram. Join the endpoints to that common
node by shortest paths in $L_X$ and $L_Y$. Each path has length at
most eight and cannot pass through a mark, since that would
require at least eleven edges. The extra edge would therefore
close a cycle containing only nodes of $S$, which is impossible.
Hence the union is a one-node or one-edge gluing.

The sixteen marked records of
Subsection~\ref{two:app:minimal} give ninety-two candidate
gluings. Three have incompatible labels, sixty-five have
non-positive-definite $S$, and three have a non-positive-semidefinite
angle diagram. The remaining twenty-one all have adjacent-pair
$w_0$-weight sum at least one. The precise gluing and
checking rules are recorded in the
same subsection. The paths in these gluings remain available in
the actual angle diagrams, so monotonicity of $w_0$ proves the
inequality.
\end{proof}

\section{Zero weights and facewise lower bounds}
\label{two:sec:deletions}

To preserve the quadrilateral total, we show that each zero assignment
forces another angle to have final weight one. Otherwise the Gram
matrix of the four-marked diagram would have two negative directions,
contrary to the Lorentzian signature. The next two lemmas establish
this implication for affine cycles and for trees.

For an added marked node, its \emph{support} on a diagram
consists of the nodes joined to it and the labels of those edges.
The word \emph{coupling} below denotes an off-diagonal matrix
entry, measured in the indicated Gram-matrix normalization.

\subsection{Schur complements along paths}
Eliminating a positive definite block $A$ from a symmetric
matrix means forming its \emph{Schur complement}:
\[
\begin{pmatrix}A&B\\B^T&D\end{pmatrix}
\longmapsto D-B^TA^{-1}B.
\]
Here the superscript $T$ denotes transpose. The number of negative
eigenvalues is unchanged, by block congruence. In the path calculations below,
internal nodes have no attachments outside the path; the
endpoints may have arbitrary retained attachments.
\begin{lem}\label{two:lem:path-extension}
For a path of $k\geq1$ edges all labeled three, use the matrix $2G$ and
eliminate its internal nodes. The remaining matrix on its
endpoints is
\begin{equation}\label{two:eq:path-schur}
I_2+\frac1k
\begin{pmatrix}1&-1\\-1&1\end{pmatrix}.
\end{equation}
If the endpoints and the rest of a larger diagram are fixed,
lengthening this path from $k$ to $k'\geq k$ changes the Schur
complement by the negative semidefinite matrix
\begin{equation}\label{two:eq:path-extension}
\left(\frac1{k'}-\frac1k\right)
(\mathbf e_u-\mathbf e_v)(\mathbf e_u-\mathbf e_v)^T,
\end{equation}
where $\mathbf e_u,\mathbf e_v$ are the coordinate unit vectors
of the retained endpoints.
\end{lem}
\begin{proof}
The internal Cartan matrix is $A_{k-1}$, and
\[
(A_{k-1}^{-1})_{ij}=\min(i,j)-ij/k,
\qquad 1\leq i,j\leq k-1.
\]
Substitution in the Schur complement gives
\eqref{two:eq:path-schur}; subtracting the expressions for the
two lengths gives \eqref{two:eq:path-extension}. This difference
is negative semidefinite because $1/k'-1/k\leq0$ and the outer
product is positive semidefinite. When placing
the path in a larger diagram, only this difference is used:
the endpoint diagonal entries are not counted a second time.
\end{proof}

We will retain the four marks in the order $(a,b,c,d)$ and
eliminate a positive definite unmarked subdiagram. Denote the
resulting matrix by $R_0$. For positive rational numbers $x,y$,
put
\begin{equation}\label{two:eq:witness}
U(x,y)=
\begin{pmatrix}
1&0\\0&1\\x&0\\0&y
\end{pmatrix}.
\end{equation}
If the first diagonal entry of $U^TR_0U$ is negative and its
determinant is positive, then this two-by-two symmetric matrix
is negative definite. The two columns of $U$ therefore span
a two-dimensional subspace on which $R_0$ is negative definite,
contradicting the Lorentzian signature.

Increasing the magnitudes of the opposite-pair couplings $ac$
and $bd$ by $t,s\geq0$ replaces this restricted matrix by
\begin{equation}\label{two:eq:stronger-coupling}
U^TR_0U-\operatorname{diag}(2tx,2sy).
\end{equation}
Thus these increases preserve the two negative directions.
This calculation is used only for the two specified couplings.

We retain all marks, support nodes, and branches before eliminating
the disjoint internal paths. Lengthening these paths subtracts a
positive semidefinite matrix from the Schur complement, preserving
both negative witnesses and failures of the required definiteness
conditions.
If an unmarked block itself ceases to be positive definite, the face
criterion is already violated. This reduces arbitrary path lengths
to finitely many initial configurations.

\subsection{Long angles in an affine cycle}
\begin{lem}\label{two:lem:cycle-delete}
Suppose an angle of a quadrilateral lies in an affine $A$-cycle
component of its complete angle diagram, and the cyclic distance
of its marked pair is at least nine. Then another angle has
distance at most seven, or another angle has distance eight
and its path contains an internal node of degree at least three
in its complete angle diagram. In either case, this other angle
has final weight one.
\end{lem}
\begin{proof}
Write the original angle diagram with marks $a,b$ as
$\Omega\sqcup\Omega^\perp$, where $\Omega$ is the cycle component.
The two arcs between them have lengths $k_1,k_2\geq9$.
The components of $\Omega^\perp$ are orthogonal to $\Omega$. If $c,d$ are
adjacent in the diagram, their angle has distance one. We may
therefore assume $g_{cd}=0$.

After deleting $a$, the diagram $\Omega\setminus\{a\}$ is a long
$A$-path. The positive semidefinite classification permits the
following eight ways for the new node $c$ to attach to this path:
\begin{center}
\begin{tabular}{lc}
\toprule
Attachment to the remaining path&Number of choices\\
\midrule
No attachment&1\\
A label-three edge to either endpoint&2\\
A label-four edge to either endpoint&2\\
A label-three edge to the neighbor of either endpoint&2\\
Label-three edges to both endpoints&1\\
\bottomrule
\end{tabular}
\end{center}
The same list applies to $d$, with the roles of $a,b$ exchanged.
If $c$ has no attachment to this path, its coupling to $a$ must
be nonzero.
Indeed, $\Omega^\perp\cup\{c\}$ is positive definite, being a principal
subdiagram of $S_F\cup\{c\}$. If $g_{ac}=0$, the actual
opposite-pair diagram is
$(\Omega\setminus\{b\})\sqcup(\Omega^\perp\cup\{c\})$, a positive definite
diagram. This contradicts the opposite-pair criterion.

Use $2G$, whose diagonal entries are two and whose simple-edge
entries are minus one. A positive diagonal congruence makes
the new node attached by a label-four edge have diagonal four
and the corresponding off-diagonal entry minus two. Take the
opposite-pair magnitude to be one for empty support and zero for
the other supports; stronger magnitudes
are handled by \eqref{two:eq:stronger-coupling}.

It suffices first to consider
$(k_1,k_2)=(9,9),(9,10),(10,10)$, with $8^2=64$ support pairs
for each choice. Subsection~\ref{two:app:cycle-witnesses} gives
the rational choices of $x,y$ in \eqref{two:eq:witness}.
Among the 192 cases, 185 have a negative definite restricted
matrix. In each of the other seven cases, the angle with marks
$c,d$ has distance seven or eight. All four distance-eight paths
contain an internal branch; the finite support checks are
specified in Subsection~\ref{two:app:cycle-witnesses}. The counts and minimum
positive determinants are
\[
\begin{array}{c|ccc}
\toprule
(k_1,k_2)&(9,9)&(9,10)&(10,10)\\
\midrule
\text{Two-negative-direction witnesses}&58&63&64\\
\text{Minimum determinant}&1/81&8/1125&24/625\\
\bottomrule
\end{array}.
\]
Choose the angle with marks $c,d$ in each exceptional case.
When a new mark attaches next to an endpoint, the attachment node
retains both of its
original cyclic neighbors after $a,b$ are deleted, and is also
joined to the new mark. It therefore has degree three in this
angle diagram and lies strictly between $c,d$ on the path.
Additional original components and their attachments can only
increase its degree. Thus the complete positive semidefinite
component is a tree, and the retained path remains its unique
path between $c,d$. A distance-eight path with this internal
branch cannot have the extension required by the tree zero
assignment. The cycle zero assignment is also excluded.
Distance-seven paths are excluded from both zero assignments
by definition. Every exceptional case therefore gives another
angle of final weight one. The path in the $(9,10)$ case uses
only the length-nine arc.

For arbitrary $k_1,k_2\geq9$, retain the marks, the nodes within
two steps of each mark on each arc, and $c,d$. The two remaining
internal paths have lengths $k_1-4$ and $k_2-4$. Eliminate these
paths first, and apply Lemma~\ref{two:lem:path-extension}.
Lengthening them preserves every negative witness. It also
preserves the failures of positive definiteness or semidefiniteness
used to reject a support. The exceptional short paths do not
pass through an arc that is lengthened. If both arcs have length
nine, use $(9,9)$. If just one has length nine, exchange the
arcs if necessary and start from $(9,10)$. If both have length
at least ten, start from $(10,10)$. Thus the three initial
length pairs cover all $k_1,k_2\geq9$.

Additional orthogonal components can be omitted: the negative
directions already lie in the retained principal submatrix.
\end{proof}

\subsection{Tree angles with extensions at both ends}
\begin{lem}\label{two:lem:tree-delete}
Suppose an angle of a quadrilateral lies in a tree component,
has a two-sided extension of length five, and has distance at
least eight. Then another angle has distance at most seven
in its actual angle diagram.
\end{lem}
\begin{proof}
Let $\Omega$ be the original component and $a,b$ the marks of
the indicated angle. Again, if $c,d$ are adjacent, the conclusion
holds, so assume $g_{cd}=0$.
The marked path and the five added edges at each end together
contain at least nineteen nodes, so $\Omega$ cannot be an exceptional
finite or affine diagram. Its main chain consists entirely of
simple edges. At either endpoint, the remaining decoration is
one of three types: $A$ has no additional leaf, $B$ has one leaf
joined by a label-four edge, and $D$ has two leaves joined by
simple edges. Up to reversing the chain, the six possibilities
are
\[
AA,\ BA,\ DA,\ BB,\ BD,\ DD.
\]
The first three are finite diagrams of types $A,B,D$; the last
three are affine diagrams of types
$\widetilde C,\widetilde B,\widetilde D$.

Number the main chain $0,\ldots,s_1+k+s_2$, omitting the extra
end leaves from this numbering, and put
\[
a=s_1,\qquad b=s_1+k,\qquad s_1,s_2\geq5,\qquad k\geq8.
\]
Thus $s_1$, $k$, and $s_2$ count the simple edges from the left
endpoint to $a$, from $a$ to $b$, and from $b$ to the right
endpoint, respectively.

Write the original angle diagram as $\Omega\sqcup\Omega^\perp$,
where $\Omega$ is the component containing $a,b$ and $\Omega^\perp$
is the union of the other components. Since $\Omega^\perp\subset S_F$,
it is elliptic, and $S_F=(\Omega\setminus\{a,b\})\sqcup\Omega^\perp$.
The other adjacent angles require
\[
(\Omega\setminus\{a\})\cup\Omega^\perp\cup\{c\},\qquad
(\Omega\setminus\{b\})\cup\Omega^\perp\cup\{d\}
\]
to be positive semidefinite. In particular, their restrictions
to the original chain give the support list below. The
attachments to $\Omega^\perp$ must then be checked as well, because
they can supply the indefiniteness required by an opposite
pair even when the opposite-pair coupling itself is zero.

\smallskip\noindent\emph{Allowed supports on the original chain.}
Deleting $a$ separates $\Omega$ into two finite classical
components. On its near side, the nodes $0,\ldots,s_1-1$,
together with the left end decoration, form $A_{s_1}$, $B_{s_1+1}$,
or $D_{s_1+2}$. On the other side, replace $s_1$ by $k+s_2$
and use the right decoration. Deleting $b$ gives the analogous
decomposition with near length $s_2$ and far length $s_1+k$.
This explains the lower orders $5,6,7$ in the following table.
The complete elliptic and parabolic classifications give the
following nonempty supports of one new node on each component:
\begin{center}
\begin{tabular}{c|p{.75\linewidth}}
\toprule
Component&Allowed support\\
\midrule
$A_r$, $r\geq5$&
A label-three or label-four edge at an endpoint; a label-three
edge at a next-to-endpoint node; simple edges at both endpoints;
or a simple leaf giving arm lengths
$(1,2,2)$, $(1,2,3)$, $(1,2,4)$, $(1,2,5)$, $(1,3,3)$.\\
$B_r$, $r\geq6$&
A label-three or label-four edge at the ordinary endpoint,
or a label-three edge at its neighbor.\\
$D_r$, $r\geq7$&
A label-three or label-four edge at the long-arm endpoint,
or a label-three edge at its neighbor. For $r=7,8$, a simple
edge at one short leaf is also possible.\\
\bottomrule
\end{tabular}
\end{center}
Empty support is allowed. A new node may also attach by one
simple edge to each of the two components: at an endpoint of
an $A$-component, the ordinary endpoint of a $B$-component,
or the long-arm endpoint of a $D$-component. The short-leaf
exceptions in $D_7,D_8$ cannot be used to join two such
components, since the resulting branch would have two long arms.
If both attachments lie on the same component and form a
cycle, the whole connected component must be an affine
$A$-cycle. The classification excludes every other repeated
support. The lower bounds on the component orders exclude
additional labels five and six and additional $F$- or $H$-type
possibilities.

\smallskip\noindent\emph{Attachments to orthogonal components.}
Temporarily set $g_{ac}=0$. If $\Omega\cup\{c\}$ is already
indefinite, we allow this zero coupling in the finite test, enlarging
the set of cases to be excluded. Otherwise, we must determine whether
attachments to $\Omega^\perp$ can make $\Omega\cup\Omega^\perp\cup\{c\}$
indefinite. This diagram contains the actual opposite-pair diagram
$S_F\cup\{a,c\}$, which omits $b$. If the larger diagram remains
positive semidefinite, the opposite-pair condition forces $g_{ac}\ne0$.
The same argument applies to $d$ and $g_{bd}$.

If $c$ has no attachment to $\Omega$ and $g_{ac}=0$, the diagram
$\Omega^\perp\cup\{c\}$ is positive definite because it is a principal
subdiagram of $S_F\cup\{c\}$. It is orthogonal to $\Omega$.
Empty support therefore cannot make the opposite-pair diagram
indefinite, and a nonzero coupling $ac$ is required.

If $\Omega$ is affine, any nonempty attachment already makes
$\Omega\cup\{c\}$ indefinite. If $\Omega$ is finite, the remaining
supports are a simple extension at an ordinary endpoint,
a heavy extension at an ordinary endpoint, a simple leaf
next to an ordinary endpoint, and an $AA$-attachment joining
the two original endpoints.

Suppose $c$ extends an ordinary endpoint by a simple edge.
After the near mark is deleted, the path from $c$ to the cut
end has length at least five. If an attachment from $\Omega^\perp$
introduces a branch at $c$, this is one of its arms; if the
branch lies farther into $\Omega^\perp$, that arm is longer still.
No exceptional finite diagram has an arm of this length.
The only affine exceptional possibility has arm lengths
$(1,2,5)$ and type $\widetilde E_8$. It lies entirely in
$S_F\cup\{c\}$, contradicting its positive definiteness.
The remaining attachments extend the chain or give a classical
$B$- or $D$-end decoration. Restoring the deleted marks gives
one of the six finite or affine classical types above, still
positive semidefinite. Thus such an attachment cannot provide
the indefiniteness required by the opposite pair when $g_{ac}=0$.
A heavy extension followed by an attachment in $\Omega^\perp$ puts the
heavy edge in the interior of a path beyond the lengths
allowed for $F_4$ or $\widetilde F_4$, so the resulting diagram is not positive
semidefinite.
If the new node joins the two original $A$-endpoints, deleting
a mark gives a long path through the new node. Its two arms
are too long to permit a further attachment from $\Omega^\perp$, since
that attachment would create a branch excluded by the elliptic
and parabolic classifications. Without such an attachment,
restoring the mark closes an affine $A$-cycle, so the
opposite-pair diagram is not indefinite.

The remaining exception is a simple leaf $c$ attached next
to an ordinary endpoint. Number that endpoint zero and
suppose the near mark is $a=s_1$. After deleting $a$, the
branch has arm lengths $(1,1,s_1-2)$. An attachment from
$\Omega^\perp$ can continue the arm through $c$ only as a simple
path. Adding $p\geq1$ edges changes these lengths to
$(1,1+p,s_1-2)$. The only positive definite possibilities
are $p=1$, $s_1=5$ and $p=1$, $s_1=6$, giving $E_7$ and
$E_8$. The semidefinite boundary cases are $(p,s_1)=(2,5)$ and
$(1,7)$, with arm lengths $(1,3,3)$ and $(1,2,5)$, of types
$\widetilde E_7$ and $\widetilde E_8$. These branches use
only nodes before the near mark and the external arm through
$c$, so they remain in $S_F\cup\{c\}$ and violate its
positive definiteness. One can check these boundaries directly:
eliminating the three arms in $2G$ leaves the central entry
\[
2-\frac12-\frac{p+1}{p+2}-\frac{s_1-2}{s_1-1}
=\frac12+\frac1{p+2}+\frac1{s_1-1}-1.
\]
For $s_1\geq5$ and $p\geq1$, it is positive exactly in the
two finite cases above, zero in the two boundary cases, and
negative otherwise.
Larger extensions, a second branch, or a heavy edge are
excluded by the classification. At the far ordinary
endpoint, the corresponding long arm has length at least
$k+s_2-2\geq11$, so neither exceptional finite case occurs.

Thus, except for the external $E_7,E_8$ attachments, a support for
which $\Omega\cup\{c\}$ remains positive semidefinite requires a nonzero
coupling $ac$. We test its minimum magnitude one in the rational
normalization; \eqref{two:eq:stronger-coupling} covers larger
magnitudes, including those arising from diagonal congruence at
a heavy node. Supports making the larger diagram indefinite are
tested with zero. In an $E_7$ or $E_8$ case, retain one external node
joined to $c$ by a simple edge and allow $g_{ac}=0$.
After omitting the other external nodes, the retained node
is a leaf in the principal diagram. Its diagonal entry in
$2G$ is two and its coupling to $c$ is minus one.
Eliminating it changes the diagonal entry at $c$ from two
to $2-1/2=3/2$. Treat $d$ in the same way when needed.
If the selected external nodes coincide, are adjacent, or
one is adjacent to the other new mark, they give an actual
path between $c,d$ of length at most three. Otherwise they
can be eliminated independently. A two-dimensional negative
definite subspace in the retained principal diagram already
contradicts the Gram matrix of $P$, so no further external
nodes need to enter the calculation.

\smallskip\noindent\emph{The finite rational calculation.}
Put $S_0=\Omega\setminus\{a,b\}$, which is positive definite,
and eliminate it to obtain the four-by-four matrix $R_0$.
Take
\[
s_1,s_2\in\{5,6,7,8,9\},\qquad
k\in\{8,9,10\}.
\]
The six end types give 450 initial chains.
For each allowed support pair, first reject a failure of
positive definiteness of $S_F\cup\{c\}$ or $S_F\cup\{d\}$,
or a failure of positive semidefiniteness of an adjacent-pair
angle diagram. Next retain any alternative angle with a
path of length at most seven. For each remaining case,
use \eqref{two:eq:witness} with $x,y$ selected from
\[
\left\{
\frac13,\frac25,\frac12,\frac35,\frac23,\frac34,\frac45,
1,\frac54,\frac32,2
\right\}.
\]
Every support record yields one of three certificates: a failed
definiteness condition, a short alternative path, or a rational
two-dimensional negative restriction.
Subsection~\ref{two:app:tree-records} records the finite census,
certificate format, and a worked example.

\smallskip\noindent\emph{Arbitrary lengths.}
Retain the end decorations, the nodes within two
steps of the marks and the original endpoints, the new nodes,
and the necessary leaves in $\Omega^\perp$. The remaining internal
paths are disjoint. Once $s_1,s_2\geq9$, the support list acquires
no further exceptional cases. Increasing the tails or the
distance between the marks amounts to lengthening these
internal paths. Lemma~\ref{two:lem:path-extension} preserves
the negative witnesses and the recorded failures of
positive definiteness or semidefiniteness.
Every retained short-angle record also specifies that its
path avoids the internal paths that may be lengthened.
Keep $s_1$ and $s_2$ unchanged when they are at most eight,
and otherwise start from nine. Keep $k=8$ or $9$ unchanged,
and otherwise start from ten. Every larger configuration with
the same retained support is obtained by lengthening the
unmarked internal paths. Finally,
\eqref{two:eq:stronger-coupling} covers all stronger
opposite-pair couplings. Thus the finite initial lengths
cover all $s_1,s_2\geq5$ and $k\geq8$.
\end{proof}

\subsection{The final lower bounds}
\begin{prop}\label{two:prop:facewise}
The final weights satisfy
\begin{equation}\label{two:eq:facewise-target}
\sum_{v\in F}w(F,v)\geq5-f_0(F)
\end{equation}
for every 2-face $F$. If $F$ is a triangle, then
\[
\sum_{v\in F}q(F,v)\geq1.
\]
\end{prop}
\begin{proof}
For a triangle, Lemma~\ref{two:lem:triangle} gives two
angles of distance at most seven, both having final weight
one, and auxiliary weight sum at least one.
For a quadrilateral on which no weight is changed,
Lemma~\ref{two:lem:quad-preliminary} applies.
If a tree angle is assigned zero,
Lemma~\ref{two:lem:tree-delete} gives another angle of
distance at most seven, whose final weight is one.
If an affine-cycle angle is assigned zero,
Lemma~\ref{two:lem:cycle-delete} directly gives another
angle of final weight one. Thus the two zero assignments
preserve the quadrilateral lower bound.
For $f_0(F)\geq5$, nonnegativity gives
\eqref{two:eq:facewise-target}.
\end{proof}

\section{Weight sums at vertices and retained cusp terms}
\label{two:sec:budgets}
We estimate $J_v/n$ at finite vertices and
$[J_u-(6/c_{n,2})\Theta_n^{(2)}(L_u)]/V_u$ at ideal vertices,
as required by \eqref{two:eq:vertex-targets}. The latter estimate
has three steps: count the planar angles represented by each marked
pair, subtract the cusp contribution, and bound the resulting ratio
by a function of the simplex-factor dimensions. Concavity then reduces
the optimization of this upper bound to endpoint comparisons.

\subsection{Weight sums and angle multiplicities}
For a positive semidefinite diagram $T$, define the total weights
by summing the rules of Subsection~\ref{two:subsec:weight-rules}
over its unordered pairs of distinct nodes:
\begin{equation}\label{two:eq:diagram-budgets}
\begin{aligned}
W(T)&=\sum_{\substack{\{a,b\}\subset V(T)\\a\ne b}}w(T;a,b),\\
W_q(T)&=\sum_{\substack{\{a,b\}\subset V(T)\\a\ne b}}q(T;a,b),\\
J(T)&=W(T)+\frac14W_q(T).
\end{aligned}
\end{equation}
Pairs in different connected components have weight zero,
so all three quantities are additive over the components.
At a finite vertex $v$, the complete vertex diagram is
elliptic of order $n$. Every pair of its nodes determines
a 2-face through $v$, and the angle diagram is the complete
vertex diagram. Hence its contribution is exactly $J(T_v)$.

At an ideal vertex, an edge of the product link varies in one factor
and fixes a vertex in each other factor. Different choices give
distinct edges of the link, hence distinct planar angles at $u$,
although the marked pair in the varying factor is unchanged.

\begin{lem}\label{two:lem:cusp-additivity}
Let the complete diagram at an ideal vertex $u$ be
$\Sigma_1\sqcup\cdots\sqcup\Sigma_r$, where $\Sigma_i$
is affine of rank $d_i$, and write
$V_u=\prod_i(d_i+1)$ for the number of link vertices.
Then
\begin{equation}\label{two:eq:cusp-additivity}
\frac{J_u}{V_u}
=\sum_{i=1}^r\frac{J(\Sigma_i)}{d_i+1}.
\end{equation}
\end{lem}
\begin{proof}
A 2-face through $u$ corresponds to an edge of the product
link $L_u=\prod_i\Delta_i$. Such an edge varies in one factor
and fixes a vertex in every other factor.
If it varies in factor $i$, its two marks lie in $\Sigma_i$.
The complete angle diagram contains all of $\Sigma_i$ and
one elliptic node deletion from each other factor, orthogonal
to $\Sigma_i$. Its weight is therefore exactly the weight
of that marked pair in the complete $\Sigma_i$.
Thus each pair represents $\prod_{j\ne i}(d_j+1)=V_u/(d_i+1)$
distinct angles, and
\[
J_u=\sum_i\frac{V_u}{d_i+1}J(\Sigma_i).
\]
Dividing by $V_u$ gives \eqref{two:eq:cusp-additivity}.
The marked component is exactly $\Sigma_i$, so no monotonicity of
the final weight is needed.
\end{proof}

For a product of $n-1$ intervals, each bold pair has $w=q=1$
and is repeated $2^{n-2}$ times. Hence
$J_u/V_u=5(n-1)/8$: normalization by the number of incident
edges alone does not give a dimension-independent upper bound.
The retained cusp term is essential.

\subsection{Bounds for elliptic and affine diagrams}
In this subsection $r$ denotes the number of nodes.
An affine diagram with $r$ nodes has rank $r-1$.
For $r\geq28$, direct counting gives
\begin{equation}\label{two:eq:stable-budgets}
\begin{array}{c|ccc}
\toprule
T&W(T)&W_q(T)&J(T)\\
\midrule
A_r&7r+59/2&2r-5&15r/2+113/4\\
B_r&7r+73/2&2r-2&15r/2+36\\
D_r&7r+183/4&2r+5&15r/2+47\\
\widetilde C_{r-1}&7r+87/2&2r+1&15r/2+175/4\\
\widetilde B_{r-1}&7r+211/4&2r+8&15r/2+219/4\\
\widetilde D_{r-1}&7r+62&2r+15&15r/2+263/4\\
\bottomrule
\end{array}.
\end{equation}

Before the zero assignments, the sum of the main weights $w_0$ is
$51r/4-b_T$, where
\[
\begin{array}{c|rrrrrr}
\toprule
T&A&B&D&\widetilde C&\widetilde B&\widetilde D\\
\midrule
b_T&88&347/4&333/4&171/2&82&157/2\\
k_T&0&1&2&2&3&4\\
\bottomrule
\end{array}.
\]
Here $k_T$ is the number of decorated nodes outside the
main chain. The number of distance-$s$ pairs having a
two-sided extension of length five is
\[
\max\{0,r-k_T-10-s\}.
\]
Subtract this number with coefficients
$1,1,1,1,1,1/2,1/4$ for $s=8,\ldots,14$, respectively.
This gives the $W(T)$ column. The auxiliary sum for a chain of simple edges is $2r-5$; the fixed corrections at heavy or branched
ends give the $W_q(T)$ column.
For $r\geq28$, every count lies in its linear range.
The end corrections involve only bounded distances, so
the same formulas hold for every larger order.

For an affine cycle with $d\geq2$, cyclic pairs are retained
only at distances at most eight. Consequently
\begin{equation}\label{two:eq:cycle-budget}
\frac{J(\widetilde A_d)}{d+1}
=\min(d/2,8)+\frac14\min(d/4,2).
\end{equation}
The diagram $\widetilde A_1$ is treated separately and has
normalized cost $5/8$.

For small classical diagrams, the sums are obtained directly
from the definition. Every connected finite exceptional
diagram has at most eight nodes and is covered by
\[
J(T)\leq\frac54\binom r2.
\]
The exceptional affine data are
\begin{equation}\label{two:eq:exceptional-budgets}
\begin{array}{c|ccccc}
\toprule
\Sigma&\widetilde E_6&\widetilde E_7&\widetilde E_8&
\widetilde F_4&\widetilde G_2\\
\midrule
J(\Sigma)/|V(\Sigma)|&
51/14&133/32&337/72&19/8&29/24\\
\bottomrule
\end{array}.
\end{equation}

For one simplex factor of dimension $d$, the normalized weight
contribution is $J(\Sigma)/(d+1)$, where $\Sigma$ is its affine
diagram. The following function bounds this contribution in terms
of $d$ alone.
\begin{defn}\label{two:def:factor-budget}
For a nonnegative integer $d$, define
\[
\psi(d)=
\begin{cases}
\dfrac{5d}{8},&0\leq d\leq5,\\[1ex]
\displaystyle\max_{\Sigma\in\{\widetilde A_d,\widetilde B_d,
\widetilde C_d,\widetilde D_d\}}
\frac{J(\Sigma)}{d+1},&d\geq6.
\end{cases}
\]
\end{defn}
The exceptional affine values in
\eqref{two:eq:exceptional-budgets} are also bounded by
$\psi$. In particular, every ideal vertex satisfies
\[
\frac{J_u}{V_u}\leq\sum_i\psi(d_i).
\]
For $d\geq27$, equations
\eqref{two:eq:stable-budgets} and \eqref{two:eq:cycle-budget}
give
\begin{equation}\label{two:eq:stable-factor-budget}
\psi(d)=\max\left\{\frac{17}{2},
\frac{15}{2}+\frac{263}{4(d+1)}\right\}.
\end{equation}
These factor costs will be combined with the retained
ideal-vertex terms in the face-number inequality.

\begin{lem}\label{two:lem:finite-vertex}
At every finite vertex, $J_v/n<83/8$.
\end{lem}
\begin{proof}
For a connected component of order $r\leq8$, the pair bound gives
$J(T)/r\leq5(r-1)/8\leq35/8$.
For $9\leq r\leq27$, the finite pair sums specified in
Subsection~\ref{two:app:budgets} give $J(T)/r<83/8$. For $r\geq28$, the largest value in
\eqref{two:eq:stable-budgets} satisfies
\[
\frac{J(T)}r\leq\frac{15}{2}+\frac{47}{r}
\leq\frac{15}{2}+\frac{47}{28}<\frac{83}{8}.
\]
The ratio for the whole elliptic vertex diagram is the average
of these component ratios weighted by their orders. Their orders
sum to $n$, proving the result.
\end{proof}

\subsection{Normalization of the retained cusp term}\label{two:sec:cusp-h}

We express the retained cusp term using the weighted average $m$
defined in Subsection~\ref{common:subsec:product-links}.

\begin{lem}\label{two:lem:Xi-mean}
If $L$ has dimension $n-1$, then
\begin{equation}\label{two:eq:Xi-mean}
\frac{6\Theta_n^{(2)}(L)}{c_{n,2}V}\geq\frac32m(d_1,\ldots,d_r),
\end{equation}
with equality for even $n$.
\end{lem}
\begin{proof}
For even $n$, direct simplification gives
\[
\frac1{c_{n,2}}\left(j-c_{n,2}\binom j2\right)=\frac{j(n-2j)}4.
\]
For odd $n$ the expression has the additional nonnegative term
$j/(4n)$. Substitute $g_j(L)/V=\pi_j/(n-2j)$ from
\eqref{ext:eq:normalized-g} into \eqref{two:eq:Xi}.
\end{proof}

Define the function of the simplex-factor dimensions
\begin{equation}\label{two:eq:net-cost}
\begin{aligned}
\mathcal F(d_1,\ldots,d_r)
&:=\sum_i\psi(d_i)-\frac32m(d_1,\ldots,d_r),\\
\frac{J_u}{V_u}-\frac{6\Theta_n^{(2)}(L_u)}{c_{n,2}V_u}
&\leq\mathcal F(d_1,\ldots,d_r).
\end{aligned}
\end{equation}
The inequality follows from the factor bound and
Lemma~\ref{two:lem:Xi-mean}: the sum of $\psi(d_i)$ bounds the
normalized angle weight, and $3m/2$ is a lower bound for the
normalized cusp term being subtracted. Thus $\mathcal F$ bounds
the remaining vertex contribution using only the dimensions,
independently of the affine diagram types. For an actual cusp,
these dimensions satisfy $\sum_i d_i=n-1$.

\subsection{Concavity in one factor dimension}\label{two:sec:concavity}

For the concavity arguments, allow any finite list of nonnegative
factor dimensions and put $R=\sum_i d_i$ and $V=\prod_i(d_i+1)$.
The polynomial definitions and identities of
Subsection~\ref{common:subsec:product-links} apply with $n=R+1$,
since their proofs use only the product formula for $h_L$.
Zero-dimensional factors may be omitted; for the empty product,
$h_L(t)=V=1$ and $m=0$.

We use concavity of $m$ to bound
$\mathcal F=\sum_i\psi(d_i)-3m/2$. Wherever an upper bound for
$\psi$ is linear, subtracting the concave term $3m/2$ gives a convex
sequence whose maximum occurs at an endpoint.

We use this observation in two forms. Varying one dimension
gives a bound independent of the total dimension of the link.
Varying two dimensions while preserving their sum gives the
sharper bounds needed for $82\leq n\leq89$.

\begin{lem}\label{two:lem:two-factor-mean}
For integers $a,b\geq0$,
\begin{equation}\label{two:eq:two-factor-mean}
m(a,b)=f(a,b):=\frac{\min(a,b)}2
-\frac{\min(a,b)(\min(a,b)+2)}{6(\max(a,b)+1)}.
\end{equation}
\end{lem}
\begin{proof}
Suppose $a\leq b$. In \eqref{ext:eq:two-factor-intervals}, the summand
$t^k(1+t+\cdots+t^{a+b-2k})$ occurs once for each $0\leq k\leq a$
and has $a+b+1-2k$ monomials. Thus
\[
m(a,b)=\frac{\sum_{k=0}^a k(a+b+1-2k)}{(a+1)(b+1)}.
\]
Using $\sum_{k=0}^a k=a(a+1)/2$ and
$\sum_{k=0}^a k^2=a(a+1)(2a+1)/6$ gives
\eqref{two:eq:two-factor-mean}.
\end{proof}

\begin{lem}\label{two:lem:coordinate-concavity}
With the other dimensions fixed, $m(d_1,\ldots,d_r)$ is
nondecreasing and concave in each dimension. Concavity on the
integers means that its successive differences are nonincreasing.
\end{lem}
\begin{proof}
For a fixed $c\geq0$, extend $f(c,d)$ to real $d\geq0$ by
\[
f(c,d)=\begin{cases}
d/2-d(d+2)/(6(c+1)),&d\leq c,\\[2pt]
c/2-c(c+2)/(6(d+1)),&d\geq c.
\end{cases}
\]
Each part is nondecreasing and concave, and its derivative lies
between zero and $1/2$. For $c>0$, at $d=c$ the left
derivative exceeds the right derivative by $1/(6(c+1)^2)$.
For $c=0$ the function is identically zero. Hence in every
case it is nondecreasing and concave.

Let $R_{\mathrm{rest}}$ be the sum of the other factor
dimensions, and put $c_j=R_{\mathrm{rest}}-2j$.
Here $g_j$ denotes the consecutive $h$-number difference
for this remaining product. Write its polynomial in the form
\[
\sum_j g_jt^j(1+t+\cdots+t^{c_j}),
\qquad V_{\mathrm{rest}}=\sum_j(c_j+1)g_j.
\]
By \eqref{ext:eq:two-factor-intervals}, the product of
$t^j(1+t+\cdots+t^{c_j})$ and $1+t+\cdots+t^d$ decomposes into
summands whose lowest exponents have weighted average $j+f(c_j,d)$,
with weights proportional to their numbers of monomials. Therefore
\begin{equation}\label{two:eq:coordinate-mixture}
m(d_1,\ldots,d_{r-1},d)
=\sum_j\frac{(c_j+1)g_j}{V_{\mathrm{rest}}}
\bigl(j+f(c_j,d)\bigr).
\end{equation}
The coefficients are nonnegative, sum to one, and do not depend
on $d$. Taking this weighted average preserves the stated
properties. When there are no other factors, the empty product
has $R_{\mathrm{rest}}=0$, $V_{\mathrm{rest}}=g_0=1$,
and $c_0=0$, so the same formula applies.
\end{proof}

\subsection{Concavity with fixed total dimension}
\begin{lem}\label{two:lem:fixed-sum-concavity}
Fix the other dimensions and an integer $s\geq0$. Then
\[
b\longmapsto m(d_1,\ldots,d_r,s-b,b),\qquad 0\leq b\leq s,
\]
is a concave integer sequence.
\end{lem}
\begin{proof}
Work first on the integer half-range
$0\leq b\leq\lfloor s/2\rfloor$.
By \eqref{ext:eq:two-factor-intervals}, the product of the factors
$1+t+\cdots+t^{s-b}$ and $1+t+\cdots+t^b$ is the sum of
$t^k(1+t+\cdots+t^{s-2k})$ over $0\leq k\leq b$.
The normalized coefficient sums of these summands are
\begin{equation}\label{two:eq:fixed-sum-weights}
\lambda_{b,k}=\frac{s+1-2k}{(b+1)(s-b+1)},\qquad
\sum_{k=0}^b\lambda_{b,k}=1.
\end{equation}
If this product is multiplied by $1+t+\cdots+t^c$, the resulting
decomposition has weighted average of lowest exponents
\[
\sum_{k=0}^b\lambda_{b,k}t_k,\qquad
t_k=k+f(c,s-2k).
\]
Since $0\leq\partial f/\partial d\leq1/2$ and $f(c,d)$ is
concave in $d$, the sequence $t_k$ is nondecreasing and concave.
Every finite sequence with these properties can be written as
\[
t_k=t_0+\alpha k+\sum_j\beta_j\min(k,j),\qquad
\alpha,\beta_j\geq0.
\]
To see this, decompose its nonnegative, nonincreasing successive
differences into constant steps. It is therefore enough to prove
concavity in $b$ of the averages of $k$ and of $\min(k,j)$
under \eqref{two:eq:fixed-sum-weights}.

The first average is
\[
F_s(b)=\frac b2-\frac{b(b+2)}{6(s-b+1)},\qquad
F_s''(b)=-\frac{(s+1)(s+3)}{3(s-b+1)^3}<0.
\]
For the truncated average, when $b\geq j$, summing explicitly gives
\[
\sum_{k=0}^b\lambda_{b,k}\min(k,j)
=j-\frac{j(j+1)(3s-2j+5)}{6(b+1)(s-b+1)}.
\]
Its second derivative is
\[
-\frac{j(j+1)(3s-2j+5)(3b^2-3bs+s^2+s+1)}
{3(b+1)^3(s-b+1)^3}\leq0.
\]
For $b\leq j$, the truncated average equals $F_s(b)$.
At $b=j$, its previous and current values agree with the
untruncated average. If $j+1\leq\lfloor s/2\rfloor$, the next
value is the untruncated average minus
$\lambda_{j+1,j+1}\geq0$. Thus its second difference at the join
is no larger than that of $F_s$. If $j$ is already at or beyond
the end of this half-range, no join needs to be checked.

These averages are also nondecreasing on the half-range: for fixed
$\ell<b$, the cumulative weight
\[
\sum_{k=0}^{\ell}\lambda_{b,k}
=\frac{\sum_{k=0}^{\ell}(s+1-2k)}{(b+1)(s-b+1)}
\]
decreases with $b$. Summation by parts gives the assertion for any
nondecreasing $t_k$. Symmetry under $b\leftrightarrow s-b$ then
extends concavity to the full range. Indeed, if $B$ denotes the
half-range average, the central second differences are
$2[B(\ell-1)-B(\ell)]\leq0$ for $s=2\ell$, and
$B(\ell-1)-B(\ell)\leq0$ for $s=2\ell+1$; all others have
already been checked. The cases $s=0,1$ are immediate.

Finally, let $R_{\mathrm{rest}}$ be the total dimension of the
remaining factors and put $c_j=R_{\mathrm{rest}}-2j$.
Write their product polynomial as
$\sum_jg_jt^j(1+t+\cdots+t^{c_j})$, where $g_j$ denotes its
consecutive $h$-number difference, and set
$V_{\mathrm{rest}}=\sum_j(c_j+1)g_j$.
The full weighted average is
\[
\sum_j\frac{(c_j+1)g_j}{V_{\mathrm{rest}}}
\left[j+\sum_k\lambda_{b,k}\bigl(k+f(c_j,s-2k)\bigr)\right].
\]
The outer coefficients are nonnegative, sum to one, and do not
depend on $b$. Their fixed weighted average lifts the result to
arbitrary remaining factors.
\end{proof}

\section{The bound \texorpdfstring{$n\leq81$}{n <= 81}}\label{two:sec:conclusion}

We now compare the vertex bounds with the necessary lower bound
$C_n\geq\tau_n$ from Proposition~\ref{two:prop:criterion}.
Since $\tau_n$ is increasing for $n\geq82$, a uniform vertex bound suffices for
all large dimensions. We first prove that both finite vertices
and ideal vertices, after subtracting the retained cusp term,
have normalized total strictly less than $83/8$. As
$\tau_{90}=83/8$, this excludes $n\geq90$.
For the remaining eight dimensions, we retain
$\sum_i d_i=n-1$ and compute the sharper maximum for each $n$.

\subsection{Uniform bounds for \texorpdfstring{$n\geq90$}{n >= 90}}
We first obtain a bound independent of the total dimension of the
cusp. For this purpose, replace $\psi$ by the simpler piecewise
linear upper bound $\bar\psi$, defined for real $d\geq0$ by
\begin{equation}\label{two:eq:envelope}
\bar\psi(d)=\min\left\{
\frac{5d}{8},\ \frac{19d}{40}+\frac65,\
\frac d4+\frac{24}{5},\
\frac{13d}{160}+\frac{327}{40},\ \frac{81}{8}
\right\}.
\end{equation}
Its breakpoints are
\[
(0,0),\quad(8,5),\quad(16,44/5),\quad
(20,49/5),\quad(24,81/8).
\]
The finite pair sums specified in Subsection~\ref{two:app:budgets}
verify $\psi(d)\leq\bar\psi(d)$ for $d\leq27$.
For $d\geq27$, the stable formula \eqref{two:eq:stable-factor-budget}
gives $\psi(d)<81/8$. Hence the inequality holds for every
integer $d\geq0$. The purpose of $\bar\psi$ is to make the
endpoint argument below possible on each of its linear pieces.

\begin{lem}\label{two:lem:uniform-cusp}
For every finite list of nonnegative integer dimensions,
\begin{equation}\label{two:eq:uniform-cusp}
\mathcal F(d_1,\ldots,d_r)<\frac{83}{8}.
\end{equation}
\end{lem}
\begin{proof}
Replace $\sum_i\psi(d_i)$ by $\sum_i\bar\psi(d_i)$ in
\eqref{two:eq:net-cost}. If a dimension exceeds $24$, reduce
it to $24$. The sum of upper bounds is unchanged, while $m$
cannot increase by Lemma~\ref{two:lem:coordinate-concavity}.
The resulting upper bound for $\mathcal F$ therefore does not
decrease.

On each interval where $\bar\psi$ is linear, the expression
$\sum_i\bar\psi(d_i)-3m/2$, with all other coordinates fixed,
is a convex integer sequence. Its largest value on that interval
is attained at an endpoint. Moving the coordinates successively
to these endpoints gives a list
\begin{equation}\label{two:eq:canonical-list}
(8^a,16^b,20^c,24^d),\qquad a,b,c,d\geq0,
\end{equation}
whose upper bound is at least the original one. Here $8^a$
means $a$ factors of dimension eight, and similarly for the
other entries. Zero-dimensional factors can be omitted.
This reduction does not preserve the total dimension, which is
appropriate for the present uniform estimate.

For the list \eqref{two:eq:canonical-list}, let
\[
\eta=a+\frac{16}{5}b+\frac{26}{5}c+\frac{63}{8}d
=\frac{3R}{4}-\sum_i\bar\psi(d_i).
\]
At each nonzero breakpoint, direct substitution gives
\[
d_i(d_i+2)\leq90\left(\frac{3d_i}{4}-\bar\psi(d_i)\right).
\]
By \eqref{two:eq:centered-identity} and Cauchy--Schwarz,
\[
(R+1-2m)^2
=\left(\sum_j(R+1-2j)\pi_j\right)^2
\leq\sum_j(R+1-2j)^2\pi_j\leq1+90\eta.
\]
Every $R+1-2j$ in the weighted average is positive, so
$R+1-2m>0$. Taking the positive square root gives
\[
m\geq\frac{R+1-\sqrt{1+90\eta}}2.
\]
Substituting this inequality and
$\sum_i\bar\psi(d_i)=3R/4-\eta$ cancels the terms linear in
$R$. We obtain a bound depending only on $\eta$:

\begin{equation}\label{two:eq:eta-bound}
\sum_i\bar\psi(d_i)-\frac32m
\leq-\eta-\frac34+\frac34\sqrt{1+90\eta}.
\end{equation}
For $\eta\geq24$, the right-hand side is strictly less than
$83/8$. Squaring the equivalent comparison between positive
quantities reduces it to
\[
64\eta^2-1816\eta+7885>0.
\]
Writing $\eta=24+t$ gives $64t^2+1256t+1165>0$.

For $\eta<24$, there are only $341$ quadruples in
\eqref{two:eq:canonical-list}. They satisfy
\[
0\leq a\leq23,\qquad 0\leq b\leq7,\qquad
0\leq c\leq4,\qquad 0\leq d\leq3,
\]
together with the stated inequality for $\eta$.
Multiplying the factors $1+t+\cdots+t^{d_i}$ and computing $m$ by
\eqref{common:eq:weighted-averages} gives the exact maximum
\begin{equation}\label{two:eq:341-maximum}
\max_{\eta<24}
\left(5a+\frac{44}{5}b+\frac{49}{5}c+\frac{81}{8}d
-\frac32m\right)=\frac{254166}{24565}.
\end{equation}
It occurs at $(a,b,c,d)=(0,3,0,0)$, and
\[
\frac{83}{8}-\frac{254166}{24565}=\frac{5567}{196520}>0.
\]
The finite comparison and \eqref{two:eq:eta-bound} prove
\eqref{two:eq:uniform-cusp} for every list.
\end{proof}

The finite-vertex estimate in Lemma~\ref{two:lem:finite-vertex}
and \eqref{two:eq:net-cost} now show that the normalized total
at every vertex, after subtracting the retained cusp contribution
at ideal vertices, is strictly less than $83/8$.
For a fixed polytope, the maximum of these finitely many values
is still strictly less than $83/8$. Since $\tau_{90}=83/8$ and
$\tau_n\geq83/8$ for all $n\geq90$, Proposition~\ref{two:prop:criterion}
excludes every dimension in this range.

\subsection{Bounds for \texorpdfstring{$82\leq n\leq89$}{82 <= n <= 89}}
For $82\leq n\leq89$, retain the relation $\sum_i d_i=n-1$
and use the sharper factor bound $\psi$.

\begin{lem}\label{two:lem:small-factors}
At fixed total dimension, the largest value of
$\mathcal F(d_1,\ldots,d_r)$ can be attained on a list in which
every dimension is at least five, except possibly one dimension
in $\{1,2,3,4\}$.
\end{lem}
\begin{proof}
Suppose two dimensions $a,b$ belong to $\{1,2,3,4\}$.
Keep their sum $s=a+b$ fixed and vary $b$ between
$\max(0,s-5)$ and $\min(5,s)$.
Both dimensions stay in $[0,5]$, so
$\psi(a)+\psi(b)=5s/8$ remains constant.
By Lemma~\ref{two:lem:fixed-sum-concavity}, the remaining term
$-3m/2$ is convex on this integer interval. At least one endpoint
therefore gives a value of $\mathcal F$ no smaller than the
original value. At that endpoint one of the two dimensions is
zero or five. Remove a zero-dimensional factor and repeat.
Each step reduces the number of dimensions in $\{1,2,3,4\}$,
without changing the total dimension.
\end{proof}

Define the cusp bound
\begin{equation}\label{two:eq:cusp-Cn}
C_n^{\mathrm{cusp}}=
\max_{\substack{5\leq d_1\leq\cdots\leq d_r,\ 0\leq s\leq4\\
d_1+\cdots+d_r+s=n-1}}
\mathcal F(d_1,\ldots,d_r,s).
\end{equation}
A zero final dimension is omitted. By the preceding lemma,
this finite set covers all possible cusp links; the large factors
are not required to have equal dimensions.

For finite vertices, let $\nu_r$ be an upper bound for the total
$J(T)$ of a connected elliptic diagram with $r$ nodes. For
$r\leq8$ use $\nu_r=\tfrac54\binom r2$, and for $r\geq9$
take the largest value among $A_r,B_r,D_r$.
For $9\leq r\leq27$, compute the values by the pair-sum rules in
Subsection~\ref{two:app:budgets}; for $r\geq28$,
$\nu_r=15r/2+47$. Set
\begin{equation}\label{two:eq:finite-DP}
\Lambda_0=0,\qquad
\Lambda_N=\max_{1\leq r\leq N}(\Lambda_{N-r}+\nu_r),\qquad
C_n^{\mathrm{fin}}=\Lambda_n/n.
\end{equation}
Splitting an elliptic diagram into its connected components shows
that this recurrence bounds every finite vertex of $P$.

\begin{lem}\label{two:lem:endpoint-table}
Let $C_n=\max(C_n^{\mathrm{cusp}},C_n^{\mathrm{fin}})$.
For $82\leq n\leq89$, the exact calculations
\eqref{two:eq:cusp-Cn} and \eqref{two:eq:finite-DP} give the
values in Table~\ref{two:tab:endpoints}. In particular $C_n<\tau_n$.
\end{lem}
\begin{proof}
For each list in \eqref{two:eq:cusp-Cn}, form the polynomial
\eqref{ext:eq:link-product-h} using integer coefficient convolution.
Then obtain $g_j$ by consecutive differences and compute $m$
using \eqref{common:eq:weighted-averages}. All comparisons in
\eqref{two:eq:cusp-Cn} are rational. The numbers of lists for
$n=82,\ldots,89$ are, respectively,
\[
226774,\ 253863,\ 283993,\ 317579,\
354899,\ 396446,\ 442587,\ 493897.
\]
The total is $2\,770\,038$. The recurrence
\eqref{two:eq:finite-DP} is another finite rational calculation.
The generating rules for the input values and the enumeration
procedure are specified in Subsection~\ref{two:app:budgets}.
\end{proof}

\begin{table}[!htb]
\centering\small
\setlength{\tabcolsep}{5pt}
\caption{Exact vertex bounds in the last eight dimensions.}
\label{two:tab:endpoints}
\begin{tabular}{rccc}
\toprule
$n$&$C_n^{\mathrm{cusp}}$&$C_n^{\mathrm{fin}}$&$\tau_n-C_n$\\
\midrule
82&$105177423/11358856$&$1513/164$&$327963/2839714$\\[3pt]
83&$26064463/2839714$&$3075/332$&$41/166$\\[3pt]
84&$103537961/11358856$&$781/84$&$55/168$\\[3pt]
85&$51508227/5679428$&$3169/340$&$149/340$\\[3pt]
86&$102692899/11358856$&$1607/172$&$183/344$\\[3pt]
87&$10474701/1169294$&$3259/348$&$56/87$\\[3pt]
88&$808290543/90870848$&$413/44$&$65/88$\\[3pt]
89&$451668671/51114852$&$3345/356$&$307/356$\\
\bottomrule
\end{tabular}
\end{table}

\begin{proof}[Proof of Theorem~\ref{thm:dimension-bound}]
Proposition~\ref{two:prop:facewise} establishes the facewise
bounds \eqref{two:eq:face-targets}. Lemma~\ref{two:lem:retained}
gives the retained cusp contribution, and
\eqref{two:eq:net-cost} bounds
$[J_u-(6/c_{n,2})\Theta_n^{(2)}(L_u)]/V_u$ at ideal vertices. The finite-vertex totals are bounded by
\eqref{two:eq:finite-DP}. For $n\geq90$, the uniform estimates
contradict Proposition~\ref{two:prop:criterion}. For
$82\leq n\leq89$, Table~\ref{two:tab:endpoints} gives the same
strict contradiction. Hence $n<82$, or equivalently $n\leq81$.
\end{proof}

At the endpoint $n=82$, the largest cusp bound is attained by
the list $(16,16,16,16,16,1)$. The largest finite-vertex bound
comes from $D_{20}\sqcup D_{20}\sqcup D_{21}\sqcup D_{21}$.
Consequently
\[
C_{82}=\frac{105177423}{11358856}<\frac{75}{8}=\tau_{82},
\qquad \tau_{82}-C_{82}=\frac{327963}{2839714}>0.
\]
These maxima range over local diagram configurations and therefore
bound every vertex of $P$, whether or not a maximizing configuration
is realizable in a hyperbolic polytope.

\begin{rem}\label{two:rem:81}
The present vertexwise comparison does not exclude dimension $81$. The list of
five factors of dimension $16$ gives
\[
\mathcal F(16,16,16,16,16)=\frac{53148337}{5679428}
>\frac{250}{27}=\tau_{81}.
\]
Further progress could come from additional restrictions
on the possible vertex diagrams or a sharper weight comparison.
This limitation of the present estimates does not imply existence
in dimension $81$.
\end{rem}

\clearpage
\appendix
\section{Lann\'er and quasi-Lann\'er diagrams}\label{app:lanner}
The diagrams below are grouped by the dimension $n$ of the corresponding
hyperbolic simplex. Each has $n+1$ nodes. We include all diagrams for
$n\geq2$, using the classifications in
\cite{Lanner1950,Chein1969,Kos67,Vinberg1985}.
For $n=1$ the only diagram is a dotted edge, which is not drawn.
In the triangular families, an edge labeled $2$ is omitted and an edge
labeled $\infty$ is bold. We set $1/\infty=0$.
\setcounter{figure}{0}
\renewcommand{\thefigure}{\thesection.\arabic{figure}}
\begin{figure}[!ht]\centering\small
\noindent\textbf{Dimension $2$.}\par\smallskip
\begin{minipage}[t]{.47\linewidth}\centering
\textit{Lann\'er.}\par
\end{minipage}
\par\vspace{.5mm}
\caption{Lann\'er and quasi-Lann\'er diagrams in dimensions $2$ and $3$.}
\label{figure:lanner}
\end{figure}
\clearpage
\begin{figure}[!ht]\centering\small
\noindent\textbf{Dimension $4$.}
\par\smallskip\noindent\textit{Lann\'er.}\par\smallskip
\begin{minipage}[c]{.245\linewidth}\centering%
\end{minipage}
\par\vspace{.5mm}
\caption{Lann\'er and quasi-Lann\'er diagrams in dimensions $4$--$9$.}
\label{figure:lanner-higher}
\end{figure}
\FloatBarrier

\clearpage\clearpage
\FloatBarrier
\section{Pseudocode for the extension check}\label{ext:app:local}
For a reduced diagram $\overline L$, an eligible marking $M$ of type I--III,
and a nonspecial scored path
$(p_0,\ldots,p_k)$, $p_0=x$, $p_k=y$, $3\leq k\leq7$, take
$\mathcal D=M\setminus\{x,y\}$ and the complete deletion tests
$\mathcal C$ defined in Subsection~\ref{ext:sec:local-algorithm-input}.
All original nodes and edges remain in the full diagram used to test
$\mathcal C$. A fresh node has only the edge just added. An old target
has degree at most one before the final step; its degree at that step
is unrestricted.

The two ends are tested independently from the same $\overline L$. A return value of false rules
out the required two-sided extension; true means only that neither end
was excluded. In the latter case the reduced weight is a safe lower
bound, even if the separate extensions cannot coexist.

\begin{algorithm}[!htb]
\caption{Necessary check for a two-sided ordinary extension}\label{ext:alg:arm}
\begingroup\small\linespread{1}\selectfont
\begin{algorithmic}[1]
\Function{TwoSidedNecessary}{$\overline L,x,y,\mathcal D,\mathcal C$}
\State $(p_0,\ldots,p_k)\gets$ the scored path, $s\gets9-k$
\State $l\gets\Call{OneSide}{\overline L,x,p_1,s,\overline L,\mathcal D,\mathcal C}$
\State $r\gets\Call{OneSide}{\overline L,y,p_{k-1},s,\overline L,\mathcal D,\mathcal C}$
\State \Return $l\ \text{and}\ r$
\EndFunction
\Statex
\Function{OneSide}{$H,u,v,s,\overline L,\mathcal D,\mathcal C$}
\If{$s=0$}\State \Return true\EndIf
\State $T\gets H\setminus\mathcal D$, $N\gets\{z:z\sim_Tu,\ z\ne v\}$
\If{$|N|>1$}\State \Return false\EndIf
\If{$N=\{z\}$}
\If{the label of $uz$ is not three}\State \Return false\EndIf
\State \Return \Call{OneSide}{$H,z,u,s-1,\overline L,\mathcal D,\mathcal C$}
\EndIf
\State $\mathcal V\gets\{\text{one fresh node}\}$
\If{$u\notin\overline L$}
\State Add the nodes in other components of $T$ to $\mathcal V$
\State If $s>1$, keep only old targets of degree at most one in $T$
\EndIf
\For{$z\in\mathcal V$}
\State $H'\gets H$ with a new edge $uz$ of label three
\If{$H'$ satisfies $\mathcal C$ and
\Call{OneSide}{$H',z,u,s-1,\overline L,\mathcal D,\mathcal C$}}
\State \Return true
\EndIf
\EndFor
\State \Return false
\EndFunction
\end{algorithmic}
\endgroup
\end{algorithm}

% Complete finite records and implementation remain in the local supplementary files.

\FloatBarrier
\section{Finite checks for planar angles}\label{two:app:finite-checks}

The standard Lann\'er and quasi-Lann\'er diagrams are listed in
Appendix~\ref{app:lanner}. This appendix specifies the finite tests
used in the planar-angle argument. All matrix operations use rational
Schur complements, after the indicated positive diagonal congruences;
all weight and endpoint comparisons are rational. The accompanying
programs give the numbered diagram inputs and the individual records.

\subsection{Minimal diagrams and long tails}
\label{two:app:minimal}
For the triangle calculation, use every minimal indefinite diagram of
orders six through ten, and every three-element marked subset. Delete
each mark in turn and apply \eqref{two:eq:weight-table} to the other
two. The $1265$ records each contain at least two marked pairs of
distance at most seven and have auxiliary total at least one.
Orders at most five are handled directly in
Lemma~\ref{two:lem:triangle} and require no enumeration of unbounded
Coxeter labels.

For the tail calculation, take each minimal indefinite Coxeter
diagram $L$ and enumerate the two-node subsets $X\subseteq V(L)$
for which both single-mark deletions are positive definite. For each node $u$,
compute $d_X(u)$ and, at distance at least six, the two rational
quantities
\[
(G_L^{-1})_{uu},\qquad
\max_{x\in X}(G_{L\setminus\{x\}}^{-1})_{uu},
\]
in unit-diagonal Gram normalization. This gives the $21$ attachment
positions in $16$ marked-diagram records and verifies every assertion
of Lemma~\ref{two:lem:tails}. Only diagrams of order at least six can
enter this calculation.

For the quadrilateral calculation, take two of these marked records,
allowing repetition, and identify one node or one edge. At every
identified node require the two nearest-mark distances to sum to at
least thirteen. Interchanging different records is counted once;
within one record, keep both orientations of a one-node gluing and
identify inverse maps for an edge gluing. Reject conflicting labels,
then test positive definiteness of the unmarked diagram $S$ and
positive semidefiniteness of every adjacent-pair diagram. Among the
$92$ candidates, the successive rejection counts are $3,65,3$; the
$21$ survivors have total $w_0$-weight at least one. A zero pivot is
accepted in a semidefinite test only when the remaining row is zero.

\subsection{Cyclic supports}
\label{two:app:cycle-witnesses}
Number the cycle modulo $k_1+k_2$, with $a=0$ and $b=k_1$.
For a deleted mark $z$, the eight supports listed in
Lemma~\ref{two:lem:cycle-delete} are: empty support; a label-three
or label-four edge to either $z-1$ or $z+1$; a label-three edge
to either $z-2$ or $z+2$; and label-three edges to both $z-1,z+1$.
Use the rational normalization from
Lemma~\ref{two:lem:cycle-delete}, taking the opposite-pair coupling
magnitude to be one for empty support and zero otherwise. Eliminate
all unmarked cycle nodes and test \eqref{two:eq:witness}, with
$x,y\in\{1/2,3/5\}$.

The three length pairs $(9,9),(9,10),(10,10)$ give $192$ cases.
Exactly $185$ have a negative definite restricted matrix. The
other seven have a $c$--$d$ path of length seven or eight.
The four length-eight cases have $(k_1,k_2)=(9,9)$. One new mark
attaches by a simple edge at an endpoint, and the other next to
an endpoint, on opposite sides of their respective deleted marks.
The latter attachment node is an internal degree-three node of
the $c$--$d$ path. These are the branch
conditions used in Lemma~\ref{two:lem:cycle-delete}; the
individual short paths and negative witnesses are retained in
the calculation records. Arbitrary arc lengths are covered by
that lemma's path comparison, rather than by the three finite
length checks alone.

\subsection{Tree-support certificates}\label{two:app:tree-records}
The supports, the $E_7,E_8$ external-leaf alternatives, and the
initial path lengths are specified in
Lemma~\ref{two:lem:tree-delete}. The verification has the following
form; each output includes the corresponding exact certificate.
\begin{verbatim}
for end_type in AA, BA, DA, BB, BD, DD:
  for s1,s2 in {5,...,9} and k in {8,9,10}:
    for each allowed support and external-leaf choice for c,d:
      form the rational Gram matrix and eliminate S_0
      if a required PD or PSD condition fails: record its failure
      elif an alternative path of length <=7 survives after
           deleting every variable segment: record that path
      else: find positive x,y in the stated rational list
            with (U^T R_0 U)_11<0 and det(U^T R_0 U)>0
            record x,y and this restricted matrix
\end{verbatim}
The program \texttt{tree\_audit.py} generates all $79\,113$ records:
$39\,389$ negative restrictions, $11\,167$ unchanged short paths,
and $28\,557$ failures of a required definiteness condition.
There are no unclassified records. The smallest positive
determinant of a negative restriction is $1/32400$.
The file \texttt{records.jsonl} retains every support, normalization,
outcome, and certificate; independent replay reconstructs the
matrices and checks the certificates. The short-path test removes
all variable internal segments before searching, so its certificate
also proves that the same short path survives arbitrary lengthening.

\begin{rem}\label{two:rem:explicit-witness}
One certificate can be checked directly. In type $AA$, take
$(s_1,k,s_2)=(9,10,9)$, give both new nodes empty support on the
original chain, and keep the minimum nonzero opposite-pair couplings
$ac,bd$. Eliminating the unmarked paths and taking $x=y=1/3$ gives
\[
U^TR_0U=
\begin{pmatrix}-11/45&-1/10\\-1/10&-11/45\end{pmatrix},
\qquad\det(U^TR_0U)=\frac{403}{8100}>0.
\]
The negative first diagonal entry and positive determinant exhibit
two negative directions. The remaining negative certificates have
the same form, with different rational entries.
\end{rem}

\subsection{Vertex bounds and endpoint checks}
\label{two:app:budgets}
The finite input bounds are generated from the weight rules.
For $9\leq r\leq27$, sum $w(T;a,b)+q(T;a,b)/4$ over all
unordered node pairs in each of $A_r,B_r,D_r$, then take the maximum
of these three diagram totals to obtain $\nu_r$.
For $6\leq d\leq26$, evaluate the four classical affine diagrams
of rank $d$ and take the maximum of $J(\Sigma)/(d+1)$ to obtain
$\psi(d)$. The bounds for smaller orders and the stable formulas
are given in Section~\ref{two:sec:budgets}. The exceptional affine
diagrams are tested separately as in
\eqref{two:eq:exceptional-budgets}. These rules specify every input
to the recurrence \eqref{two:eq:finite-DP} and the cusp search
\eqref{two:eq:cusp-Cn} without a table of small-order values.

For the fixed-dimension cusp search, choose $s\in\{0,1,2,3,4\}$
and enumerate nondecreasing lists of dimensions at least five
whose sum is $n-1-s$. At each recursion step the next entry ranges
from the preceding entry to the remaining sum; a zero remaining sum
completes one list. This produces each eligible list exactly once.
For each completed list $(d_1,\ldots,d_r)$, form the product
polynomial for $(d_1,\ldots,d_r,s)$ by integer coefficient convolution
(the factor for $s=0$ is one). Take consecutive coefficient differences
and compute $m$ from
\eqref{common:eq:weighted-averages}. Compare
$\sum_i\psi(d_i)+\psi(s)-3m/2$ by integer cross-products.
The same operations check the $341$ reduced lists in
\eqref{two:eq:341-maximum} and all $2\,770\,038$ lists for
$82\leq n\leq89$. The exact maxima and strict margins are
displayed in Table~\ref{two:tab:endpoints}. The unbounded ranks and
numbers of factors are covered by the stable formulas and concavity
arguments in the body of the paper.

\end{document}